\documentclass[preprint,3p]{elsarticle}
\journal{J. Comput. Phys.}
\usepackage{epstopdf}
\usepackage{pdfpages}
\usepackage{amssymb}
\usepackage{empheq}
\usepackage{cases}
\usepackage{amsthm,amsmath}
\usepackage{caption,lipsum}
\usepackage{stmaryrd}
\usepackage{graphicx,epsfig}
\usepackage{tabularx}
\usepackage{color}
\usepackage{empheq,float}
\DeclareGraphicsExtensions{-eps-converted-to.pdf}
\input psfig.sty
\newcommand{\sech}{\mathrm{sech}}

\newcommand{\eps}{\varepsilon}
\newcommand{\bR}{\mathbb{R}}
\newcommand{\bx}{\mathbf{x}}

\newcommand{\fe}{e}

\newcommand{\bC}{\mathbb{C}}
\newcommand{\bT}{\mathbb{T}}

\newcommand{\ud}{\mathrm{d}}

\newcommand{\nn}{\nonumber}
\newcommand{\dps}{\displaystyle}

\newcommand{\be}{\begin{equation}}
\newcommand{\ee}{\end{equation}}
\newcommand{\ba}{\begin{array}}
\newcommand{\ea}{\end{array}}
\newcommand{\bea}{\begin{eqnarray}}
\newcommand{\eea}{\end{eqnarray}}
\newcommand{\beas}{\begin{eqnarray*}}
\newcommand{\eeas}{\end{eqnarray*}}
\newtheorem{remark}{Remark}[section]

\numberwithin{equation}{section}

\begin{document}

\begin{frontmatter}

\title{Comparison of numerical methods for the nonlinear Klein-Gordon
equation in the nonrelativistic limit regime}

\author[NUS]{Weizhu Bao}
 \ead{matbaowz@nus.edu.sg}
\address[NUS]{Department of Mathematics, National
              University of Singapore, Singapore 119076, Singapore}
\address[Whu]{School of Mathematics and Statistics, Wuhan University, 430072 Wuhan, China}


\author[Whu]{Xiaofei Zhao\corref{5}}
\ead{matzhxf@whu.edu.cn}

\cortext[5]{Corresponding author.}

\begin{abstract}
Different efficient and accurate numerical methods
have recently been  proposed and analyzed
for the nonlinear Klein-Gordon  equation (NKGE)
with a dimensionless parameter $\varepsilon\in (0,1]$,
which is inversely proportional to the speed of light.
In the nonrelativestic limit regime, i.e. $0<\varepsilon\ll1$,
the solution of the NKGE  propagates waves with
wavelength at $O(1)$ and $O(\varepsilon^2)$ in space and time, respectively, which brings significantly numerical burdens in
designing numerical methods. We compare systematically
spatial/temporal efficiency and accuracy as well as
$\varepsilon$-resolution (or $\varepsilon$-scalability)
of different numerical methods including finite difference
time domain methods, time-splitting method, exponential wave integrator, limit integrator,
multiscale time integrator, two-scale formulation method and iterative exponential integrator.
Finally, we adopt the multiscale time integrator to study the convergence
rates from the NKGE to its limiting models when $\varepsilon\to0^+$.
\end{abstract}


\end{frontmatter}

\section{Introduction}
\label{sec:intro}
Consider the dimensionless nonlinear Klein-Gordon  equation (NKGE) in $d$-dimensions ($d=1,2,3$) with cubic nonlinearity \cite{BZ,BD,Cla,Davydov,Gru}:
\begin{equation}\label{KG}
\left\{
  \begin{split}
    & \eps^2\partial_{tt}u(\bx,t)-\Delta u(\bx,t)+\frac{1}{\eps^2}u(\bx,t)
    +\lambda\,|u(\bx,t)|^2u(\bx,t)=0,\quad \mathbf{x}\in\bR^d,\quad t>0,\\
    & u(\mathbf{x},0)= \phi_1(\mathbf{x}),\quad\partial_tu(\mathbf{x},0)=
    \frac{1}{\eps^2}\phi_2(\mathbf{x}),\quad \bx\in\bR^d,
  \end{split}
\right.
\end{equation}
where $t$ is time, $\mathbf{x}\in\bR^d$ is the spatial coordinate,
$u:=u(\mathbf{x},t)$ is a complex-valued scalar field, $0<\eps\leq1$ is a dimensionless parameter inversely proportional to the speed of light,
$\lambda\in \bR$ is a given dimensionless parameter (positive and negative for defocusing and focusing self-interaction, respectively), and
 $\phi_1$ and $\phi_2$ are given complex-valued $\eps$-independent initial data.

 When $\lambda=0$, the above Klein-Gordon equation is known as the relativistic version of the Schr\"{o}dinger equation for correctly describing the spinless relativistic composite particles, like the pion and the Higgs boson \cite{Davydov}. When $\lambda\ne0$, the NKGE  was widely
adapted in plasma physics for modeling interaction between Langmuir and ion sound waves \cite{Bellan,Dendy} and in cosmology as a phonological model for dark-matter and/or black-hole evaporation \cite{Huang,Xiong}.
The NKGE (\ref{KG}) is time symmetric and conserves the {\sl energy}
\begin{align}\label{energy}
E(t)&:=\int_{\bR^d}\left[\eps^2|\partial_t u(\bx,t)|^2+|\nabla u(\bx,t)|^2+\frac{1}{\eps^2}|u(\bx,t)|^2+\frac{\lambda}{2}|u(\bx,t)|^4\right]d \bx\nn\\
&\equiv\int_{\bR^d}\left[\frac{1}{\eps^2}|\phi_2(\bx)|^2+|\nabla \phi_1(\bx)|^2+\frac{1}{\eps^2}|\phi_1(\bx)|^2+\frac{\lambda}{2}|\phi_1(\bx)
|^4\right]d \bx=E(0),\qquad t\geq0.
\end{align}
For the derivation and nondimensionlization of (\ref{KG}), we refer to
\cite{BD,Davydov,Machihara,Masmoudi} and references therein.
For the well-posedness of the Cauchy problem (\ref{KG}) with a fixed $\eps\in(0,1]$, e.g. $\eps=1$, we refer to \cite{Davydov,Machihara,Masmoudi}
and references therein. We remark here that when the initial data
$\phi_1$ and $\phi_2$ are real-valued, then the solution $u$ of
(\ref{KG}) is also real-valued; and this case has been widely
studied analytically and numerically in the literature \cite{Adomian1,Bainov,google2,Cla,Deeba,google1,Duncan,Ginibre_KG,Ginibre_KG_2,Gru,Ibrahim,Jimenez,
Khalifa,Li-Guo,Add6,Add7,Strauss}.
For simplicity of notations and without loss of generality,
from now on, we assume that $\phi_1$ and $\phi_2$ are real-valued,
and thus the solution $u$ of
(\ref{KG}) is real-valued too.
Our methods and results can be straightforwardly extended to the case when $\phi_1$ and $\phi_2$ are complex-valued and/or general nonlinearity in (\ref{KG})
\cite{BZ,BD},
and the conclusion will be remained the same.

\begin{figure}[t!]
\centerline{\psfig{figure=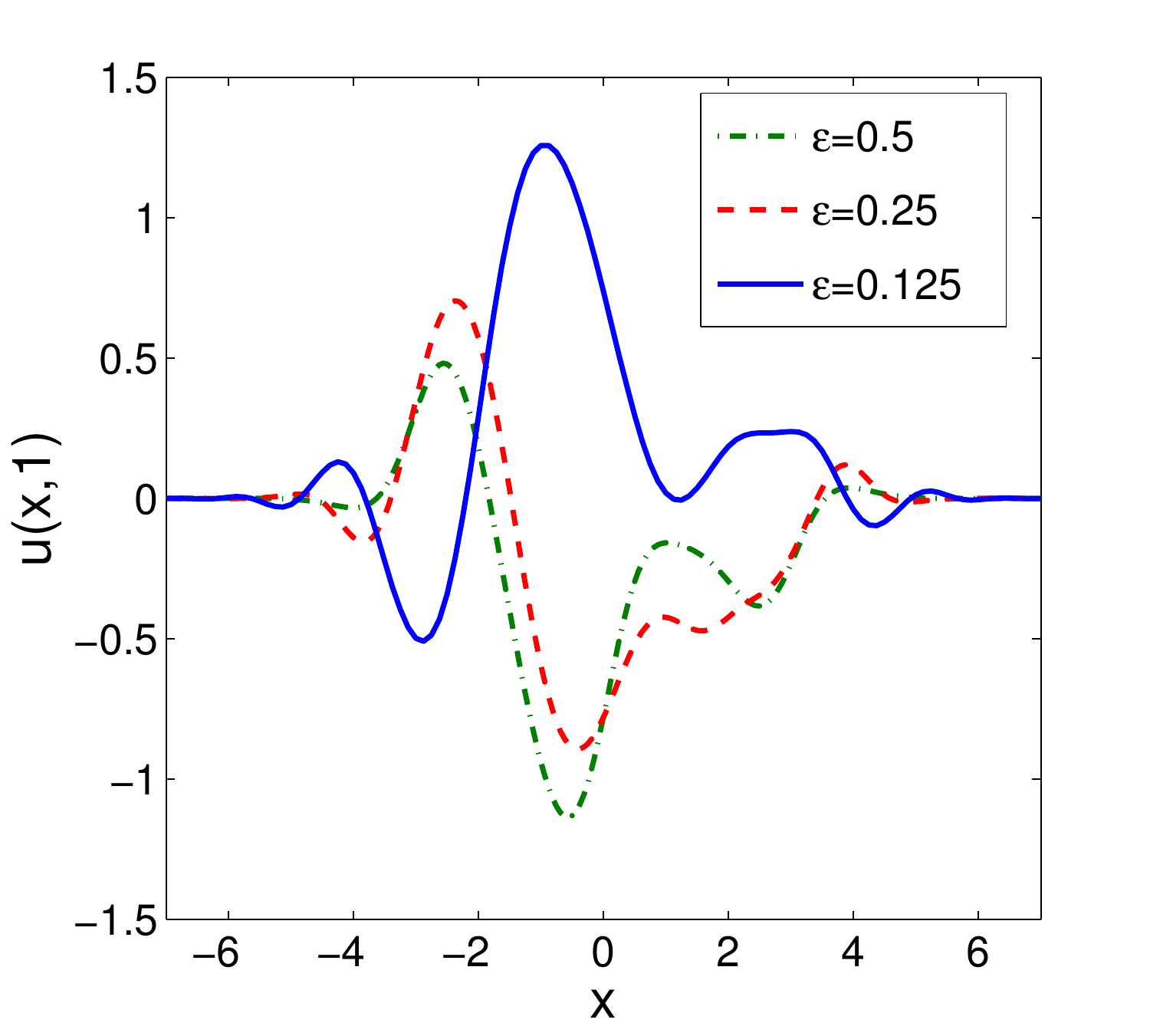,height=6cm,width=8cm}
\psfig{figure=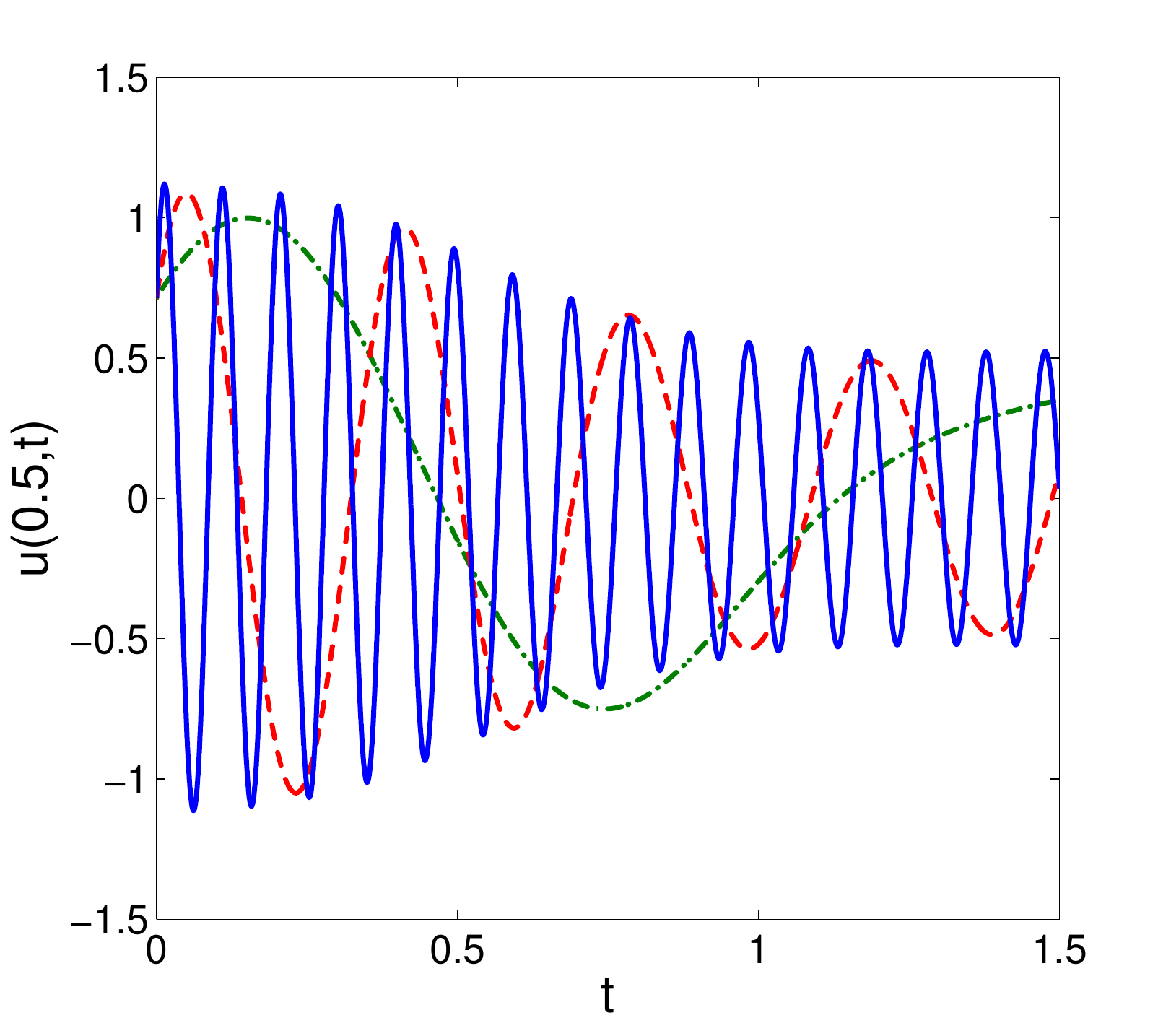,height=6cm,width=8cm}}
\caption{Solution of the NKGE (\ref{KG})  with
$d=1$, $\lambda=1$ and the initial data taken as \eqref{initKG} for different $\eps$.}\label{fig:intro}
\end{figure}

When $\eps\to0^+$ in (\ref{KG}), due to that the energy $E(t)=O(\eps^{-2})$
in (\ref{energy}) becomes unbounded, the analysis of the nonrelativistic limit of the solution $u$ becomes challenging and quite complicated.
Fortunately, convergence from the NKGE (\ref{KG}) to
a nonlinear Schr\"{o}dinger equation  has been extensively studied
in the mathematics literature \cite{BZ,Machihara,Masmoudi,Najman}.
Based on their results, the solution $u$ of the NKGE (\ref{KG})
propagates waves with wavelength at $O(\eps^2)$ and $O(1)$ in time and space, respectively, in the nonrelativistic limit regime, i.e. $0<\eps\ll1$.
To illustrate this, Figure \ref{fig:intro} shows the solution of (\ref{KG})
with $d=1$ and $\lambda=1$ and the initial data
\be\label{initKG}
\phi_1(x)=\frac{3\sin(x)}{\fe^{x^2/2}+\fe^{-x^2/2}},\qquad \phi_2(x)=\frac{2\fe^{-x^2}}{\sqrt{\pi}}, \qquad x\in{\mathbb R}.
\ee

In fact, when $0<\eps\ll1$, formally by taking the ansatz \cite{Machihara,Masmoudi}
\begin{equation}\label{ansatz intr}
u(\bx,t)=\fe^{it/\eps^2}z(\bx,t)+\fe^{-it/\eps^2}\overline{z}(\bx,t)+O(\eps^2),
\qquad \mathbf{x}\in\bR^d,\quad t\ge0,
\end{equation}
where $z:=z(\bx,t)$ is a complex-valued function and $\overline{z}$ denotes the complex conjugate of $z$, the NKGE (\ref{KG}) can be formally reduced to -- semi-limiting model -- the nonlinear Schr\"{o}dinger equation with wave operator (NLSW) under well-prepared initial data \cite{Cai1,Cai2}
\begin{equation}\label{nlsw}
\left\{
  \begin{split}
&2i\partial_tz(\bx,t)+\eps^2\partial_{tt}z(\bx,t)-\Delta z(\bx,t)+3\lambda|z(\bx,t)|^2z(\bx,t)=0,\qquad  \bx\in\bR^d, \quad t>0,\\
&z(\bx,0)=\frac{1}{2}\left[\phi_1(\bx)-i\phi_2(\bx)\right]=:z_0(\bx),\qquad \bx\in\bR^d,\\
&\partial_tz(\bx,0)=\frac{i}{2}\left[-\Delta z_0(\bx)+3\lambda|z_0(\bx)|^2z_0(\bx)\right].
\end{split}\right.
\end{equation}
In addition, by dropping the small term $\eps^2\partial_{tt}z$ in (\ref{nlsw}), one gets -- limiting model -- the nonlinear Schr\"{o}dinger  equation (NLSE) \cite{Cai1,Cai2,Machihara,Masmoudi}
\begin{equation}\label{nls}
\left\{
  \begin{split}
&2i\partial_tz(\bx,t)-\Delta z(\bx,t)+3\lambda|z(\bx,t)|^2z(\bx,t)=0,\qquad \bx\in\bR^d,\quad t>0,\\ &z(\bx,0)=\frac{1}{2}\left[\phi_1(\bx)-i\phi_2(\bx)\right]:=z_0(\bx),\qquad \bx\in\bR^d.
\end{split}\right.
\end{equation}

When $\eps=1$ in (\ref{KG}), i.e. $O(1)$-wave speed regime,
several numerical methods have been proposed and analyzed
for the Cauchy problem (\ref{KG}) in the literature, see
\cite{google2,Duncan,Iserles,Jimenez,Add6,Add7} and references therein.
Specifically, the finite difference time domain (FDTD) methods \cite{Duncan,Jimenez,Strauss} have been
demonstrated excellent performance in terms of efficiency and accuracy
for  (\ref{KG}) when $\eps=1$.
However, when $0<\eps\ll1$, i.e. in the nonrelativistic limit regime,
it becomes much more challenging in
designing and analyzing efficient and accurate numerical methods
for  (\ref{KG}) due to the highly oscillatory nature of the solution in time
(cf. Fig. \ref{fig:intro}). To address this issue,  Bao and Dong \cite{BD}
established rigorous error bounds of
the FDTD methods for (\ref{KG}), which depends
explicitly on the mesh size $h$ and time step $\tau$ as well as the small
parameter $\eps\in(0,1]$. Based on their results \cite{BD}, in order to obtain the `correct' numerical solution of (\ref{KG}), the $\eps$-resolution (or $\eps$-scalability or meshing strategy) of the FDTD methods is $\tau=O(\eps^3)$ and $h=O(1)$,
which is  {\bf under-resolution} in time with
respect to $\eps\in(0,1]$ regarding to the
Shannon's information theory \cite{Lan,Shan1,Shan2} -- to resolve a wave one needs a few points per wavelength --
since the wavelength in time is at $O(\eps^2)$. To overcome
the temporal under-resolution of the FDTD methods,
they \cite{BD} proposed to adapt the exponential wave integrator (EWI) \cite{Hochbruck}
for discretizing temporal derivatives in (\ref{KG}) and showed rigorously
the $\eps$-resolution of EWI is $\tau=O(\eps^2)$ and $h=O(1)$,
which is {\bf optimal-resolution} in time with
respect to $\eps\in (0,1]$. Later, the time-splitting (TS) method \cite{Mc} was also applied to discretize the NKGE (\ref{KG}), and the method was shown as equivalent to one type EWI and thus it retains the same $\eps$-resolution as that of the EWI but with an improved error bound regarding to the small parameter $\eps\in(0,1]$ \cite{Zhaocicp}. In fact, FDTD, EWI and TS methods perform very well when $\tau\to0$ under $\eps=\eps_0$ fixed and they lose accuracy
when $\eps\to0$ under $\tau=\tau_0$ fixed.
At the meantime, Faou and Schratz \cite{Faou} presented
a class of limit integrators (LI) for (\ref{KG}) via solving numerically
the limiting model NLSE (\ref{nls}) and obtained their error bounds.
On the contrary, the LI methods perform very well when $\eps\to0$ under $\tau=\tau_0$ fixed and they lose accuracy
when $\tau\to0$ under $\eps=\eps_0$ fixed.

It is a natural question to ask on whether one can design a numerical method
for the NKGE \eqref{KG} such that it is uniformly accurate for
$\eps\in(0,1]$, i.e. {\bf super-resolution} in time, especially in the nonrelativistic limit regime, since we have the solution structure
\eqref{ansatz intr} of the NKGE \eqref{KG} via the limiting model NLSW \eqref{nlsw} or NLSE \eqref{nls}.
Recently, different uniformly accurate (UA) numerical methods have been
designed and analyzed for the NKGE
\eqref{KG} including
a multiscale time integrator (MTI) via a multiscale decomposition of the
solution \cite{BZ} and a two-scale formulation (TSF) method \cite{Chartier}
as well as two uniformly and optimally accurate (UOA) methods \cite{BZnew,Bau}.
The main aim of this paper is to carry out a systematical comparison
of different numerical methods which have been proposed for the NKGE \eqref{KG}
in terms of temporal/spatial accuracy and efficiency as well as
$\eps$-resolution for $\eps\in(0,1]$, especially in the nonrelativistic limit regime.

The rest of the paper is organized as follows. The FDTD, EWI and TS methods as well as the LI schemes for the NKGE \eqref{KG} are briefly reviewed in Section \ref{sec:classical};  the uniformly accurate methods are
briefly reviewed in Section \ref{sec:UA}; and the uniformly and optimally accurate methods are briefly reviewed in Section \ref{set:UOA}. In Section \ref{sec:result}, we present detailed comparison of different numerical methods; and in Section \ref{appl9}, we report convergence
rates of the NKGE \eqref{KG} to its limiting models NLSW \eqref{nlsw} and
NLSE \eqref{nls} and show wave interactions of NKGE in two dimensions (2D). Finally, some conclusions are drawn in Section \ref{sec:conc}.
Throughout this paper, we adopt the notation $A\lesssim B$ to represent that there exists a generic constant $C>0$, which is independent of $\tau$ (or $n$), $h$ and $\eps$, such that $|A|\leq CB$.

\section{Non-uniformly accurate numerical methods}\label{sec:classical}
In this section, we briefly review the FDTD, EWI and TS methods \cite{BD}
as well as the LI methods \cite{Faou}
which have been proposed in the literature for discretizing
the NKGE \eqref{KG trun} (or (\ref{KG})).

For simplicity of notation and without loss of generality,
we only present the numerical methods in one dimension (1D).
Generalization to high dimensions is straightforward by tensor product.
As adapted in the literature \cite{BZ,Bau,Chartier}, the NKGE
(\ref{KG}) with $d=1$ is usually truncated
onto a bounded interval $\Omega=(a,b)$ ($|a|$ and $b$ are usually taken large enough such that the truncation error is negligible) with periodic boundary condition
\begin{equation}\label{KG trun}
\left\{
  \begin{split}
&\eps^2\partial_{tt}u(x,t)-\partial_{xx}u(x,t)+\frac{1}
{\eps^2}u(x,t)+\lambda( u(x,t))^3=0,\qquad x\in\Omega,\quad t>0,\\
&u(x,0)=\phi_1(x),\qquad\partial_{t}u(x,0)=\frac{1}{\eps^2}\phi_2(x),\qquad x\in\overline{\Omega},\\
&u(a,t)=u(b,t),\qquad \partial_xu(a,t)=\partial_xu(b,t).
\end{split}\right.
\end{equation}
Choose $\tau>0$ be the time step and $h=(b-a)/N$ be the mesh size with $N$ an even positive integer, denote the grid points as $x_j=a+jh$ for $j=0,1,\ldots,N$ and time steps as $t_n=n\tau$ for $n\ge0$.
Let $u_j^n$ be the numerical approximation of $u(x_j,t_n)$ for
$0\le j\le N$ and $n\ge0$ and $u^n=(u_0^n,u_1^n,\ldots,u_N^n)^T$ be the solution vector at $t=t_n$, and define
\be
\|u^n\|_{l^2}^2=h\sum_{j=0}^{N-1} |u_j^n|^2, \qquad n\ge0.
\ee

\subsection{Finite difference time domain (FDTD) methods}
 Introduce the finite difference operators as
$$\delta_t^2u_j^n=\frac{u_j^{n+1}-2u_j^n+u_j^{n-1}}{\tau^2},\quad
\delta_x^+u_j^n=\frac{u_{j+1}^n-u_{j}^n}{h},\quad
\delta_x^2u_j^n=\frac{u_{j+1}^n-2u_j^n+u_{j-1}^n}{h^2}.
$$
As used in \cite{BD}, the \emph{Crank-Nicolson finite difference (CNFD) method} for
discretizing \eqref{KG trun} reads
\begin{align}\label{cnfd}
\eps^2\delta_t^2u_j^n+\left[-\frac{1}{2}\delta_x^2
+\frac{1}{2\eps^2}+
\frac{\lambda}{4}[(u_j^{n+1})^2+(u_j^{n-1})^2]\right](u_j^{n+1}+u_j^{n-1})=0,
\quad 0\le j\le N-1,\quad n\ge0.
\end{align}
Similarly, the \emph{semi-implicit finite difference (SIFD) method} is \cite{BD}
\begin{align}\label{sifd}
\eps^2\delta_t^2u_j^n-\frac{1}{2}\delta_x^2(u_j^{n+1}+u_j^{n-1})
+\frac{1}{2\eps^2}(u_j^{n+1}+u_j^{n-1})+\lambda (u_j^{n})^3=0,\quad 0\le j\le N-1,\quad  n\ge0;
\end{align}
and the \emph{leap-frog finite difference (LFFD) method} is \cite{BD}
\begin{align}\label{lffd}
\eps^2\delta_t^2u_j^n-\delta_x^2u_j^{n}
+\frac{1}{\eps^2}u_j^{n}+\lambda(u_j^{n})^3=0,\quad 0\le j\le N-1,\quad n\ge0.
\end{align}
The initial and boundary conditions in \eqref{KG trun}  are discretized as \cite{BD,BDZJMS}
\begin{align}\label{init-dis}
\begin{split}
&u_0^n=u_N^n,\quad u_{-1}^n=u_{N-1}^n,\quad n\geq0;\qquad u_j^0=\phi_1(x_j),\quad j=0,\ldots,N,\\
&u_j^1=\phi_1(x_j)+\sin\left(\frac{\tau}{\eps^2}\right)
\phi_2(x_j)+\frac{\tau}{2}\sin\left(\frac{\tau}{\eps^2}\right)
\left[\delta_x^2\phi_1(x_j)-\frac{1}{\tau}\sin\left(\frac{\tau}
{\eps^2}\right)\phi_1(x_j)-\lambda(\phi_1(x_j))^3\right].
\end{split}
\end{align}
We remark here that we adapt \eqref{init-dis} to compute the approximation
at $t=t_1$ instead of the classical method
\be
u_j^1=\phi_1(x_j)+\frac{\tau}{\eps^2}\phi_2(x_j)+\frac{\tau^2}{2\eps^2}
\left[\delta_x^2\phi_1(x_j)-\frac{1}{\eps^2}\phi_1(x_j)-
\lambda(\phi_1(x_j))^3\right],\qquad j=1,2,\ldots, N-1,
\ee
i.e. replacing $\tau/\eps^2$ by $\sin(\tau/\eps^2)$, such that
the numerical solution $u^1$ is uniformly bounded for $\eps\in(0,1]$.

As observed and stated in \cite{BD}, the above CNFD, SIFD and LFFD methods
are time symmetric and their memory cost is $O(N)$. The LFFD method is explicit and its computational cost per step is $O(N)$. It is conditionally stable and  there is a severe
stability condition which depends on both $h$ and $\eps$, especially when
$0<\eps\ll1$ \cite{BD}. In fact, it is the most efficient and accurate method among all FDTD methods for the NKGE when $\eps=1$.
The SIFD method is implicit, but it can be solved efficiently via the fast Fourier transform (FFT) and thus its computational cost per step is $O(N\,\ln N)$. It is conditionally stable and  the stability condition depends on
$\eps$ and is independent of $h$ \cite{BD}.
The CNFD method is implicit and at every time step a fully nonlinear coupled
system needs to be solved. One main advantage is that it
conserves the energy (\ref{energy}) in the discrete level as \cite{BD}
\begin{align}
E^n&:=\eps^2\|\delta_t^+u^n\|_{l^2}^2+\frac{1}{2}\left(\|\delta_x^+u^n\|_{l^2}^2+
\|\delta_x^+u^{n+1}\|_{l^2}^2\right)+\frac{1}{2\eps^2}\left(\|u^n\|_{l^2}^2
+\|u^{n+1}\|_{l^2}^2\right)
+\frac{h\lambda}{4}\sum_{j=0}^{N-1}\left[(u_j^n)^4+(u_j^{n+1})^4\right]\nonumber\\
&\equiv E^0,\qquad n\ge0,
\end{align}
which immediately implies that it is unconditionally stable when $\lambda\ge0$. In addition, under proper regularity of the solution
$u$ of the NKGE \eqref{KG trun} and stability conditions for
the SIFD and LFFD methods, the following rigorous error bound
was established for the three FDTD methods \cite{BD}
\be
\|e^n\|_{l^2}+\|\delta_x^+e^n\|_{l^2}\lesssim h^2+\frac{\tau^2}{\eps^6},\qquad 0\le n\le \frac{T}{\tau},
\ee
where $T>0$ is a fixed time and  the error function $e^n$ is defined as $e_j^n = u(x_j,t_n)-u_j^n$ for
$0\le j\le N$ and $n\ge0$.
This error bound suggests that the FDTD methods are second order in both space and time discretization for any fixed $\eps=\eps_0$ and the $\eps$-resolution of the FDTD methods
is $h=O(1)$ and $\tau=O(\eps^3)$ in the nonrelativistic limit regime, i.e.
$0<\eps\ll1$, which immediately show that the temporal resolution with respect to $\eps\in(0,1]$ of the FDTD methods is {\bf under-resolution} in time
since the wavelength in time is at $O(\eps^2)$.

\subsection{Exponential wave integrator (EWI)}
As it has been proposed in \cite{BD}, the NKGE \eqref{KG trun}
is discretized in space by the Fourier (pseudo)spectral method
and followed by adapting an exponential wave integrator (EWI)
in time which has been widely used
for discretizing second order oscillatory differential equations
in the literature
\cite{Cohen,Deuflhard,Gautschi,Add3,Add4,Lubich,Add5,Hochbruck}.

Let $u_j^n$ and $\dot{u}_j^n$ be the approximations of $u(x_j,t_n)$ and $\partial_tu(x_j,t_n)$, respectively,  for
$0\le j\le N$ and $n\ge0$ and take $u_j^0=\phi_1(x_j),\dot{u}_j^0=\phi_2(x_j)/\eps^2$ for $j=0,1,\ldots,N$.
When the EWI is taken as the
Gautschi's quadrature \cite{BD,Gautschi,Add3,Add4,Add5},
a Gautschi-type exponential wave integrator Fourier pseudospectral (EWI-FP) method
\cite{BD} reads as:
\begin{equation}\label{GauFP}
  u_j^{n+1}=\sum_{l=-N/2}^{N/2-1}\widetilde{(u^{n+1})}_l\;
\fe^{i\mu_l(x_j-a)}=\sum_{l=-N/2}^{N/2-1}\widetilde{(u^{n+1})}_l\;
\fe^{2ijl\pi/N},\qquad j=0,1,\ldots,N, \quad n\ge0,
\end{equation}
where
\[ 
\widetilde{(u^{n+1})}_l =\left\{\begin{array}{ll}
\left[\cos(\omega_l^0\tau)+\displaystyle\frac{\alpha^0\left(1-\cos(\omega_l^0\tau)
  \right)}{(\eps\omega_l^0)^2}\right] \widetilde{(u^0)}_l
  +\dps\frac{\sin(\omega_l^0\tau)}{\omega_l^0} \widetilde{(\dot{u}^0)}_l
 +\dps \frac{\cos(\omega_l^0\tau)-1}{(\eps\omega_l^0)^2}\widetilde{(f^0)}_l, &n=0,\\
 \\
 -\widetilde{(u^{n-1})}_l+2\left[\cos(\omega_l^n\tau)
 +\dps\frac{\alpha^n\left(1-\cos(\omega_l^n\tau)\right)}{(\eps\omega_l^n)^2}\right]
 \widetilde{(u^n)}_l+ \dps \frac{2\left(\cos(\omega_l^n\tau)-1\right)}{(\eps\omega_l^n)^2}
 \widetilde{(f^n)}_l, &n\ge1,\\
 \ea\right.
\]
with $f(u)=\lambda u^3$, $f^n:=f(u^n)=(f(u_0^n),f(u_1^n),\ldots,f(u_N^n))^T$,
$\mu_l=2l\pi/(b-a)$, $\omega_l^n= \frac{1}{\eps^2}\sqrt{1+\eps^2(\mu_l^2+\alpha^n)}$ ($l=-N/2,\ldots, N/2-1$) with
$\alpha^n = \max\left\{\alpha^{n-1},\;\max_{0\leq j\leq N}\{\lambda(u_j^n)^2\}\right\}$ for $n\ge0$ and $\alpha^{-1}=0$ being
the stabilization constants, and $\widetilde{v}_l$ ($-N/2\le l\le N/2-1$)
being the discrete Fourier transform coefficients of the vector
$v=(v_0,v_1,\ldots,v_N)^T$ with $v_0=v_N$ defined as
\be
\widetilde{v}_l=
\frac{1}{N}\sum_{j=0}^{N-1}v_j\;\fe^{-i\mu_l(x_j-a)}=
\frac{1}{N}\sum_{j=0}^{N-1}v_j\;\fe^{-2ijl\pi/N}, \qquad
l=-\frac{N}{2},-\frac{N}{2}+1,\ldots,\frac{N}{2}-1.
\ee
Of course, in practice if the approximation of the first order derivative in time is needed, then they can be obtained as \cite{BD}
\begin{equation}\label{GauFP1}
\dot{u}_j^{n+1}=\sum_{l=-N/2}^{N/2-1}\widetilde{(\dot{u}^{n+1})}_l\
\fe^{i\mu_l(x_j-a)}=\sum_{l=-N/2}^{N/2-1}\widetilde{(\dot{u}^{n+1})}_l\;
\fe^{2ijl\pi/N},\qquad j=0,1,\ldots,N, \quad n\ge0,
\end{equation}
where
\[
\widetilde{(\dot{u}^{n+1})}_l =\left\{\ba{ll}
-\omega_l\sin(\omega_l\tau)\widetilde{(u^0)}_l+\cos(\omega_l\tau)
  \widetilde{(\dot{u}^0)}_l-\dps\frac{\sin(\omega_l\tau)}{\eps^2\omega_l}
  \widetilde{(f^0)}_l, &n=0, \\
  \\
\widetilde{(\dot{u}^{n-1})}_l-2\omega_l
\sin(\omega_l\tau)\widetilde{(u^{n})}_l-2\dps\frac{\sin(\omega_l\tau)}
{\eps^2\omega_l}\widetilde{(f^n)}_l, &n\geq1.\\
\ea\right.
\]

As it can be seen, EWI-FP is explicit and
time symmetric. The memory cost is $O(N)$ and
computational cost per step is $O(N\,\ln N)$.
In addition, the EWI-FP is unconditionally stable due to
the stabilization constant $\alpha^n$ \cite{BD}.
Under proper regularity of the solution
$u$ of the NKGE \eqref{KG trun} and the assumption $\tau\lesssim \eps^2$,
the following rigorous error bound
was established for the EWI-FP method \cite{BD}
\be\label{bdEWI}
\|u(\cdot,t_n)-I_Nu^n\|_{L^2}\lesssim h^{m_0}+\frac{\tau^2}{\eps^4}, \quad \|\partial_x [u(\cdot,t_n)-I_Nu^n]\|_{L^2}\lesssim h^{m_0-1}+\frac{\tau^2}{\eps^4},\quad 0\le n\le \frac{T}{\tau},
\ee
where $m_0\ge2$ depends on the regularity of the solution $u$ of (\ref{KG trun}) and $I_N$ is the standard interpolation operator \cite{ST}. The error bounds suggest that EWI-FP is spectral order in space if the solution is smooth and is second order in time for any fixed $\eps=\eps_0$ and the $\eps$-resolution
is $h=O(1)$ and $\tau=O(\eps^2)$ in the nonrelativistic limit regime, i.e.
$0<\eps\ll1$, which immediately show that EWI-FP is {\bf optimal-resolution} in time with respect to $\eps\in(0,1]$
since the wavelength in time is at $O(\eps^2)$.  Recently, the EWI-FP method
has been extended to arbitrary even order in time \cite{Iserles,WangZhao}.

\subsection{Time-splitting (TS) method }
The time-splitting method \cite{Mc} has been widely used to solve different (partial) differential equations and it has shown great advantages in many cases, such as for the (nonlinear) Schr\"{o}dinger equation \cite{BaoC,Mc}.
As proposed in \cite{Zhaocicp}, in order to adapt the TS method for solving the NKGE, the NKGE \eqref{KG trun} is first re-formulated into a first order system
by introducing $v:=v(x,t)=\partial_tu(x,t)$. Then the  first order system is split into
\begin{equation*}
\left\{\begin{split}
         &\partial_tu=0,\\
         &\partial_tv+\frac{\lambda}{\eps^2}u^3=0,
       \end{split}\right.\quad\mbox{and}\quad
\left\{\begin{split}
 &\partial_tu-v=0,\\
&\partial_tv-\frac{1}{\eps^2}\partial_{xx}u
+\frac{1}{\eps^4}u=0.
       \end{split}\right.
\end{equation*}

Let $u_j^n$ and $v_j^n$ be the approximations of $u(x_j,t_n)$  and $v(x_j,t_n)$, respectively, for
$0\le j\le N$ and $n\ge0$ and take $u_j^0=\phi_1(x_j)$, $v_j^0=\phi_2(x_j)/\eps^2$ for $j=0,1,\ldots,N$.
Then a second-order time splitting Fourier pseudospectral
method (TS-FP) \cite{Zhaocicp} reads as:
\begin{equation}\label{DeuFP}
\begin{split}
&v_j^{(1)}=v_j^n-\frac{\lambda\tau}{2\eps^2}(u_j^n)^3,\\
&v_j^{(2)}=\sum_{l=-N/2}^{N/2-1}\left[-\omega_l\sin\left(\omega_l\tau\right)
\widetilde{(u^n)}_l+\cos\left(\omega_l\tau\right)
\widetilde{(v^{(1)}) }_l\right] \fe^{i\mu_l(x_j-a)},\\
&u_j^{n+1}=\sum_{l=-N/2}^{N/2-1}\left[\cos\left(\omega_l\tau\right)
\widetilde{(u^n)}_l    +\frac{\sin\left(\omega_l\tau\right)}{\omega_l}
    \widetilde{(v^{(1)}) }_l\right] \fe^{i\mu_l(x_j-a)},\\
&v_j^{n+1}=v_j^{(2)}-\frac{\lambda\tau}{2\eps^2}(u_j^{n+1})^3,
\end{split}
\qquad j=0,1,\ldots,N, \quad n\ge0,
\end{equation}
where $\omega_l= \frac{1}{\eps^2}\sqrt{1+\eps^2\mu_l^2}$ for
$l=-N/2,\ldots,N/2-1$.

Again, the TS-FP \eqref{DeuFP} is explicit and time symmetric.
Its memory cost is $O(N)$ and
computational cost per time step is $O(N\,\ln N)$.
We remark here that the TS-FP \eqref{DeuFP}
is mathematically equivalent to an EWI via trapezoidal quadrature (or known as Deuflhard-type exponential integrator \cite{Deuflhard}) for solving the NKGE \eqref{KG trun} (or \eqref{KG}) \cite{Zhaocicp,Lubich}.
Under the condition $\tau\lesssim\eps^2$, the following error bound was observed for the TS-FP in
\cite{Zhaocicp}:
\be\label{bdtsfp}
\|u(\cdot,t_n)-I_Nu^n\|_{L^2}\lesssim h^{m_0}+\frac{\tau^2}{\eps^2}, \quad \|\partial_x [u(\cdot,t_n)-I_Nu^n]\|_{L^2}\lesssim h^{m_0-1}+\frac{\tau^2}{\eps^2},\quad 0\le n\le \frac{T}{\tau}.
\ee
The above error bound could be rigorously obtained by the super-convergence analysis in \cite{BaoC5,BaoC6, super-TS}.
It can be seen that the error bound  \eqref{bdtsfp} is an improved
error bound compared to the error bound \eqref{bdEWI}
regarding to the small parameter $\eps$ when $0<\tau\lesssim \eps^2$ and $0<\eps\ll1$ (cf. Tabs. \ref{tab:GIFP}\&\ref{tab:DIFP}). Of course, due to the convergence restriction $\tau\lesssim \eps^2$,  the $\eps$-resolution of TS-FP is still $h=O(1)$ and $\tau=O(\eps^2)$ in the nonrelativistic limit regime. Thus the TS-FP
is also {\bf optimal-resolution} in time with respect to $\eps\in(0,1]$.

\subsection{Limit integrators (LIs)}
As presented in \cite{Faou}, a class of limit integrators (LIs)
has been designed for the NKGE \eqref{KG trun} (or \eqref{KG}) with different order of accuracy in terms of $\eps$ when $0<\eps\ll1$.
In the LIs, the limiting equation of the NKGE \eqref{KG},
e.g. \eqref{nls}, is solved numerically and the numerical solution
of the NKGE \eqref{KG} is constructed via the ansatz
\eqref{ansatz intr}.

 In practice, the NLSE \eqref{nls} in 1D is truncated on a bounded computational domain
$\Omega=(a,b)$ with periodic boundary condition and then it is discretized by the second order time-splitting Fourier pseudospectral (TSFP) method \cite{Anto,BaoC,BaoJ,BaoJ1}.  Let $u_j^n$ and $z_j^n$ be the approximations of $u(x_j,t_n)$  and $z(x_j,t_n)$, respectively, for
$0\le j\le N$ and $n\ge0$ and take $u_j^0=\phi_1(x_j)$, $z_j^0=z_0(x_j)$ for $j=0,1,\ldots,N$. Then
a first order  (with respect to the small parameter $\eps$) limit integrator Fourier pseudospectral (LI-FP1) method
was proposed in \cite{Faou} as:
\be\label{LI-FP1}
u^{n+1}_j=\fe^{it_{n+1}/\eps^2}\,z_{j}^{n+1}+
\fe^{-it_{n+1}/\eps^2}\;\overline{z_{j}^{n+1}}, \qquad j=0,1,\ldots, N,\quad
n\ge0,
\ee
where $z^{n+1}$ is a numerical approximation of \eqref{nls} by a TSFP method \cite{Anto,BaoC,BaoJ,BaoJ1} and is given as
\be\label{LI-FP11}
\begin{split}
&z_{j}^{(1)}=\sum_{l=-N/2}^{N/2-1}\fe^{i\mu_l^2\tau/4}\;\widetilde{(z^n)}_l
\;\fe^{i\mu_l(x_j-a)},\\
&z_j^{(2)}=\fe^{3i\lambda\tau|z_j^{(1)}|^2/2}\;z_j^{(1)},\\
&z_{j}^{n+1}=\sum_{l=-N/2}^{N/2-1}\fe^{i\mu_l^2\tau/4}\; \widetilde{(z^{(2)})}_l\;\fe^{i\mu_l(x_j-a)},
\end{split}
\qquad j=0,1,\ldots, N,\quad
n\ge0.
\ee

Again, as presented in \cite{Faou}, when $0<\eps\ll1$, formally by taking
the following ansatz (found by the modulated Fourier expansion \cite{Cohen,Cohen1,ICM,Lubich}) which is more accurate than \eqref{ansatz intr}
for approximating the solution of NKGE \eqref{KG trun},
\be
u(x,t)=\fe^{it/\eps^2}z(x,t)
+\fe^{-it/\eps^2}\overline{z}(x,t)
+\eps^2w(x,t)+O(\eps^4),
\qquad x\in\Omega,\ t\ge0,
\ee
one can obtain $z:=z(x,t)$ still satisfies the NLSE \eqref{nls} and $w:=w(x,t)$ is given by \cite{Faou}
\begin{align}\label{wxt}
w(x,t)=&-\frac{3\lambda}{4}|z(x,t)|^2\left[z(x,t)
\fe^{it/\eps^2}+\overline{z}(x,t)\fe^{-it/\eps^2}\right]+\frac{\lambda}{8}
\left[z(x,t)^3\fe^{3it/\eps^2}+\overline{z}(x,t)^3\fe^{-3it/\eps^2}\right]\nonumber\\
&+\frac{1}{2}\left[v(x,t)\fe^{it/\eps^2}
+\overline{v}(x,t)\fe^{-it/\eps^2}\right],\qquad x\in\Omega,\   t>0,
\end{align}
with $v:=v(x,t)$ satisfying \cite{Faou}
\begin{align}\label{xi def}
i\partial_t v-\frac{1}{2}\partial_{xx} v+3\lambda|z|^2v+
\frac{3\lambda}{2}z^2\overline{v}=\frac{1}{4}\partial_{xxxx} z
+\frac{51\lambda^2}{8}|z|^4z-\frac{3\lambda}{2}\partial_{xx}  (|z|^2z),\qquad x\in\Omega,\  t>0;
\end{align}
and the initial condition \cite{Faou}
\be\label{xi def-in}
v(x,0)=-\frac{\lambda}{2}z(x,0)^3
+\frac{\lambda}{4}\overline{z}(x,0)^3+\frac{3\lambda}{2}|z(x,0)|^2
\overline{z}(x,0)
+\frac{1}{2}\partial_{xx} (z(x,0)-\overline{z}(x,0))=:v_0(x),\quad x\in\overline{\Omega}.
\ee
The NLSE \eqref{nls} can be solved by the second order TSFP method \cite{Anto,BaoC,BaoJ,BaoJ1}
(cf. \eqref{LI-FP11} for the  case of 1D) as before.
In order to solve (\ref{xi def}) numerically \cite{Faou},  it is split into a kinetic part
$$\Phi_k(t):\qquad i\partial_tv=\frac{1}{2}\partial_{xx}  v,\qquad  x\in\Omega,\quad t>0,$$
and a potential part
\begin{equation}\label{pontential}
\Phi_p(t):\qquad i\partial_tv+3\lambda|z|^2v+\frac{3\lambda}{2}z^2\overline{v}=
\frac{1}{4}\partial_{xxxx}z
+\frac{51\lambda^2}{8}|z|^4z-\frac{3\lambda}{2}\partial_{xx} (|z|^2z),\qquad  x\in\Omega, \quad t>0,
\end{equation}
and then the flow is composed by a second order splitting scheme as
$\Phi(\tau)\approx\Phi_k\left(\frac{\tau}{2}\right)\Phi_p\left(\tau\right)\Phi_k\left(\frac{\tau}{2}\right).$
The kinetic part can be integrated as usual, while the potential part is integrated in its vector form by an exponential trapezoidal rule in \cite{Faou}.  

Here we present the method in 1D and truncate
(\ref{pontential}) on the bounded domain $\Omega= (a,b)$ with periodic boundary condition.
Let $u_j^n$, $z_j^n$, $w_j^n$ and $v_j^n$ be the approximations of $u(x_j,t_n)$, $z(x_j,t_n)$, $w(x_j,t_n)$  and $v(x_j,t_n)$, respectively, for
$0\le j\le N$ and $n\ge0$ and take $u_j^0=\phi_1(x_j)$, $z_j^0=z_0(x_j)$, $v_j^0=v_0(x_j)\approx -\frac{\lambda}{2}(z_j^0)^3+\frac{\lambda}{4}(\overline{z_j^0})^3+
\frac{3\lambda}{2}|z_j^0|^2 \overline{z_j^0}
+\frac{1}{2}\left(\partial_{xx}^{\mathcal{F}}z^0_j-\partial_{xx}^\mathcal{F}
\overline{z^0_j}\right)$ for $j=0,1,\ldots,N$, where $\partial_{xx}^{\mathcal F}$  is the standard Fourier pseudospectral approximation of the operator $\partial_{xx}$ on the bounded domain
$\Omega=(a,b)$, e.g. $\partial_{xx}^{\mathcal F}z_j^0=\partial_{xx}(I_Nz^0)(x_j)$ \cite{ST}.
Then a second order (with respect to the small parameter $\eps\in(0,1]$) limit integrator Fourier pseudospectral (LI-FP2) method is given \cite{Faou} as
\be
u_j^{n+1}=\fe^{it_{n+1}/\eps^2}z_j^{n+1}+\fe^{-it_{n+1}/\eps^2}
\overline{z_j^{n+1}} +\eps^2w_j^{n+1},\qquad j=0,1,\ldots, N, \quad n\ge0,\label{LI2}
\ee
where $z_j^{n+1}$ ($j=0,1,\ldots,N$) are given in \eqref{LI-FP11} and
$w^{n+1}$ is an approximation of \eqref{wxt} as
\begin{align*}
w_j^{n+1}=&-\frac{3\lambda}{4}|z_j^{n+1}|^2\left[z_j^{n+1}
\fe^{it_{n+1}/\eps^2}+\overline{z_j^{n+1}}\fe^{-it_{n+1}
/\eps^2}\right]+\frac{\lambda}{8}
\left[(z_j^{n+1})^3\fe^{3it_{n+1}/\eps^2}+
(\overline{z_j^{n+1}})^3\fe^{-3it_{n+1}/\eps^2}\right]\\
&+\frac{1}{2}\left[v_j^{n+1}\fe^{it_{n+1}/\eps^2}
+\overline{v_j^{n+1}}\fe^{-it_{n+1}/\eps^2}\right],\qquad j=0,1,\ldots, N, \quad n\ge0.
\end{align*}
Here $v^{n+1}$ is a numerical solution of \eqref{xi def} by a TSFP
method \cite{Faou} and is given as:
\[
\begin{split}
&v_{j}^{(1)}=\sum_{l=-N/2}^{N/2-1}\fe^{i\mu_l^2\tau/4}\;\widetilde{(v^n)}_l \;\fe^{i\mu_l(x_j-a)},\\
&v_j^{(2)}=\alpha_j^{(2)}+i\,\beta_j^{(2)},\\
&v_{j}^{n+1}=\sum_{l=-N/2}^{N/2-1}\fe^{i\mu_l^2\tau/4}\; \widetilde{(v^{(2)})}_l\;\fe^{i\mu_l(x_j-a)},
\end{split}
\qquad j=0,1,\ldots, N,\quad
n\ge0,
\]
where
\be\label{alpbt}
\left[\begin{array}{cc}\alpha_j^{(2)}\\ \beta_j^{(2)}\end{array}\right]=\fe^{\frac{\tau}{2}\left(
A(z_{R,j}^{n+1},z_{I,j}^{n+1})+A(z_{R,j}^{n},z_{I,j}^{n})\right)}
\left(\left[\begin{array}{cc}\alpha_j^{(1)}\\ \beta_j^{(1)}\end{array}\right]+\frac{\tau}{2}\left[\begin{array}{cc}Im(\chi_j^n)\\ -Re(\chi_j^n)\end{array}\right]\right)+\frac{\tau}{2}
\left[\begin{array}{cc}Im(\chi_j^{n+1})\\ -Re(\chi_j^{n+1})\end{array}\right],
\ee
with
\begin{align*}
&\alpha_j^{(1)}=Re\left(v_j^{(1)}\right),\quad \beta_j^{(1)}=Im\left(v_j^{(1)}\right),\quad
z_{R,j}^n=Re(z_j^n),\quad z_{I,j}^n=Im(z_j^n),\\
&A(z_R,z_I)=-\frac{3\lambda}{2}\left[\begin{array}{cc}2z_Rz_I & z_R^2+3z_I^2 \\ -3z_R^2-z_I^2 & -2z_Rz_I\end{array}\right],\\
&\chi_j^n=\frac{1}{4}\partial_{xxxx}^{\mathcal F}z^n_j
+\frac{51\lambda^2}{8}|z^n_j|^4z^n_j-\frac{3\lambda}{2}\partial_{xx}^{\mathcal F}
(|z_j^n|^2z_j^n).
\end{align*}
Here $Re (f)$ and $Im(f)$ denote the real and imaginary parts of $f$, respectively, and $\partial_{xxxx}^{\mathcal F}$ is the Fourier pseudospectral
approximation of $\partial_{xxxx}$ on the bounded domain $\Omega=(a,b)$ \cite{ST}.

As stated and proved in \cite{Faou},
both LI-FP1 and LI-FP2 are explicit, unconditionally stable and
time symmetric, and their memory cost is $O(N)$ and
computational cost per step is $O(N\,\ln N)$.
In addition, under proper regularity of the solution
$u$ of the NKGE \eqref{KG trun} and $z$ of the NLSE \eqref{nls},
the following rigorous error bound
was established for the LI-FP1 method \cite{Faou}
\be
\|u(\cdot,t_n)-I_Nu^n\|_{H^1}\lesssim h^{m_1}+\tau^2+\eps^2,\qquad n=0,1,\ldots,\frac{T}{\tau},
\ee
where $m_1\ge1$ depends on the regularity of the solution $u$ of
(\ref{KG trun}).
Similarly, the following rigorous error bound
was established for the LI-FP2 method \cite{Faou}
\be
\|u(\cdot,t_n)-I_Nu^n\|_{H^1}\lesssim h^{m_1}+\tau^2+\eps^4,\qquad n=0,1,\ldots,\frac{T}{\tau}.
\ee

These error bounds suggest that both LI-FP1 and LI-FP2 methods are spectral order in space if the solution is smooth and when $\eps\to0^+$, and
LI-FP1 and LI-FP2 are second order in time when $0<\eps\lesssim \tau$ and
$0<\eps\lesssim \tau^{1/2}$, respectively.
The $\eps$-resolution of the two methods
is $h=O(1)$ and $\tau=O(1)$ in the nonrelativistic limit regime, i.e.
$0<\eps\ll1$, which immediately show that both LI-FP1 and LI-FP2
are {\bf super-resolution} in time  with respect to $0<\eps\ll1$
since the wavelength in time is at $O(\eps^2)$.
On the contrary, when $\eps=\eps_0$ is fixed, e.g. $\eps=1$,
there is no convergence of LI-FP1 and LI-FP2 for the NKGE
\eqref{KG trun} (or \eqref{KG}).

\section{Uniformly accurate (UA) methods} \label{sec:UA}
In this section, we review the uniformly accurate
MTI \cite{BDZJMS,BZ} and  TSF method \cite{Chartier}
which have been proposed in the literature for discretizing
the NKGE \eqref{KG trun} (or (\ref{KG})).

\subsection{A multiscale time integrator (MTI)}
As proposed in \cite{BZ}, the MTI was designed via a multiscale decomposition of the solution of the NKGE (\ref{KG}) and
adapting the EWI-FP method for discretizing the decomposed sub-problems.

For any fixed $n\ge0$, by assuming that
the initial data at $t=t_n$  is given as
\be
u(\bx,t_n)=\phi_1^n(\bx)=O(1),\qquad \partial_tu(\bx,t_n)=\frac{1}{\eps^2}\phi_2^n(\bx)=O(\eps^{-2}),\qquad \bx\in{\mathbb R}^d,
\ee
and decomposing the solution $u(\bx,t)=u(\bx,t_n+s)$ of the NKGE \eqref{KG}
on the time interval $t\in[t_n, t_{n+1}]$ as \cite{BZ,Machihara, Masmoudi}
\begin{equation}\label{ansatz}
u(\bx,t_n+s)=\fe^{is/\eps^2}z^n(\bx,s)+\fe^{-is/\eps^2}\overline{z^n}(\bx,s)
+r^n(\bx,s), \qquad \bx\in{\mathbb R}^d,\quad  0\leq s\leq\tau,
\end{equation}
then a multiscale decomposition by the $\eps$-frequency (MDF) of the NKGE (\ref{KG}) can be given  as \cite{BZ,BDZSISC}
\be\label{MDF}
\left\{\begin{split}
&2i\partial_s z^n(\bx,s)+\eps^2\partial_{ss}z^n(\bx,s)-\Delta z^n(\bx,s)+3\lambda |z^n(\bx,s)|^2z^n(\bx,s)=0,\\
&\eps^2\partial_{ss}r^n(\bx,s)-\Delta r^n(\bx,s)+\frac{1}{\eps^2}r^n(\bx,s)+f_r\left(z^n(\bx,s),r^n(\bx,s);s\right)=0,
\end{split}
\right.
\quad x\in {\mathbb R}^d,\quad  0\leq s\leq\tau,
\ee
with the well-prepared initial data for $z^n$ and small initial data for $r^n$ as \cite{BZ,BDZSISC}
\begin{equation}\label{FSW-i21}
\left\{
  \begin{split}
&z^n(\bx,0)=\frac{1}{2}\left[\phi_1^n(\bx)-i\phi_2^n(\bx)\right],\quad
\partial_s z^n(x,0)=\frac{i}{2}\left[-\Delta z^n(\bx,0)+3\lambda |z^n(\bx,0)|^2z^n(\bx,0)\right],\\
&r^n(\bx,0)=0, \qquad \partial_sr^n(\bx,0)=-\partial_sz^n(\bx,0)-\partial_s
\overline{z^n}(\bx,0),\qquad \qquad \bx\in{\mathbb R}^d,
\end{split}
  \right.
\end{equation}
where
\be
\begin{split}
f_r\left(z,r;s\right)=&\lambda\fe^{3is/\eps^2}z^3+
\lambda\fe^{-3is/\eps^2}\overline{z}^3+3\lambda \left(\fe^{2is/\eps^2}z^2+ \fe^{-2is/\eps^2}\overline{z}^2\right)r
+3\lambda \left(\fe^{is/\eps^2}z+\fe^{-is/\eps^2}\overline{z}\right)r^2\nonumber\\
&+6\lambda|z|^2r+\lambda r^3.
\end{split}
\ee
\begin{figure}[t!]
\centerline{\psfig{figure=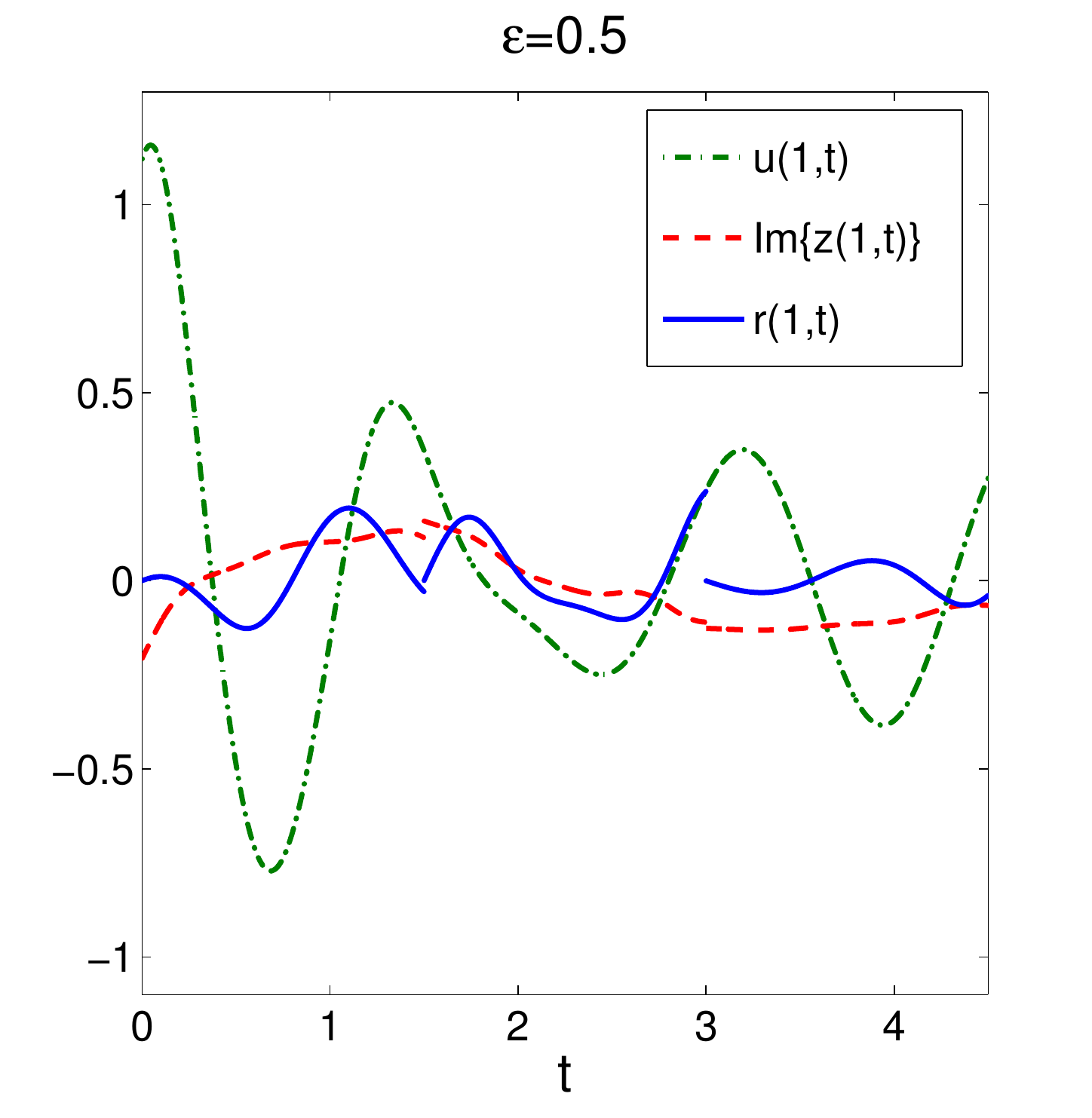,height=6cm,width=8cm}
\psfig{figure=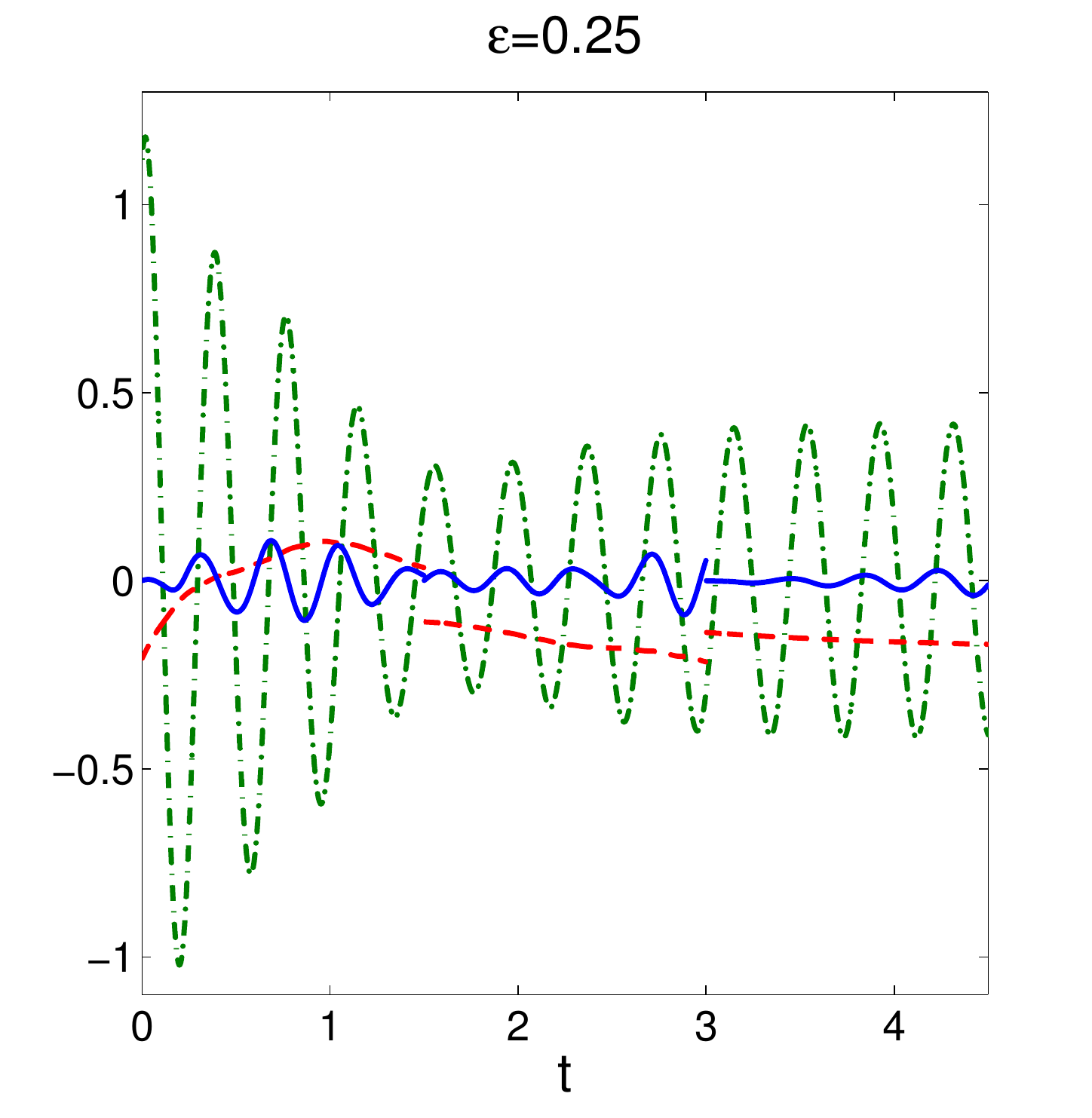,height=6cm,width=8cm}}
\caption{Decomposition of the solution $u(1,t)$ of \eqref{KG} and \eqref{initKG} with $d=1$ and $\lambda=1$ via the MDF \eqref{MDF} with $\tau=0.5$ for different $\eps$: $z(1,t)=z^n(1,s)$ and
$r(1,t)=r^n(1,s)$ if $t=n\tau+s$.}\label{fig:MTI}
\end{figure}

Then the problem \eqref{MDF} with \eqref{FSW-i21} is truncated on a bounded domain with periodic boundary condition and then discretized by the EWI-FP method \cite{BZ} with details omitted here for brevity.
After solving numerically the decomposed problem \eqref{MDF} with \eqref{FSW-i21}, the solution of the NKGE \eqref{KG} at $t=t_{n+1}$ is reconstructed by the ansatz (\ref{ansatz}) by setting $s=\tau$ \cite{BZ}.

For the convenience of the reader and simplicity of notation,
here we only present the method in 1D on $\Omega=(a,b)$ with the periodic  boundary condition.
Let $u_j^n$ and $\dot{u}_j^n$ be the approximations of $u(x_j,t_n)$ and $\partial_tu(x_j,t_n)$, respectively; and let $z_{j}^{n+1},$ $\dot{z}_{j}^{n+1},$ $r_j^{n+1}$ and $\dot{r}_j^{n+1}$ be the approximations of $z^n(x_j,\tau),$ $\partial_sz^n(x_j,\tau),$ $r^n(x_j,\tau)$ and $\partial_s r^n(x_j,\tau)$,  respectively, for $j=0,1, \ldots,N$ and $n\ge0$. Choosing $u_j^0=\phi_1(x_j)$ and $\dot{u}_j^0=\phi_2(x_j)/\eps^2$,
then a multiscale time integrator Fourier pseudospectral (MTI-FP) method
for discretizing the NKGE \eqref{KG trun} is given as \cite{BZ}
\be\label{MTI-FP}
\left\{
\begin{split}
&u^{n+1}_j=\fe^{i\tau/\eps^2}z_{j}^{n+1}+\fe^{-i\tau/\eps^2}
\overline{z_{j}^{n+1}}+r^{n+1}_j,\qquad \qquad \qquad j=0,1,\ldots,N, \quad n\ge0, \\
&\dot{u}^{n+1}_j=\fe^{i\tau/\eps^2}\left(\dot{z}_{j}^{n+1}+\frac{i}{\eps^2}z_{j}^{n+1}\right)
+\fe^{-i\tau/\eps^2}\left(\overline{\dot{z}_{j}^{n+1}}-
\frac{i}{\eps^2}\overline{z_{j}^{n+1}}\right)+\dot{r}^{n+1}_j,
\end{split}
\right.
\ee
where $z^{n+1}, $ $\dot{z}^{n+1}$, $r^{n+1}$ and $\dot{r}^{n+1}$
are numerical approximations of \eqref{MDF} with \eqref{FSW-i21} as
\be
\left\{
\begin{split}
&z_{j}^{n+1}=\sum_{l=-N/2}^{N/2-1}\widetilde{(z^{n+1})}_l\,\fe^{i\mu_l(x_j-a)},
\quad
r_{j}^{n+1}=\sum_{l=-N/2}^{N/2-1}\widetilde{(r^{n+1})}_l\,\fe^{i\mu_l(x_j-a)},\\
&\dot{z}_{j}^{n+1}=\sum_{l=-N/2}^{N/2-1}\widetilde{(\dot{z}^{n+1})}_l\, \fe^{i\mu_l(x_j-a)},\quad
\dot{r}_{j}^{n+1}=\sum_{l=-N/2}^{N/2-1}\widetilde{(\dot{r}^{n+1})}_l\,
\fe^{i\mu_l(x_j-a)},
\end{split}
\right.
\ee
with
\be\label{MTI-FP E}
\left\{
\begin{split}
&\widetilde{(z^{n+1})}_l=a_l(\tau)\widetilde{(z^{(1)})}_l+\eps^2b_l(\tau)
\widetilde{(\dot{z}^{(1)})}_l
-c_l(\tau)\widetilde{(\eta^{(1)})}_l-d_l(\tau)\widetilde{(\dot{\eta}^{(1)})}_l,\\
&\widetilde{(\dot{z}^{n+1})}_l= a_l'(\tau)\widetilde{(z^{(1)})}_l+\eps^2b_l'(\tau)\widetilde{(\dot{z}^{(1)})}_l
-c_l'(\tau)\widetilde{(\eta^{(1)})}_l-d_l'(\tau)\widetilde{(\dot{\eta}^{(1)})}_l,\quad l=-N/2,\ldots,N/2-1,\\
&\widetilde{(r^{n+1})}_l=\frac{\sin(\omega_l \tau)}{\omega_l}\widetilde{(\dot{r}^{(1)})}_l
-p_l(\tau)\widetilde{(g^{(1)})}_l-q_l(\tau)\widetilde{(\dot{g}^{(1)})}_l
-\overline{p_l}(\tau)\widetilde{\left(\overline{g^{(1)}}\right)}_l
-\overline{q_l}(\tau)\widetilde{\left(\overline{\dot{g}^{(1)}}\right)}_l,\\
&\widetilde{(\dot{r}^{n+1})}_l=\cos(\omega_l \tau)\widetilde{(\dot{r}^{(1)})}_l-p_l'(\tau)\widetilde{(g^{(1)})}_l
-q_l'(\tau)\widetilde{(\dot{g}^{(1)})}_l
-\overline{p_l'}(\tau)\widetilde{\left(\overline{g^{(1)}}\right)}_l
-\overline{q_l'}(\tau)\widetilde{\left(\overline{\dot{g}^{(1)}}\right)}_l
-\frac{\tau}{2\eps^2}\widetilde{(w^{n+1})}_l,
\end{split}
\right.
\ee
and
\be
\left\{
\begin{split}
&z_{j}^{(1)}=\frac{1}{2}\left(u_j^n-i\eps^2\dot{u}_j^n\right),\quad
\eta_{j}^{(1)}=3\lambda\left|z_{j}^{(1)}\right|^2z_{j}^{(1)},\quad
g_{j}^{(1)}=\lambda \left(z_{j}^{(1)}\right)^3,\quad
\dot{z}_{j}^{(1)}=\frac{i}{2}\left[-\partial_{xx}^{\mathcal F}z^{(1)}_j+\eta_{j}^{(1)}\right],\\ &\dot{r}_j^{(1)}=-\dot{z}_{j}^{(1)}-\overline{\dot{z}_{j}^{(1)}},\quad
\dot{\eta}_{j}^{(1)}=6\lambda z_{j}^{(1)}\cdot Re\left(\overline{z_{j}^{(1)}}\dot{z}_{j}^{(1)}\right)
+3\lambda\dot{z}_{j}^{(1)}|z_{j}^{(1)}|^2,\quad
\dot{g}_{j}^{(1)}=3\lambda \left(z_{j}^{(1)}\right)^2\dot{z}_{j}^{(1)},\\
&w_j^{n+1}=3\lambda r_j^{n+1}\left(\fe^{2i\tau/\eps^2}(z_j^{n+1})^2+ \fe^{-2i\tau/\eps^2}(\overline{z_j^{n+1}})^2\right)
+3\lambda (r_j^{n+1})^2\left(\fe^{i\tau/\eps^2}z_j^{n+1}+\fe^{-i\tau/\eps^2}\overline{z_j^{n+1}}\right)\\
&\qquad\quad\ +6\lambda|z_j^{n+1}|^2r_j^{n+1}
+\lambda (r_j^{n+1})^3,\quad
j=0,1,\ldots,N, \quad n\ge0.
\end{split}
\right.
\ee
Here we adopt the following functions \cite{BZ}
\begin{equation*}
\begin{split}
&a_l(s):=\frac{\lambda^+_l\fe^{is\lambda^-_l}-\lambda^-_l
\fe^{is\lambda^+_l}}{\lambda^+_l-\lambda^-_l},\qquad b_l(s):=i\frac{\fe^{is\lambda^+_l}-
\fe^{is\lambda^-_l}}{\eps^2(\lambda^-_l-\lambda^+_l)},\qquad
\lambda^\pm_l=-\frac{1\pm\sqrt{1+\mu_l^2\eps^2}}{\eps^2},\\
&c_l(s):=\int_0^{s}b_l(s-\theta)\, d\theta,\qquad d_l(s):=\int_0^{s}b_l(s-\theta)\theta\, d\theta,\qquad 0\le s\le \tau,\\
&p_l(s):=\int_0^s\frac{\sin\left(\omega_l(s-\theta)\right)}
{\eps^2\omega_l}\fe^{3i\theta/\eps^2}\,d\theta,\qquad
q_l(s):=\int_0^s\frac{\sin\left(\omega_l(s-\theta)\right)}
{\eps^2\omega_l}\fe^{3i\theta/\eps^2}\theta\,d\theta.
\end{split}
\end{equation*}

The MTI-FP is explicit, unconditionally stable,
and its memory cost is $O(N)$ and
computational cost per step is $O(N\,\ln N)$.
As established in \cite{BZ},
under proper regularity of the solution
$u$ of the NKGE (\ref{KG trun}),
the following two error bounds
were established by using two different techniques
for the MTI-FP method \cite{BZ}
\be
\|u(\cdot,t_n)-I_Nu^n\|_{H^2}\lesssim h^{m_2}+\tau^2+\eps^2,\quad
\|u(\cdot,t_n)-I_Nu^n\|_{H^2}\lesssim h^{m_2}+\frac{\tau^2}{\eps^2},\qquad n=0,1,\ldots,\frac{T}{\tau},
\ee
which imply a uniform error bound for $\eps\in (0,1]$ \cite{BZ}
\be
\|u(\cdot,t_n)-I_Nu^n\|_{H^2}\lesssim h^{m_2}+\max_{0<\eps\le1}\min\left\{\eps^2,\frac{\tau^2}{\eps^2}\right\}
\lesssim h^{m_0}+\tau,\qquad 0\le n\le \frac{T}{\tau},
\ee
where $m_2\ge1$ depends on the regularity of the solution $u$ of
(\ref{KG trun}).

These error bounds suggest that the MTI-FP method is uniformly
spectral order in space if the solution is smooth and is uniformly
first order in time for $0<\eps\le 1$. The $\eps$-resolution
is $h=O(1)$ and $\tau=O(1)$ in the nonrelativistic limit regime,
which immediately show that the MTI-FP is {\bf super-resolution} in time
with respect to $\eps\in(0,1]$
since the time step can be chosen as independently of $\eps$ although
the solution is highly oscillatory in time at wavelength $O(\eps^2)$.

\subsection{Two-scale formulation (TSF) method}
As presented in \cite{Chartier}, the TSF method was constructed
by separating the slow time scale $t$ and the fast time scale $\xi=t/\eps^2$ and re-formulating the NKGE \eqref{KG} into a two-scale formulation. This approach offers a general strategy to design uniformly accurate schemes for highly oscillatory differential equations and PDEs which contain general nonlinearity or strong couplings \cite{VPzhao0,VPzhao1,VPzhao2}.

Introducing \cite{Chartier}
\be\label{vbxt}
v:=v(\bx,t)=u(\bx,t)-i\eps^2(1-\eps^2\Delta)^{-1/2}\partial_tu(\bx,t)
\Rightarrow  u(\bx,t)=\frac{1}{2}\left[v(\bx,t)+\overline{v}(\bx,t)\right],
\quad \bx\in{\mathbb R}^d,\  t\ge0,
\ee
then the NKGE (\ref{KG}) can be re-written as a first order PDE
\begin{equation}\label{KG1}
\left\{\begin{split}
&i\partial_tv(\bx,t)=-\frac{1}{\eps^2}(1-\eps^2\Delta)^{1/2}v(\bx,t)
-\frac{\lambda}{8}(1-\eps^2\Delta)^{-1/2}
\left[v(\bx,t)+\overline{v}(\bx,t)\right]^3,\qquad \bx\in{\mathbb R}^d,\quad t>0,\\
&v(\bx,0)=v_0(\bx):=\phi_1(\bx)-i(1-\eps^2\Delta)^{-1/2}\phi_2(\bx),\qquad \bx\in{\mathbb R}^d.
\end{split}\right.
\end{equation}

Let
\be \label{wbxt}
w:=w(\bx,t)=\fe^{-it/\eps^2\sqrt{1-\eps^2\Delta}}v(\bx,t) \Leftrightarrow
v(\bx,t)=\fe^{it/\eps^2\sqrt{1-\eps^2\Delta}}w(\bx,t),
\qquad \bx\in{\mathbb R}^d,\quad t\ge0,
\ee
so as to filter out the main oscillation in the above PDE. Then one gets \cite{Chartier}
\begin{equation}\label{kg two scale pre}
\left\{
\begin{split}
&\partial_t w=\frac{i\lambda}{8}(1-\eps^2\Delta)^{-1/2}
\fe^{-it/\eps^2\sqrt{1-\eps^2\Delta}}
\left[\fe^{it/\eps^2\sqrt{1-\eps^2\Delta}}
w+\fe^{-it/\eps^2\sqrt{1-\eps^2\Delta}}
\overline{w}\right]^3, \quad \bx\in{\mathbb R}^d,\ t>0,\\
&w(\bx,0)=v_0(\bx),\qquad \bx\in{\mathbb R}^d.
\end{split}\right.
\end{equation}

Introduce $U:=U(\bx,t,\xi)$ with $\xi$ interpreted as another `space' variable on torus $\mathbb{T}=\mathbb{R}/(2\pi\mathbb{Z})$ such that
\be \label{wbxt1}
w(\bx,t)=U\left(\bx,t,\frac{t}{\eps^2}\right), \qquad \bx\in{\mathbb R}^d,\quad t\ge0,
\ee
with $t$ the slow time variable and
$\xi=t/\eps^2$ the fast time variable \cite{Chartier}.
Noticing \eqref{kg two scale pre}, one needs to request $U$ satisfies
the following PDE \cite{Chartier}
\be\label{KG two scale}
\left\{
\begin{split}
&\partial_t{U}(\bx,t,\xi)+\frac{1}{\eps^2}\partial_\xi U(\bx,t,\xi)=F(t,\xi,U(x,t,\xi)),\qquad  \bx\in{\mathbb R}^d,\quad t>0,\quad  \xi\in\mathbb{T},\\
&{U}(\bx,0,\xi)=U_0(\bx,\xi),\qquad \bx\in{\mathbb R}^d,\quad \xi\in\bT,\\
&{U}(\bx,t,\xi)=U(\bx,t,\xi+2\pi),\qquad \bx\in{\mathbb R}^d,\quad t\geq0,
 \quad \xi\in \bT,\\
\end{split}\right.
\ee
where $U_0(\bx,\xi)$ to be determined later satisfying $U_0(\bx,0)=v_0(\bx)$ and
\begin{equation}\label{F def}
F(t,\xi,\phi):=\frac{i\lambda}{8}(1-\eps^2\Delta)^{-1/2}\fe^{-i\xi}
\fe^{-itD_\eps}\left[\fe^{i\xi}\fe^{itD_\eps}\phi+\fe^{-i\xi}\fe^{
-itD_\eps}\overline{\phi}\right]^3,\ D_\eps=\frac{1}{\eps^2}\left[\sqrt{1-\eps^2\Delta}-1\right],
\end{equation}
with $\phi:=\phi(\bx)$.

 The initial data $U_0(\bx,\xi)$ in (\ref{KG two scale}) is only prescribed
at one point, i.e. $\xi=0$, so there is some freedom to choose the initial data in order to bound the time derivatives of $U$. By using
the Chapman-Enskog expansion, the initial data  $U_0(\bx,\xi)$ was obtained at different order of accuracy in term of
$\eps$ \cite{Chartier}. For example,
the initial data at first order of accuracy was given as
\cite{Chartier}
\be\label{1st data}
U_0(\bx,\xi):=U_0^{\rm 1st}(\bx,\xi)
=v_0(\bx)+G_1(\xi,v_0)-G_1(0,v_0), \qquad \bx\in{\mathbb R}^d, \quad 0\le
\xi\le 2\pi,
\end{equation}
such that $\partial_t^2U(x,t,\xi)=O(1)$ for fixed $t\geq0$ as $\eps\to0$ \cite{Chartier}, where
\[
G_1(\xi,\phi)=\eps^2 {\mathcal A}\,F(0,\xi,\phi) \quad \hbox{with}
\quad {\mathcal A}:={\mathcal L}^{-1}({\mathcal I}-\Pi),
\]
with  the operators ${\mathcal L}$ and $\Pi$ defined as
\[
{\mathcal L}\varphi(\xi)=\partial_\xi \varphi(\xi),\quad \Pi \varphi=\frac{1}{2\pi}\int_0^{2\pi}\varphi(\xi)\ud\xi, \quad \mbox{and}
\quad {\mathcal L}^{-1}\varphi=({\mathcal I}- \Pi)\int_0^{\xi}\varphi(\theta)\ud \theta\ \ \mbox{when}
\ \ \Pi \varphi=0,
\]
for some periodic function $\varphi:=\varphi(\xi)$ on $\bT$.

Let $U^n(\bx,\xi)$ be the numerical approximation of $U(\bx,t_n,\xi)$
for $n\ge0$. Then \eqref{KG two scale} with \eqref{1st data} can be discretized in time as
\begin{equation}\label{1st scheme}
U^{n+1}(\bx,\xi)=U^n(\bx,\xi)+\tau\, F\left(t_n,\xi,U^n(\bx,\xi)\right)-\frac{\tau}{\eps^2}\partial_\xi U^{n+1}(\bx,\xi),\quad \bx\in{\mathbb R}^d,\quad
\xi\in\bT, \quad n\ge0,
\end{equation}
with $U^n(\bx,\xi)=U_0^{\rm 1st}(\bx,\xi)$.
In the $\xi$-direction, thanks to the periodicity,
one can further discretize (\ref{1st scheme}) by the Fourier pseudospectral method as: Let $h_\xi=2\pi/N_\xi$ with $N_\xi$ an even positive integer,   $\xi_m=mh_\xi$ and $U^{n}_{m}(\bx)$ be the approximation of  $U^{n}(\bx,\xi_m)$ for $m=0,1,\ldots,N_\xi$, denote
$U^n(\bx)=(U_0^{n}(\bx),U^{n}_1(\bx),\ldots,$ $U^{n}_{N_\xi}(\bx))^T$,
take $U^{0}_{m}(\bx)=U_0^{\rm 1st}(\bx,\xi_m)$ for $m=0,1,\ldots,N_\xi$, then
one can get
\be \label{Unp1}
U^{n+1}_{m}(\bx)=\sum_{l=-N_\xi/2}^{M_\xi/2-1} \widetilde{U}^{n+1}_l(\bx)
\,\fe^{il\xi_m}, \quad m=0,1,\ldots,N_\xi,\qquad n\ge0,
\ee
with  $\widetilde{U}^n_l(\bx)=\sum_{m=0}^{N_\xi-1}U^n_{m}(\bx)\,\fe^{-il\xi_m}$ and
\be\label{Unp2}
\widetilde{U}^{n+1}_l(\bx)=\frac{\widetilde{U}^n_l(\bx)+\tau \widetilde{F}^n_l(\bx)}{1+il\tau/\eps^2},\quad
\widetilde{F}^n_l(\bx)=
\sum_{m=0}^{N_\xi-1}F\left(t_n,\xi_m,U^n_m(\bx)\right)\,\fe^{-il\xi_m},
\quad l=-\frac{N_\xi}{2},\ldots, \frac{N_\xi}{2}-1.
\ee
Then \eqref{Unp1} with \eqref{Unp2} will be first truncated (in $\bx$) on a bounded computational domain with periodic boundary condition and
then discretized by the standard Fourier pseudospectral method with details omitted here for brevity \cite{Chartier}.
Finally, noticing \eqref{vbxt}, \eqref{wbxt} and \eqref{wbxt1},
one can reconstruct the approximation of the solution $u$ of the NKGE \eqref{KG} (or \eqref{KG trun}). For the simplicity of notations, here we only present a first order two-scale formulation Fourier pseudospectral (TSF-FP1)
method in 1D as \cite{Chartier}:
\be \label{ujn1rc}
u_j^{n+1}=\frac{1}{2}\left[v_{j}^{n+1}+\overline{v}_{j}^{n+1}\right],\quad v_j^{n+1}=\fe^{\frac{it_{n+1}}{\eps^2}\sqrt{1-\eps^2\partial_{xx}^{\mathcal F}}}\, \left(I_{N_\xi} \left. U^{n+1}(x_j)\right)\right|_{\xi=\frac{t_{n+1}}{\eps^2}},\quad 0\le j\le N, \
n\ge0,
\ee
where $u_j^0=\phi_1(x_j)$ for $j=0,1,\ldots,N$.

Similarly, by taking the initial data as
\begin{align}\label{2nd data}
U_0(\bx,\xi)=U_0^{\rm 2nd}(\bx,\xi):=v_0(\bx)+G_1\left(\xi,U_0^{\rm 1st}(\bx)\right)
-G_1\left(0,U_0^{\rm 1st}(\bx)\right)+G_2(\xi,v_0(\bx))-G_2(0,v_0(\bx)),
\end{align}
such that $\partial_t^3U(x,t,\xi)=O(1)$ for fixed $t\geq0$ as $\eps\to0$ \cite{Chartier}, where
\[G_2(\xi,\phi):=-\eps^2 {\mathcal A}^2\left[\partial_t F(0,\xi,\phi)+\partial_\phi F(0,\xi,\phi)
 \Pi\left(F(0,\xi,\phi)\right)\right].\]
Then \eqref{KG two scale} with \eqref{2nd data} can be discretized by a second order scheme in time as
\be\label{2nd scheme}
\left\{
\begin{split}
&U^{n+1/2}(\bx,\xi)=U^n(\bx,\xi)+\frac{\tau}{2}\mathcal{F}\left(t_n,\xi,
U^n(\bx,\xi)\right)-\frac{\tau}{2\eps^2}\partial_\xi U^{n+1/2}(\bx,\xi),
\quad \bx\in{\mathbb R}^d,\ \xi\in\bT,\\
&U^{n+1}(\bx,\xi)=U^n(\bx,\xi)+\tau\mathcal{F}\left(t_{n+1/2},
\xi,U^{n+1/2}(\bx,\xi)\right)
-\frac{\tau}{2\eps^2}\partial_\xi\left(U^{n+1}(\bx,\xi)+U^{n}(\bx,\xi)\right),\\
\end{split}
\right.
\ee
with $U^n(\bx,\xi)=U_0^{\rm 2nd}(\bx,\xi)$. Similarly, \eqref{2nd scheme}
can be discretized in $\xi$-direction via the Fourier pseudospectral method,
truncated in $\bx$-direction onto a bounded computational domain with periodic boundary condition and then discretized via the Fourier pseudospecral method with details omitted here for brevity \cite{Chartier}.
Finally one can obtain a second order two-scale formulation Fourier pseudospectral (TSF-FP2) for the NKGE \eqref{KG} (or \eqref{KG trun}) via the reconstruction
\eqref{ujn1rc}.

As shown in \cite{Chartier},
both TSF-FP1 and TSF-FP2 are explicit, unconditionally stable,
and its memory cost is $O(N_\xi N)$ and
computational cost per step is $O(N_\xi N\,\ln (N_\xi N))$.
Under proper regularity of the solution
$U$ of the PDE (\ref{KG two scale}),
the following error bound
was established for TSF-FP1 \cite{Chartier}
\be
\|u(\cdot,t_n)-I_Nu^n\|_{H^1}\lesssim h^{m_0}+h_\xi^{m_1}+\tau,\qquad
n=0,1,\ldots,\frac{T}{\tau},
\ee
and respectively, for TSF-FP2 \cite{Chartier}
\be
\|u(\cdot,t_n)-I_Nu^n\|_{H^1}\lesssim h^{m_0}+h_\xi^{m_1}+\tau^2,\qquad
n=0,1,\ldots,\frac{T}{\tau},
\ee
where $m_0$ and $m_1$ are two positive integers
which depend on the regularity of solution $U$
of (\ref{KG two scale}) in $\bx$- and $\xi$-direction,
respectively.

These error bounds suggest that both TSF-FP1 and TSF-FP2 methods are uniformly
spectral order in space and in $\xi$-direction
if the solution is smooth, and the TSF-FP1 and TSF-FP2 are uniformly
first and second order, respectively, in time. Again, the $\eps$-resolution
is $h=O(1)$ and $\tau=O(1)$ in the nonrelativistic limit regime,
which immediately show that both TSF-FP1 and TSF-FP2 are {\bf super-resolution} in time with respect to $\eps\in(0,1]$
since the time step can be chosen independently on $\eps$ although
the solution is highly oscillatory in time at wavelength $O(\eps^2)$. We remark that
the finite difference integrator (\ref{1st scheme}) or (\ref{2nd scheme}) from \cite{Chartier} is not the unique choice for discritizating the two-scale system (\ref{KG two scale}). The formulation (\ref{KG two scale}) and the well-prepared initial data (\ref{1st data}) or (\ref{2nd data}) are essential for the TSF  approach, and other numerical discretizations such as EWI-FP can also be applied to solve (\ref{KG two scale}), see \cite{VPzhao0,VPzhao2,ZhaoJCAM}.

\section{Uniformly and optimally accurate (UOA) methods}\label{set:UOA}
In this section, we review briefly two UA methods with optimal convergence rate in time and/or computational costs, i.e. uniformly second-order
in time without solving a problem in one more spatial dimension. One is the iterative exponential-type integrator in \cite{Bau}, and the other is a MTI based on higher order multiscale expansion by frequency in \cite{BZnew}.

\subsection{An iterative exponential integrator (IEI) }

Without using higher order approximations or extra dimensions, a second order UOA method for the NKGE \eqref{KG} (or \eqref{KG trun}) was very recently proposed in \cite{Bau} by using an iterative exponential integrator. By reformulating the NKGE (\ref{KG}) into the first order PDE (\ref{KG1}) and then introducing
\be\label{vstar1}
v_*(\bx,t)=\fe^{-it/\eps^2}v(\bx,t),\qquad \mathcal{A}_\eps=(1-\eps^2\Delta)^{-1/2},
\ee
one finds
$$i\partial_tv_*(\bx,t)=-D_\eps v_*(\bx,t)-\frac{\lambda}{8}\mathcal{A}_\eps\fe^{-it/\eps^2}
\left[\fe^{it/\eps^2}v_*(\bx,t)+\fe^{-it/\eps^2}
\overline{v_*}(\bx,t)\right]^3,\quad \bx\in\bR^d,\ t>0,$$
which based on the Duhammel's formula gives
\begin{align}
v_*(\bx,t_n+s)=&\fe^{isD_\eps}v_*(\bx,t_n)\label{EI UOA}\\
&+\frac{\lambda i\mathcal{A}_\eps}{8}
\int_0^s\fe^{i(s-\theta)D_\eps-i(t_n+\theta)/\eps^2}
\left[\fe^{i(t_n+\theta)/\eps^2}v_*(\bx,t_n+\theta)
+\fe^{-i(t_n+\theta)/\eps^2}\overline{v_*}(\bx,t_n+\theta)\right]^3d\theta.\nonumber
\end{align}
Here $D_\eps$ is defined the same as in (\ref{F def}).
As used in \cite{Bau}, an exponential integrator is proposed by plugging (\ref{EI UOA}) iteratively into the cubic terms (see this technique also in \cite{Schratz2}). To describe the scheme, the following functions and operators are introduced.
\begin{align*}
&\varphi_1(z)=\frac{\fe^{z}-1}{z},\quad \varphi_2(z)=\frac{z\fe^{z}-\fe^z+1}{z^2},\quad z\in\bC,\\
&
\mathfrak{D}_k^n=\tau\fe^{i(\tau D_\eps+2t_n/\eps^2)}\varphi_k\left(i\tau\left(
2\eps^{-2}-\Delta/2\right)\right),\quad k=1,2,\ n\geq0,\\
&\mathcal{D}_{k,m}^n=
\tau\fe^{i(\tau D_\eps-mt_n/\eps^2)}\varphi_k\left(-i\tau\left(
m\eps^{-2}+D_\eps\right)\right),\quad k=1,2,\ m=2,4,\ n\geq0.
\end{align*}

For simplicity of notations, here we only present the method in 1D.
In this case, $\Delta=\partial_{xx}$. Let $v_*^n$ be an approximation
of $v_*(x,t_n)$ for $n\ge0$ via time integration, and take $u^0(x)=\phi_1(x)$ and
$v_*^0=\phi_1(x)-i(1-\eps^2\partial_{xx})^{-1/2}\phi_2(x)$.
Then an iterative exponential integrator (IEI) scheme for approximating
\eqref{EI UOA} reads:
\begin{align}\label{UOA}
v_*^{n+1}=&\fe^{i\tau D_\eps/2 }\fe^{3\lambda i\tau|w_*^n|^2/8}w_*^n+\frac{3\lambda i\tau}{8}
(\mathcal{A}_\eps-1)\fe^{i\tau D_\eps/2}|w_*^n|^2w_*^n+\tau^2\lambda^2\kappa^n +\frac{\lambda i}{8}\mathcal{A}_\eps
\chi^n \nonumber\\
&+\frac{3i\tau\eps^2\lambda^2}{128}\mathcal{A}_\eps\left(2|v_*^n|^2\mathcal{A}_\eps\zeta_0^n+
(v_*^n)^2\mathcal{A}_\eps\overline{\zeta_0^n}\right),\quad n\geq0,
\end{align}
where
\begin{align*}
w_*^n=&\fe^{iD_\eps\tau/2}v_*^n,\quad
\Gamma_{j,k}^{n,m}=\frac{3\lambda i}{8}\tau^2\fe^{mit_n/\eps^2}(v_*^n)^j(\overline{v_*^n})^k\mathcal{A}_\eps,\quad
j,k=0,1,2,\ \ m=-4,-2,2,\ \ n\geq0,\\
\kappa^n=&\frac{9}{128}\fe^{iD_\eps \tau/2}\left[
\mathcal{A}_\eps (w_*^n)^2(\mathcal{A}_\eps-1)|w_*^n|^2\overline{w_*^n}-(\mathcal{A}_\eps-1)|w_*^n|^4w_*^n-
2\mathcal{A}_\eps |w_*^n|^2(\mathcal{A}_\eps -1)|w_*^n|^2w_*^n
\right],\\
\chi^n=&\mathfrak{D}_1^n(v^n_*)^3+i\tau\mathfrak{D}_2^n
\left[(\partial_{xx}/2-D_\eps)(v^n_*)^3+3(v^n_*)^2D_\eps
v_*^n\right]+3\mathcal{D}_{1,2}^n|v^n_*|^2\overline{v^n_*}+
\mathcal{D}_{1,4}^n(\overline{v^n_*})^3-\Gamma_{2,0}^{n,2}\mathcal{U}_2^n\\
&+3i\tau \mathcal{D}_{2,2}^n\left[(\overline{v^n_*})^2
D_\eps v_*^n-2|v_*^n|^2D_\eps\overline{v_*^n}\right]
-3i\tau \mathcal{D}_{2,4}^n(\overline{v^n_*})^2D_\eps\overline{v^n_*}
-\Gamma_{0,2}^{n,-2}\mathcal{U}_{-2}^n+2\Gamma_{1,1}^{n,-2}\mathcal{W}_{2}^n+\Gamma_{0,2}^{n,-4}\mathcal{W}_{4}^n,\\
\zeta_m^n=&\fe^{2it_n/\eps^2}\left(\varphi_1((m+2)i\tau/\eps^2)
-\varphi_1(mi\tau/\eps^2)\right)(v_*^n)^3
-3\fe^{-2it_n/\eps^2}\left(\varphi_1((m-2)i\tau/\eps^2)\right.\\
&\left.-\varphi_1(mi\tau/\eps^2)\right)|v_*^n|^2\overline{v_*^n}
-\fe^{-4it_n/\eps^2}\frac{\varphi_1((m-4)i\tau/\eps^2)
-\varphi_1(mi\tau/\eps^2)}{2}(\overline{v_*^n})^3,\ \  m=-2,0,2,4,\ n\ge0,
\end{align*}
with
\begin{align*}
\mathcal{U}_m^n=&3\varphi_2(mi\tau/\eps^2)|v_*^n|^2v_*^n-\frac{i\eps^2}{2\tau}\zeta^n_m,
\quad \mathcal{W}_m^n=3\varphi_2(mi\tau/\eps^2)|v_*^n|^2\overline{v_*^n}+\frac{i\eps^2}{2\tau}\overline{\zeta^n_m},\quad
m=-2,2,4,\ n\geq0.
\end{align*}
Combining \eqref{vbxt} and \eqref{vstar1}, a semi-discretized approximation of the NKGE \eqref{KG trun} is given as
$$u^{n+1}(x)=\frac{1}{2}\left[\fe^{it_{n+1}/\eps^2}v_*^{n+1}(x)+
\fe^{-it_{n+1}/\eps^2}\overline{v_*^{n+1}}(x)\right], \quad x\in \Omega, \ n\ge0.$$
In practice, performing all the differential operations in the above IEI method by the Fourier pseudospectral approximation with details omitted here for brevity \cite{Schratz2}, we obtain the iterative exponential integrator Fourier pseudospectral (IEI-FP) scheme.

The IEI-FP is explicit, unconditionally stable,
and its memory cost is $O(N)$ and
computational cost per step is $O(N\,\ln N)$.
As established in \cite{Schratz2}, under proper regularity of the solution $u$ of the NKGE \eqref{KG trun} \cite{Schratz2} (which is weaker than that is needed for MTI-FP or TSF-FP), the following error bound
was established for IEI-FP \cite{Schratz2}
\be
\|u(\cdot,t_n)-I_Nu^n\|_{H^1}\lesssim h^{m_0}+\tau^2,\qquad
n=0,1,\ldots,\frac{T}{\tau}.
\ee

\subsection{A higher order MTI}
Another uniform second-order  in time UOA method for solving the NKGE (\ref{KG}) was very recently presented in \cite{BZnew} via a higher order multiscale expansion of the solution of the NKGE. As obtained in \cite{BZnew}, the solution $u(\bx,t)$ of (\ref{KG}) is expanded as
\begin{equation}\label{ansatz2}
u(\bx,t)=\left[\fe^{it/\eps^2}v(\bx,t)+\frac{\eps^2\lambda}{8}
\fe^{3it/\eps^2}v(\bx,t)^3+c.c.\right]+\eps^2R(\bx,t),\quad
\bx\in\bR^d,\ t\geq0,
\end{equation}
where $c.c.$ represents the complex conjugate of the whole expression
before it within the bracket, and  $v:=v(\bx,t)$ solves the following nonlinear Schr\"{o}dinger equation with wave operator (NLSW) under well-prepared initial data \cite{Cai2,BZ,BZnew}
\begin{equation}\label{eq: v}
\left\{\begin{split}
&2i\partial_tv+\eps^2\partial_{tt}v-\Delta v+3\left(\lambda|v|^2+\frac{\eps^2\lambda^2}{8}|v|^4\right)v=0,\quad \bx\in\bR^d,\ t>0,\\
&v(\bx,0)=w_0(\bx)+\eps^2r_0(\bx)=:v_0(\bx),\quad  \partial_tv(\bx,0)=\frac{i}{2}\left(-\Delta w_0(\bx)+3\lambda|w_0(\bx)|^2w_0(\bx)\right)=:v_1(\bx),
\end{split}\right.
\end{equation}
where
\begin{align*}
&w_0(\bx)=\frac{1}{2}(\phi_1(\bx)-i\phi_2(\bx)),\quad
r_0(\bx)=\frac{\lambda}{8}\overline{w}_0(\bx)^3-\frac{\lambda}
{4}w_0(\bx)^3+i\,Re(v_1(\bx)),
\end{align*}
and $R:=R(\bx,t)$ solves the following NKGE
with small initial data \cite{BZnew}
\begin{equation}\label{eq: R}
\left\{\begin{split}
&\eps^2\partial_{tt}R-\Delta R+\frac{1}{\eps^2}R+\lambda(F_v+F_R)=0,\quad \bx\in\bR^d,\ t>0,\\
&R(\bx,0)=-\frac{\lambda\eps^2}{4}Re\left(r_1(\bx)\right)=:R_0(\bx)=O(\eps^2),\quad \bx\in\bR^d,\\
&\partial_tR(\bx,0)=-\frac{3\lambda}{4}Re\left(v_0(\bx)^2v_1(\bx)\right)
+\frac{3\lambda}{4}
Im\left(r_1(\bx)\right)=:R_1(\bx)=O(1),
\end{split}\right.
\end{equation}
with
\begin{align*}
r_1(\bx)=&r_0(\bx)\left[v_0(\bx)^2+v_0(\bx)w_0(\bx)+w_0(\bx)^2\right],\\
F_v(\bx,t)=&\left[\fe^{\frac{it}{\eps^2}}\frac{3\lambda^2}{32}\eps^2|v|^6v+\fe^{\frac{3it}{\eps^2}}
\left(\frac{3\lambda}{4}|v|^2v^3+\frac{9i}{4}v^2\partial_tv-\frac{1}{8}\Delta v^3+\eps^2\frac{3}{4}
(\partial_tv)^2+\eps^2\frac{3}{8}
v^2\partial_{tt}v+\eps^4\frac{3\lambda^3}{512}|v|^6v^3
\right)\right.\\
&\left.+\fe^{5it/\eps^2}\left(\frac{3\lambda}{8}v^5+\eps^2
\frac{3\lambda^2}{64}|v|^2v^5\right)+
\fe^{7it/\eps^2}\eps^2\frac{3\lambda^2}{64}v^7+\fe^{9it/\eps^2}
\eps^4\frac{\lambda^3}{512}v^9+c.c.\right],\\
F_R(\bx,t)=&\left[\fe^{2it/\eps^2}\left(3v^2+\frac{3\lambda}
{4}\eps^2|v|^2v^2\right)R
+3\fe^{it/\eps^2}\eps^2vR^2
+\frac{3\lambda^2}{64}\fe^{6it/\eps^2}\eps^4v^6R
+\frac{3\lambda}{8}\fe^{3it/\eps^2}\eps^4v^3R^2\right.\\
&\left.+\frac{3\lambda}{4}\fe^{4it/\eps^2}\eps^2v^4R+c.c.\right]+6|v|^2R+\eps^4R^3
+\frac{3\lambda^2}{32}\eps^4|v|^6R.\nonumber
\end{align*}
It was shown that $R(\bx,t)=O(\eps^2)$ \cite{BZnew} and thus the multiscale expansion (\ref{ansatz2}) without the last term is a higher order multiscale
expansion of the solution of the NKGE (\ref{KG}) at $O(\eps^4)$ \cite{BZnew}.
Then the decomposed problems (\ref{eq: v}) and (\ref{eq: R}) are solved
numerically by the EWI-FP method \cite{BZnew}.

Again, for the convenience of the reader and simplicity of notations,
here we only present the method in 1D on $\Omega=(a,b)$ with the periodic boundary condition.
Let $u_j^n$ and $\dot{u}_j^n$ be the approximations of $u(x_j,t_n)$ and $\partial_tu(x_j,t_n)$, respectively; and let $v_{j}^{n},$ $\dot{v}_{j}^{n},$ $R_j^{n}$ and $\dot{R}_j^{n}$ be the approximations of $v(x_j,t_n),$ $\partial_t v(x_j,t_n),$ $R(x_j,t_n)$ and $\partial_t R(x_j,t_n)$  respectively, for $j=0,1, \ldots,N$ and $n\ge0$. Choosing $u_j^0=\phi_1(x_j)$ and $\dot{u}_j^0=\phi_2(x_j)/\eps^2$, $v^0_j=v_0(x_j)$, $\dot{v}^0_j=v_1(x_j)$, $R^0_j=R_0(x_j)$ and $\dot{R}^0_j=R_1(x_j)$ for $j=0,1\ldots,N$, then a high-order multiscale time integrator Fourier pseudospectral (MTI-FP2) method
for discretizing the NKGE \eqref{KG trun} is given as \cite{BZnew}
\be\label{UOA}
\begin{split}
u_j^{n+1}=&\left[\fe^{it_{n+1}/\eps^2}v_j^{n+1}+
\frac{\lambda\eps^2}{8}\fe^{3it_{n+1}/\eps^2}(v_j^{n+1})^3+c.c.\right]
+\eps^2R_j^{n+1},\qquad \qquad j=0,1,\ldots,N, \quad n\ge0, \\
\dot{u}_j^{n+1}=&\left[\frac{i}{\eps^2}\fe^{\frac{it_{n+1}}{\eps^2}}v_j^{n+1}
+\fe^{\frac{it_{n+1}}{\eps^2}}\dot{v}_j^{n+1}
+\frac{3i\lambda}{8}\fe^{\frac{3it_{n+1}}{\eps^2}}(v_j^{n+1})^3+\frac{3\lambda}{8}
\eps^2\fe^{\frac{3it_{n+1}}{\eps^2}}(v_j^{n+1})^2\dot{v}_j^{n+1}+c.c.\right]\\
&+\eps^2\dot{R}_j^{n+1}.
\end{split}
\ee
Here $v^{n+1}$ and $\dot{v}^{n+1}$ are approximations of the NLSW
\eqref{eq: v} by an EWI-FP \cite{Cai2,BZ,BZnew} as
\be \label{EI: v}
v^{n+1}_j=\sum_{l=-N/2}^{N/2-1}\widetilde{(v^{n+1})}_l\;\fe^{i\mu_l(x_j-a)},
\qquad  \dot{v}^{n+1}_j=\sum_{l=-N/2}^{N/2-1}\widetilde{(\dot{v}^{n+1})}_l\;
\fe^{i\mu_l(x_j-a)},\quad j=0,1,\ldots,N,\quad n\ge0,
\ee
where
\[
\widetilde{(v^{n+1})}_l =\left\{\begin{array}{ll}
a_l(\tau)\widetilde{(v^0)}_l+\eps^2
b_l(\tau)\widetilde{(\dot{v}^0)}_l-c_l(\tau)
\widetilde{(g^0)}_l, &n=0,\\
 \\
 a_l(\tau)\widetilde{(v^n)}_l+
 \eps^2 b_l(\tau)\widetilde{(\dot{v}^n)}_l-c_l(\tau)
\widetilde{(g^n)}_l-\frac{d_l(\tau)}{\tau}\left[
\widetilde{(g^n)}_l-\widetilde{(g^{n-1})}_l\right], &n\ge1,\\
 \ea\right.
\]
\[
\widetilde{(\dot{v}^{n+1})}_l =\left\{\begin{array}{ll}
a^\prime_l(\tau)\widetilde{(v^0)}_l+\eps^2 b^\prime_l(\tau)\widetilde{(\dot{v}^0)}_l-c^\prime_l(\tau)
\widetilde{(g^0)}_l, &n=0,\\
 \\
 a^\prime_l(\tau)\widetilde{(v^n)}_l+\eps^2 b^\prime_l(\tau)\widetilde{(\dot{v}^n)}_l-c^\prime_l(\tau)
\widetilde{(g^n)}_l-\frac{d_l'(\tau)}{\tau}
\left[\widetilde{(g^n)}_l-\widetilde{(g^{n-1})}_l\right], &n\ge1,\\
 \ea\right.
\]
with $g^n=(g_0^n,g_1^n,\ldots, g_N^n)^T$ given as
\[g^n_j=3\left(\lambda|v^n_j|^2
+\frac{\eps^2\lambda^2}{8}|v^n_j|^4\right)v^n_j,\qquad j=0,1,\ldots,N, \quad n\ge0.\]
Similarly, here $R^{n+1}$ and $\dot{R}^{n+1}$ are approximations of the NKGE
\eqref{eq: R} by an EWI-FP \cite{BZ,BD,BZnew} as
\be \label{EI: R}
R^{n+1}_j=\sum_{l=-N/2}^{N/2-1}\widetilde{(R^{n+1})}_l\;\fe^{i\mu_l(x_j-a)},
\quad  \dot{R}^{n+1}_j=\sum_{l=-N/2}^{N/2-1}\widetilde{(\dot{R}^{n+1})}_l\;
\fe^{i\mu_l(x_j-a)},\quad j=0,1,\ldots,N,\ n\ge0,
\ee
where
\[
\begin{split}
\widetilde{(R^{n+1})}_l=&\cos(\omega_l \tau)\widetilde{(R^n)}_l+\frac{\sin(\omega_l \tau)}{\omega_l}\widetilde{(\dot{R}^n)}_l-
p_{l,3}^n(\tau)\widetilde{\left(G_3^n\right)}_l
-q_{l,3}^n(\tau)\widetilde{\left(\dot{G}_3^n\right)}_l-
\overline{p_{l,3}^n}(\tau)\widetilde{\left(\overline{G_3^n}\right)}_l
-\overline{q_{l,3}^n}(\tau)\widetilde{\left(\overline{\dot{G}_3^n}\right)}_l\nonumber\\
&-p_{l,5}^n(\tau)\widetilde{(G_5^n)}_l-q_{l,5}^n(\tau)
\widetilde{\left(\dot{G}_5^n\right)}_l-\overline{p_{l,5}^n}(\tau)
\widetilde{\left(\overline{G_5^n}\right)}_l
-\overline{q_{l,5}^n}(\tau)\widetilde{\left(\overline{\dot{G}_5^n}\right)}_l
-\frac{\tau\lambda\sin(\omega_l\tau)}{2\eps^2\omega_l}
\widetilde{(D^n)}_l,\\
\widetilde{(\dot{R}^{n+1})}_l=&-\omega_l\sin(\omega_l \tau)\widetilde{(R^n)}_l+\cos(\omega_l \tau)\widetilde{(\dot{R}^n)}_l-(p_{l,3}^n)'(\tau)\widetilde{\left(G_3^n\right)}_l
-(q_{l,3}^n)'(\tau)\widetilde{\left(\dot{G}_3^n\right)}_l
-\overline{(p_{l,3}^n)'}(\tau)\widetilde{\left(\overline{G_3^n}\right)}_l
\nonumber\\
&-\overline{(q_{l,3}^n)'}(\tau)\widetilde{\left(\overline{\dot{G}_3^n}
\right)}_l-(p_{l,5}^n)'(\tau)\widetilde{(G_5^n)}_l
-(q_{l,5}^n)'(\tau)\widetilde{\left(\dot{G}_5^n\right)}_l-
\overline{(p_{l,5}^n)'}(\tau)\widetilde{\left(\overline{G_5^n}\right)}_l
-\overline{(q_{l,5}^n)'}(\tau)\widetilde{\left(\overline{\dot{G}_5^n}\right)}_l\nonumber\\
&-\frac{\tau\lambda}{2\eps^2}
\left[\cos(\omega_l\tau)\widetilde{(D^n)}_l+\widetilde{(D^{n+1})}_l\right],
\end{split}
\]
with
\begin{align*}
G_{3,j}^n=&\frac{3\lambda}{4}|v_j^n|^2(v_j^n)^3+\frac{9i}{4}(v_j^n)^2\dot{v}_j^n-\frac{1}{8}\partial_{xx}^{\mathcal{F}} (v^n_j)^3,\quad
G_{5,j}^n=\frac{3\lambda}{8}(v_j^n)^5,\quad D_j^n=F_{v,j}^n+F_{R,j}^n,\\
\dot{G}_{3,j}^n=&\frac{3\lambda}{4}\left[4(v^n_j)^3\dot{v}_j^n\overline{v_j^n}
+(v_j^n)^4\overline{\dot{v}_j^n}\right]+\frac{9i}{4}
\left[2v^n_j(\dot{v}_j^n)^2+(v_j^n)^2\ddot{v}_j^n\right]
-\frac{3}{8}\partial_{xx}^{\mathcal{F}}(\dot{v}_j^n(v_j^n)^2),\\
\dot{G}_{5,j}^n=&\frac{15\lambda}{8}(v_j^n)^4\dot{v}_j^n,\quad \ddot{v}_j^n=-\frac{1}{\eps^2}\left[2i\dot{v}_j^n-\partial_{xx}^{\mathcal{F}}v_j^n
+3\left(\lambda|v_j^n|^2+\frac{\eps^2\lambda^2}{8}|v_j^n|^4\right)v_j^n\right],\\
\end{align*}
and
\begin{align*}
F_{v,j}^n=&
\left[\fe^{it_n/\eps^2}\frac{3\lambda^2}{32}\eps^2|v_j^n|^6v_j^n
+\fe^{3it_n/\eps^2}\left(\eps^2\frac{3}{4}v_j^n(\dot{v}_j^n)^2
+\eps^2\frac{3}{8}(v_j^n)^2\ddot{v}_j^n
+\eps^4\frac{3\lambda^3}{512}|v_j^n|^6(v_j^n)^3\right)\right.
\\
&+\left.\fe^{5it_n/\eps^2}\frac{3\lambda^2}{64}\eps^2|v_j^n|^2(v_j^n)^5+
\fe^{7it_n/\eps^2}\frac{3\lambda^2}{64}\eps^2(v_j^n)^7
+\fe^{9it_n/\eps^2}\frac{\lambda^3}{512}\eps^4(v_j^n)^9+c.c.\right],\\
F_{R,j}^n=&\left[\fe^{2it_n/\eps^2}\left(3(v_j^n)^2
+\frac{3\lambda}{4}\eps^2|v_j^n|^2(v_j^n)^2\right)R_j^n
+3\fe^{it_n/\eps^2}\eps^2v_j^n(R_j^n)^2
+\frac{3\lambda^2}{64}\fe^{6it_n/\eps^2}\eps^4(v_j^n)^6R_j^n\right.\\
&+\left.\frac{3\lambda}{8}\fe^{3it_n/\eps^2}\eps^4(v_j^n)^3(R_j^n)^2
+\frac{3\lambda}{4}\fe^{4it_n/\eps^2}\eps^2(v_j^n)^4R_j^n+c.c.\right]
+6|v_j^n|^2R_j^n+\eps^4(R_j^n)^3
+\frac{3\lambda^2}{32}\eps^4|v_j^n|^6R_j^n.
\end{align*}
Here we adopt the following functions from \cite{BZnew}: for $n\geq0,\,l=-N/2,\ldots,N/2-1$ and $k=3,5$,
\begin{align*}
&p_{l,k}^n(s)=\int_0^s\frac{\lambda\sin(\omega_l(s-\theta))}{\eps^2\omega_l}\fe^{ki(t_n+\theta)/\eps^2}d\theta,
\quad
q_{l,k}^n(s)=\int_0^s\frac{\lambda\sin(\omega_l(s-\theta))}{\eps^2\omega_l}\fe^{ki(t_n+\theta)/\eps^2}\theta d\theta.
\end{align*}

Again, the MTI-FP2 is explicit, unconditionally stable,
and its memory cost is $O(N)$ and
computational cost per step is $O(N\,\ln N)$.
As established in \cite{BZnew}, under proper regularity of the solution $u$ of the NKGE \eqref{KG trun}, the following error bound
was established for MTI-FP2 \cite{BZnew}
\begin{equation}\label{thm u}
\|u(\cdot,t_n)-I_Nu^n\|_{H^1}+\eps^2\|\partial_tu(\cdot,t_n)-I_N\dot{u}^n\|_{H^1}\lesssim \tau^2+h^{m_0},\quad n=0,1,\ldots\frac{T}{\tau}.
\end{equation}

\bigskip
We remark here that two new UOA schemes named as micro-macro method and pull-back method were proposed very recently in \cite{NUA} based on the averaging theory. The accuracy of the micro-macro method is very similar to TSF-FP2 and the efficiency is very similar to MTI-FP2 and IEI-FP. The pull-back is implicit but with superior long-time behaviour over other methods \cite{NUA}.

\section{Numerical comparisons and results}\label{sec:result}
In this section, we report the performance of different numerical methods reviewed in previous sections and carry out a systematical comparison.

In order to do so, we take $d=1$ and $\lambda=1$ in (\ref{KG}) and
choose the initial data as
\[\phi_1(x)= \frac{3\sin(x)}{\fe^{x^2/2}+\fe^{-x^2/2}}, \quad \phi_2(x)= \frac{2\fe^{-x^2}}{\sqrt{\pi}},\quad x\in {\mathbb R}.\]
The problem is solved numerically on a bounded computational domain $\Omega=(-16,16)$.
The `exact' solution of the NKGE \eqref{KG trun} is obtained numerically by TSF-FP2 with a very small step size, i.e. $\tau=10^{-6},h=1/64, h_\xi=\pi/64$.
Define error
$$\fe^{\tau,h}_{\eps}(t=t_n):=\left\|P_Nu(\cdot,t_n)-I_Nu^n\right\|_{H^1},$$
where $P_N$ is the standard projection operator \cite{ST}.
We depict the errors at $t=1$. The temporal error and spatial error of each numerical method are studied and shown separately in the following.

\subsection{Spatial errors}
We first test and compare the spatial discretization error of different numerical methods.
For spatial error analysis, the time step $\tau$ is chosen small enough such that the discretization error in time is negligible, e.g. $\tau=10^{-6}$.

\begin{table}[t!]
\tabcolsep 0pt
\caption{Spatial error of ECFD for different $\eps$ at time $t=1$ under $\tau=10^{-5}$.}\label{tab:spaceFD}
\begin{center}
\def\temptablewidth{1\textwidth}
{\rule{\temptablewidth}{0.75pt}}
\begin{tabular*}{\temptablewidth}{@{\extracolsep{\fill}}ccccccc}
$\fe_{\eps}^{\tau,h}(t=1)$   &  $h_0=0.5$   & $h_0/2$  & $h_0/4$ & $h_0/8$ & $h_0/16$
& $h_0/32$ \\[0.25em]
\hline
$\eps_0=1.0$       &3.10E-1            &8.37E-2          &2.09E-2        &5.31E-3 &1.24E-3 &3.11E-4 \\
rate	&---  		 &1.89		  &2.00		&1.98 		&2.09 		&2.00 \\ \hline
$\eps_0/2$       &4.26E-1	           & 1.19E-1	     & 3.11E-2       &7.92E-3 &1.83E-3 &4.58E-4 \\
rate	&---		 &1.84		  &1.94		&1.97 &2.11 &2.00 \\  \hline
$\eps_0/2^2$     & 6.67E-1	       & 2.16E-1	     & 5.71E-2       &1.41E-2 &3.30E-3 &8.26E-4 \\
rate	&---  	 &1.63		  &1.92		&2.01 &2.09 &2.00 \\  \hline
$\eps_0/2^3$    & 9.05E-1	       & 2.77E-1	     & 7.58E-2       &1.98E-2  &4.46E-3 &1.11E-3\\
rate	&--- 		 &1.70		  &1.87		&1.94 &2.14 &2.00 \\  \hline
$\eps_0/2^4$      & 8.48E-1	       & 3.00E-1	     & 8.37E-2       &2.21E-2 &5.11E-3 &1.21E-3 \\
rate	&--- 		 &1.50		  &1.84		&1.92 &2.10 &2.07 \\
\end{tabular*}
{\rule{\temptablewidth}{0.75pt}}
\end{center}
\end{table}

\begin{table}[t!]
  \caption{Spatial error of LI-FP1 for different $\eps$ at time $t=1$ under $\tau=10^{-7}$.}\label{tab:spaceAPFP1}
  \vspace*{-10pt}
\begin{center}
\def\temptablewidth{1\textwidth}
{\rule{\temptablewidth}{0.75pt}}
\begin{tabular*}{\temptablewidth}{@{\extracolsep{\fill}}cccccc}
 $\fe_{\eps}^{\tau,h}(t=1)$      & $h= 2$       & $h= 1$    &  $h/2$   & $h/4$  & $h/8$\\[0.25em]
\hline
$\eps_0=1.0$   &9.53E-1       &1.49	           &1.44	       &1.44	         &1.44\\
$\eps_0/2^2$      &9.42E-1      &4.80E-1	           &5.06E-1  	   &5.05E-1	         &5.05E-1\\
$\eps_0/2^4$     &7.48E-1     &3.04E-1	           &5.39E-2	       & 5.38E-2	     &5.38E-2\\
$\eps_0/2^6$     &1.10       &4.98E-1	           &1.98E-2	       &3.29E-3	         &3.29E-3\\
$\eps_0/2^8$     &1.08     &4.89E-1	           &1.76E-2	       & 2.07E-4	     &2.07E-4\\
$\eps_0/2^{10}$     &9.74E-1     &3.45E-1	           &2.57E-2	       &2.11E-5    	     &1.33E-5\\
$\eps_0/2^{12}$     &1.04     &4.74E-1	           &1.54E-2	       &6.77E-6   	     &8.17E-7\\
$\eps_0/2^{14}$     &7.31E-1     &2.52E-1	           &1.06E-2	       &1.34E-5   	     &7.69E-8\\
\end{tabular*}
{\rule{\temptablewidth}{0.75pt}}
\end{center}
\end{table}

\begin{table}[t!]
  \caption{Spatial error of LI-FP2 for different $\eps$ at time $t=1$ under $\tau=10^{-7}$.}\label{tab:spaceAPFP2}
  \vspace*{-10pt}
\begin{center}
\def\temptablewidth{1\textwidth}
{\rule{\temptablewidth}{0.75pt}}
\begin{tabular*}{\temptablewidth}{@{\extracolsep{\fill}}cccccc}
 $\fe_{\eps}^{\tau,h}(t=1)$     & $h= 2$       & $h= 1$    &  $h/2$   & $h/4$  & $h/8$\\[0.25em]
\hline
$\eps_0=1.0$   &7.41E-1       &10.2	           &12.9	       &14.2	         &14.2\\
$\eps_0/2^2$     &8.84E-1       &6.07E-1	           &5.86E-1  	   &6.50E-1	         &6.50E-1\\
$\eps_0/2^4$    &7.51E-1      &2.77E-1	           &2.54E-2	       &9.77E-3	         &9.77E-3\\
$\eps_0/2^6$    &1.10      &4.97E-1	           &1.90E-2	       &6.46E-5	         &6.38E-5\\
$\eps_0/2^8$     &1.08     &4.89E-1	           &1.76E-2	       &7.03E-6  	     &2.39E-7\\
$\eps_0/2^{10}$   &9.74E-1       &3.46E-1	           &2.57E-2	       &1.66E-5    	     &6.24E-9\\
$\eps_0/2^{12}$    &1.04      &4.74E-1	           &1.54E-2	       &6.73E-6   	     &7.20E-9\\
$\eps_0/2^{14}$    &7.31E-1      &2.52E-1	           &1.06E-2	       &1.34E-5   	     &4.41E-9\\
\end{tabular*}
{\rule{\temptablewidth}{0.75pt}}
\end{center}
\end{table}

\begin{table}[t!]
  \caption{Spatial error of MTI-FP for different $\eps$ at time $t=1$ under $\tau=10^{-7}$.}\label{tab:spaceMTI}
  \vspace*{-10pt}
\begin{center}
\def\temptablewidth{1\textwidth}
{\rule{\temptablewidth}{0.75pt}}
\begin{tabular*}{\temptablewidth}{@{\extracolsep{\fill}}cccccc}
$\fe_{\eps}^{\tau,h}(t=1)$     & $h= 2$        & $h= 1$    &  $h/2$   & $h/4$  & $h/8$\\[0.25em]
\hline
$\eps_0=1.0$      &5.88E-1      &2.10E-1	           &9.61E-3	       &7.58E-6	         &1.61E-11\\
$\eps_0/2$       &5.88E-1     &4.52E-1	           &2.37E-2  	   &1.70E-5	         &1.63E-11\\
$\eps_0/2^2$     &8.99E-1     &4.60E-1	           &3.05E-2	       &1.29E-5	         &1.61E-11\\
$\eps_0/2^3$    &7.07E-1      &1.57E-1	           &7.35E-3	       &6.05E-6	         &8.31E-12\\
$\eps_0/2^4$    &7.58E-1      &2.76E-1	           &2.60E-2	       &1.72E-5  	     &7.75E-12\\
$\eps_0/2^5$    &1.12      &4.66E-1	           &2.43E-2	       &1.55E-5    	     &9.33E-12\\
$\eps_0/2^8$    &1.08      &4.90E-1	           &1.64E-2	       &7.58E-6   	     &6.56E-12\\
$\eps_0/2^{11}$   &7.35E-1    &2.32E-1	           &1.64E-2	       &1.66E-5   	     &7.29E-12\\
\end{tabular*}
{\rule{\temptablewidth}{0.75pt}}
\end{center}
\end{table}

\begin{table}[t!]
  \caption{Spatial Error of TSF-FP1 in $\xi$ for different $\eps$ at time $t=1$ under $\tau=10^{-7},h= 1/16$.}\label{tab:spaceTSF1}
  \vspace*{-10pt}
\begin{center}
\def\temptablewidth{1\textwidth}
{\rule{\temptablewidth}{0.75pt}}
\begin{tabular*}{\temptablewidth}{@{\extracolsep{\fill}}cccccc}
$\fe_{\eps}^{\tau,h}(t=1)$        & $h_\xi =\pi$    &  $h_\xi/2$   & $h_\xi/4$  & $h_\xi/8$ & $h_\xi/16$\\[0.25em]
\hline
$\eps_0=1$            &4.51E-1          &2.01E-1       &1.47E-2     &4.21E-5    &9.05E-10\\
$\eps_0/2$            &4.16E-1          &1.32E-1       &4.33E-3     &1.92E-6    &1.08E-12\\
$\eps_0/2^2$          &6.33E-1          &1.27E-1       &7.46E-4     &2.59E-8    &1.34E-11\\
$\eps_0/2^3$          &6.22E-1          &1.10E-1       &1.44E-5     &2.34E-12   &4.45E-13\\
$\eps_0/2^4$          &8.14E-1          &9.84E-2       &3.70E-7     &4.12E-13   &4.12E-13\\
$\eps_0/2^5$          &9.33E-1          &1.13E-1       &3.22E-8     &5.14E-13   &3.74E-13\\
$\eps_0/2^8$          &1.07          &9.16E-2       &5.48E-12    &4.56E-13   &2.91E-13\\
$\eps_0/2^{11}$       &7.24E-1          &1.06E-1       &5.18E-13    &3.07E-13   &2.64E-13\\
\end{tabular*}
 {\rule{\temptablewidth}{0.75pt}}
\end{center}
\end{table}

The three finite difference methods share almost the same discretization error in space. Thus, we only give the spatial error of ECFD in Tab. \ref{tab:spaceFD} as a representative and omit the results of SIFD and LFFD. The errors of LI-FP1 and LI-FP2 are given in Tab. \ref{tab:spaceAPFP1} and Tab. \ref{tab:spaceAPFP2}, respectively. These errors are the spatial errors mixed with the model reduction errors. The spatial error of MTI-FP is given in Tab. \ref{tab:spaceMTI}. The results of EWI-FP, TS-FP, TSF-FP2, IEI-FP and MTI-FP2 behave similarly as that of MTI-FP since they share the same Fourier discretization, so the corresponding results have been omitted for brevity as well. The error of TSF-FP1 in the extra space direction $\xi$ is given in Tab. \ref{tab:spaceTSF1}, which represents the very similar corresponding results of TSF-FP2.

From Tabs. \ref{tab:spaceFD}-\ref{tab:spaceTSF1}, we can draw the following observations:

(i) ECFD, SIFD and LFFD have second order accuracy in space error, while EWI-FP, TS-FP, MTI-FP, TSF-FP1 and TSF-FP2 have spectral accuracy in space.  The errors are uniform in space in terms of $\eps$ with spatial $\eps$-scalability $h=O(1)$. Thus, when the initial data of NKGE is smooth enough, the Fourier pseudospectral discretization in space is obviously more efficient than finite difference discretization.

(ii) The spatial errors of the LI-FP1 and LI-FP2 are mixed with the residue of $O(\eps^2)$ and $O(\eps^4)$, respectively from the model reductions. Thus, for a fixed $0<\eps\leq1$, the spectral accuracy of space discretization is broken and the error is bounded from blow.

(iii) The errors of TSF-FP1 and TSF-FP2 in the extra space direction $\xi$ are of spectral accuracy. The error is uniformly bounded for all $\eps\in(0,1]$ which allows the use of $h_\xi=O(1)$. The computational resource needed in the $\xi$-direction is not very heavy, i.e. $N_\xi=32$ is enough to get machine accuracy for all the $\eps$ in this example. The smaller $\eps$ is, the less grid points are needed in $\xi$ direction to reach the machine accuracy.

\subsection{Temporal errors}

For temporal error analysis, the mesh size $h$ (and so is $h_\xi$ for TSF-FP) is chosen small enough such that the discretization error in space is negligible. The detailed data of the used $h,h_\xi$ is given case by case below.

The results of ECFD, SIFD, LFFD, EWI-FP, TS-FP, LI-FP1, LI-FP2, MTI-FP, TSF-FP1, TSF-FP2, IEI-FP and MTI-FP2 are shown in Tabs. \ref{tab:ECFD}-\ref{tab:MTI2}, respectively. For the LFFD method, in order to show the temporal discretization error but meanwhile to satisfy the stability condition, we choose \cite{BaoC5,BaoC6}
\begin{equation}\label{LFFD rule}
\delta_j(\eps) = \left\{\begin{array}{ll} \eps^2, &\eps_0/2^j\le \eps \le 1,\\
\eps_0^2/4^j, &0<\eps<\eps_0/2^j,\\
 \end{array} \right.  \qquad j=0,1,\ldots,
 \end{equation}
in Tab. \ref{tab:LFFD} for the temporal error. To illustrate the UA property of MTI-FP, TSF-FP1 and TSF-FP2, we also define the error
$$\fe^{\tau,h}_{\infty}(T):=\max_{\eps}\left\{\fe^{\tau,h}_\eps(T)\right\}.$$
The convergence rate of each method is shown along with the error. In the tables, we highlight the error of each classical method in the table when the time step is chosen according to its $\eps$-scalability. The efforts made here are to illustrate the convergence rate of each method and the $\eps$-scalability in the limit regime.

\begin{table}[t!]
\tabcolsep 0pt
\caption{Temporal error of ECFD for different $\eps$ at time $t=1$ under $h=1/1024$.}\label{tab:ECFD}
\begin{center}
\def\temptablewidth{1\textwidth}
{\rule{\temptablewidth}{0.75pt}}
\begin{tabular*}{\temptablewidth}{@{\extracolsep{\fill}}llllll}
 $\fe_{\eps}^{\tau,h}(t=1)$     & $\tau_0=0.2$	&$\tau_0/2^3$	&$\tau_0/2^6$	&$\tau_0/2^9$	&$\tau_0/2^{12}$\\
\hline
$\eps_0=1$	            &{\bf 2.90E-1}  &6.79E-3  &1.15E-4 &1.81E-6    &3.57E-8\\
rate                    &---	  &1.80     &1.96    &2.02       &1.89 \\  \hline
$\eps_0/2$	            &2.73  &{\bf 7.13E-2}  &1.15E-3 &2.04E-5    &3.47E-7\\
rate                    &---      &{\bf 1.76}     &1.98    &1.94      &1.95\\ \hline
$\eps_0/2^2$	        &3.16  &2.32  &{\bf 3.90E-2}  &6.28E-4     &1.51E-5\\
rate                    &---      &0.15     &{\bf 1.97}    &1.97         &1.93\\ \hline
$\eps_0/2^3$	        &6.22  &3.23  &1.73 &{\bf 2.71E-2}      &4.28E-4\\
rate                    &---      &0.32     &0.30   &{\bf 1.99}           &2.00\\ \hline
$\eps_0/2^4$	        &4.03  &7.30  &7.01 &1.61      &{\bf 2.60E-2}\\
rate                    &---      &0.29     &0.02    &0.71          &{\bf 1.98}\\
\end{tabular*}
{\rule{\temptablewidth}{0.75pt}}
\end{center}
\end{table}

\begin{table}[t!]
\tabcolsep 0pt
\caption{Temporal error of SIFD for different $\eps$ at time $t=1$ under $h=1/1024$.}\label{tab:SIFD}
\begin{center}
\def\temptablewidth{1\textwidth}
{\rule{\temptablewidth}{0.75pt}}
\begin{tabular*}{\temptablewidth}{@{\extracolsep{\fill}}llllll}
$\fe_{\eps}^{\tau,h}(t=1)$      & $\tau_0=0.2$	&$\tau_0/2^3$ &$\tau_0/2^6$	&$\tau_0/2^9$ &$\tau_0/2^{12}$\\
\hline
$\eps_0=1$	            &{\bf 2.42E-1}  &5.46E-3   &9.27E-5 &1.72E-6    &3.50E-8\\
rate                    &---	  &1.82      &1.96    &1.99       &1.87\\ \hline
$\eps_0/2$	            &2.28  &{\bf 5.83E-2}   &9.67E-4 &2.11E-5      &3.51E-7\\
rate                    &---      &{\bf 1.76}      &1.98    &1.84         &1.97\\ \hline
$\eps_0/2^2$	        &4.06  &2.07   &{\bf 3.39E-2} &5.63E-4      &8.88E-6\\
rate                    &---      &0.32      &{\bf 1.96}    &1.97         &2.00\\ \hline
$\eps_0/2^3$	        &6.05  &2.67   &1.67 &{\bf 2.66E-2}      &4.14E-4\\
rate                    &---      &0.39      &0.22    &{\bf 1.99}         &2.00\\ \hline
$\eps_0/2^4$	        &4.05  &6.78   &7.07 &1.60       &{\bf 2.60E-2}\\
rate                    &---      &-0.24     &-0.01   &0.71          &{\bf 1.99}\\
\end{tabular*}
{\rule{\temptablewidth}{0.75pt}}
\end{center}
\end{table}

\begin{table}[t!]
\tabcolsep 0pt
\caption{Temporal error of LFFD for different $\eps$ at time $t=1$ with rule (\ref{LFFD rule}).}\label{tab:LFFD}
\begin{center}
\def\temptablewidth{1\textwidth}
{\rule{\temptablewidth}{0.75pt}}
\begin{tabular*}{\temptablewidth}{@{\extracolsep{\fill}}llllll}
$\fe_{\eps}^{\tau,h}(t=1)$        & $\begin{array}{c} \tau_0=0.2\\ h_0=0.5\\ \end{array}$  & $\begin{array}{c} \tau_0/8\\ h_0/8\delta_1(\eps)\\ \end{array}$  &  $\begin{array}{c} \tau_0/8^2\\ h_0/8^2\delta_2(\eps)\\ \end{array}$   & $\begin{array}{c} \tau_0/8^3\\ h_0/8^3\delta_3(\eps)\\ \end{array}$ & $\begin{array}{c} \tau_0/8^4\\ h_0/8^4\delta_4(\eps)\\ \end{array}$\\
\hline
$\eps_0=1$	            &{\bf 2.29E-1}   &3.94E-3      &6.22E-5     &9.78E-7   &3.07E-8\\
rate                    &---	   &1.95         &1.99        &2.00      &1.66\\ \hline
$\eps_0/2$	            &6.05E-1   &{\bf 9.05E-3}      &1.45E-4     &2.27E-6   &1.68E-8\\
rate                    &---       &{\bf 2.02}         &1.99        &2.00      &2.36\\ \hline
$\eps_0/2^2$	        &unstable  &3.12E-1     &{\bf 4.93E-3}    &7.13E-5  &1.24E-6 \\
rate                    &---       &-        &{\bf 2.00}       &2.03 &1.95 \\ \hline
$\eps_0/2^3$	        &unstable      &unstable      &2.38E-1    &{\bf 3.56E-3} &6.22E-5\\
rate                    &---       &-        &-       &{\bf 2.02}  &1.95\\ \hline
$\eps_0/2^4$	        &unstable       &unstable     &2.68     &2.35E-1 &{\bf 3.64E-3}\\
rate                    &---            &-              &-          &1.17  &{\bf 2.00}\\
\end{tabular*}
{\rule{\temptablewidth}{0.75pt}}
\end{center}
\end{table}

\begin{table}[t!]
  \caption{Temporal error of EWI-FP for different $\eps$ at time $t=1$ under $h= 1/8$.}\label{tab:GIFP}
  \vspace*{-10pt}
\begin{center}
\def\temptablewidth{1\textwidth}
{\rule{\temptablewidth}{0.75pt}}
\begin{tabular*}{\temptablewidth}{@{\extracolsep{\fill}}lllllll}
$\fe_{\eps}^{\tau,h}(t=1)$   & $\tau_0=0.2$ & $\tau_0/2^2$ & $\tau_0/2^4$ & $\tau_0/2^6$ & $\tau_0/2^8$ & $\tau_0/2^{10}$ \\[0.25em]
\hline
$\eps_0=1$    &{\bf 1.41E-2}     &8.14E-4	    &5.07E-5	&3.09E-6	 &1.62E-7	    &1.06E-8\\
rate            & --             &2.05             &2.00              &2.02 &2.13 &1.96\\ \hline
$\eps_0/2$      & 1.11E-1	  &{\bf 4.40E-3}	     &2.75E-4	       &1.72E-5	      &1.07E-6	         &6.79E-8\\
rate            & --             &{\bf 2.32}	            & 2.00	           &2.00	          &2.00	             &1.99\\ \hline
$\eps_0/2^2$      & 2.47	     &6.56E-2	   &{\bf 3.90E-3}	   &2.42E-4	   &1.51E-5   	     &9.50E-7\\
rate 	        & --             &2.61	 &{\bf 2.04}	    &2.00	          &2.00    	         &2.00\\ \hline
$\eps_0/2^3$    & 6.73E-1	     &2.82	        & 6.62E-2	&{\bf 4.00E-3} &2.51E-4	 &1.56E-5\\
rate	        & --       &-1.03	 & 2.71	           &{\bf 2.02}  &2.00	             &2.00\\ \hline
$\eps_0/2^4$    &9.50E-1	     &9.28E-1	&2.67 &6.73E-2	&{\bf 4.00E-3}       	 &2.49E-4\\
rate	        & --             & 0.02	  & -0.46     	   & 2.66	  &{\bf 2.04}	             &2.00\\  \hline
$\eps_0/2^5$    &9.96E-1	     &1.05	  &1.11  	&3.87      &6.34E-2 	 &{\bf 3.70E-3}\\
rate	        & --             &-0.04	  &-0.04    	   &-0.90	  &2.97	             &{\bf 2.04}\\
\end{tabular*}
{\rule{\temptablewidth}{0.75pt}}
\end{center}
\end{table}

\begin{table}[t!]
  \caption{Temporal error of TS-FP for different $\eps$ at time $t=1$ under $h = 1/8$.}\label{tab:DIFP}
  \vspace*{-10pt}
\begin{center}
\def\temptablewidth{1\textwidth}
{\rule{\temptablewidth}{0.75pt}}
\begin{tabular*}{\temptablewidth}{@{\extracolsep{\fill}}llllllll}
$\fe_{\eps}^{\tau,h}(t=1)$    &$\tau_0=0.2$ & $\tau_0/2^2$ & $\tau_0/2^4$ & $\tau_0/2^6$ & $\tau_0/2^8$ & $\tau_0/2^{10}$ \\[0.25em]
\hline
$\eps_0=1$     &{\bf 8.49E-3	}	&5.12E-4		&3.19E-5		&2.00E-6		&1.24E-7	&7.64E-9\\
rate		       &---			&2.02		&2.00		&2.00		&2.00	&2.01\\  \hline
$\eps_0/2$  & 8.60E-2	 &{\bf 3.20E-3}	 & 1.97E-4	       & 1.23E-5         &	7.69E-7         & 4.73E-8\\
rate          & --               &{\bf 2.37}             & 2.01	           & 2.00	  & 2.00 & 2.01\\ \hline
$\eps_0/2^2$  & 7.18E-1	     & 2.15E-2       	&{\bf 1.11E-3}	& 6.90E-5	      & 4.31E-6	     & 2.65E-7\\
rate 	      & --               & 2.53	  &{\bf 2.13}	  & 2.00	          & 2.00	         & 2.01\\ \hline
$\eps_0/2^3$  & 6.39E-1	     & 6.39E-1	        & 5.05E-3 &{\bf 2.74E-4}	      & 1.70E-5	     & 1.05E-6\\
rate	      & --               & 0.00	            & 3.49	&{\bf 2.10}	          & 2.01 & 2.01\\ \hline
$\eps_0/2^4$  & 6.84E-1	     & 2.58E-1	        & 2.56E-1	& 1.32E-3 & {\bf 7.18E-5}	     & 4.39E-6\\
rate 	      & --               & 0.70	            & 0.01	  & 3.80	          &{\bf 2.10} & 2.02\\  \hline
$\eps_0/2^5$    &7.64E-1	     &5.03E-2	  &5.77E-2  	&5.88E-2     &3.89E-4  	 &{\bf 2.94E-5}\\
rate	        & --             & 1.96	  &-0.10     	   &-0.01  &3.62	             &{\bf 1.86}\\
\end{tabular*}
{\rule{\temptablewidth}{0.75pt}}
\end{center}
\end{table}

\begin{table}[t!]
  \caption{Temporal error of LI-FP1 for different $\eps$ at time $t=1$ under $h = 1/8$.}\label{tab:APFP1}
  \vspace*{-10pt}
\begin{center}
\def\temptablewidth{1\textwidth}
{\rule{\temptablewidth}{0.75pt}}
\begin{tabular*}{\temptablewidth}{@{\extracolsep{\fill}}lllllll}
$\fe_{\eps}^{\tau,h}(t=1)$    & $\tau_0=0.2$ & $\tau_0/2^2$ & $\tau_0/2^4$ & $\tau_0/2^6$ & $\tau_0/2^8$ & $\tau_0/2^{10}$\\[0.25em]
\hline
$\eps_0=1$      &1.44	 &1.44	&1.44	&1.44	&1.44	&1.44\\
rate            & --         &0.00         &0.00          &0.00         &0.00           &0.00            \\ \hline
$\eps_0/2^2$      &5.05E-1	 &5.05E-1	&5.05E-1	&5.05E-1	&5.05E-1	&5.05E-1\\
rate            & --         &0.40	        &0.00       &0.00     &0.00            &0.00             \\  \hline
$\eps_0/2^4$    &7.11E-2	 &5.37E-2	&5.37E-2	&5.37E-2	&5.37E-2	&5.37E-2\\
rate            & --         &1.68	        &0.03       &0.00    &0.00            &0.00             \\  \hline
$\eps_0/2^6$    &{\bf 4.06E-2}	 &3.40E-3	&3.29E-3	&3.29E-3	&3.29E-3	&3.29E-3	\\
rate            & --         &1.79	        &0.02       &0.00     &0.00            &0.00              \\  \hline
$\eps_0/2^8$    &4.00E-2	 &{\bf 1.07E-3}	&2.12E-4	&2.07E-4	&2.07E-4	&2.07E-4\\
rate	        & --         &{\bf 2.61}	     &1.17      &0.01	   &0.00	          &0.00      \\  \hline
$\eps_0/2^{10}$    &4.19E-2	 &1.11E-3	&{\bf 7.03E-5}	&1.44E-5	&1.33E-5	&1.33E-5\\
rate	        & --         &2.61	     &{\bf 1.99}           &1.14	 &0.06	             &0.00      \\  \hline
$\eps_0/2^{12}$    &4.00E-2	 &1.05E-3	&6.38E-5	&{\bf 4.00E-6}	&8.33-7	&8.16E-7\\
rate	        & --         &2.62	    &2.02           &{\bf 2.00}	            &1.13	     &0.01     \\  \hline
$\eps_0/2^{14}$    &4.13E-2	 &9.43E-4	&5.76E-5	&3.60E-6	&{\bf 2.38E-7}	&7.82E-8\\
rate	        & --         &2.72	    &2.01        &2.00	      &{\bf 1.96}             &0.8  \\
\end{tabular*}
{\rule{\temptablewidth}{0.75pt}}
\end{center}
\end{table}

\begin{table}[t!]
  \caption{Temporal error of LI-FP2 for different $\eps$ at time $t=1$ under $h = 1/8$.}\label{tab:APFP2}
  \vspace*{-10pt}
\begin{center}
\def\temptablewidth{1\textwidth}
{\rule{\temptablewidth}{0.75pt}}
\begin{tabular*}{\temptablewidth}{@{\extracolsep{\fill}}lllllll}
$\fe_{\eps}^{\tau,h}(t=1)$      & $\tau_0=0.2$ & $\tau_0/2^2$ & $\tau_0/2^4$ & $\tau_0/2^6$ & $\tau_0/2^{8}$ & $\tau_0/2^{10}$\\[0.25em]
\hline
$\eps_0=1$    & 12.4	 &13.3	&14.1	 &14.1	&14.1	&14.1\\
rate            & --         &-0.05       &-0.04          &0.00         &0.00           &0.00       \\  \hline
$\eps_0/2^2$    &1.18E-1	 &6.40E-1	& 6.07E-1	 &6.47E-1	&6.50E-1	&6.50E-1\\
rate            & --         &-1.21	    & 0.04	     &-0.04      &-0.01            &0.00              \\  \hline
$\eps_0/2^4$    &\textbf{7.01E-2}	 &6.52E-3	&9.46E-3	 &9.75E-3	&9.77E-3	&9.77E-3\\
rate	        & --         &1.71	    &-0.02	     &-0.01	    & -0.01     &0.00     \\ \hline
$\eps_0/2^{5}$    &4.62E-2	 &\textbf{1.01E-3}	&1.01E-3	 &1.07E-3	&1.07E-3	&1.07E-3\\
rate	        & --         &\textbf{2.75}     	&0.00	     &0.00	    &0.00	    & 0.00   \\  \hline
$\eps_0/2^{6}$    &4.11E-2	 &1.05E-3	&\textbf{4.99E-5}	 &6.11E-5	&6.36E-5	&6.38E-5\\
rate	        & --         &2.64	    &\textbf{2.19	}     &-0.15	     &-0.01     &	-0.01  \\  \hline
$\eps_0/2^{7}$    &4.26E-2	 &1.08E-3	&6.39E-5	 &\textbf{3.15E-6}	&3.63E-6	&3.80E-6\\
rate	        & --         &2.65	    &2.03	     &\textbf{2.17}	    &0.1	    &0.03\\ \hline
$\eps_0/2^{8}$    &4.00E-2	 &1.07E-3	&6.48E-5	 &3.89E-6	&\textbf{1.93E-7}	&2.31E-7\\
rate	        & --         &2.61	    &2.02	     &2.03	    &\textbf{2.16}	    &-0.13\\ \hline
$\eps_0/2^{9}$    &4.14E-2	 &1.13E-3	&6.90E-5	 &4.30E-6	&2.60E-7	&\textbf{1.41E-8}\\
rate	        & --         &2.59	    &2.01	     &2.00	    &2.02	    &\textbf{2.10}\\
\end{tabular*}
{\rule{\temptablewidth}{0.75pt}}
\end{center}
\end{table}

\begin{table}[t!]
  \caption{Temporal error of MTI-FP for different $\eps$ at time $t=1$ under $h = 1/8$.}\label{tab:MTI}
  \vspace*{-10pt}
\begin{center}
\def\temptablewidth{1\textwidth}
{\rule{\temptablewidth}{0.75pt}}
\begin{tabular*}{\temptablewidth}{@{\extracolsep{\fill}}lllllll}
$\fe_{\eps}^{\tau,h}(t=1)$    & $\tau_0=0.2$ & $\tau_0/2^2$ & $\tau_0/2^4$ & $\tau_0/2^6$ & $\tau_0/2^8$ & $\tau_0/2^{10}$ \\[0.25em]
\hline
$\eps_0=1$    &1.90E-1	    &{\bf 1.98E-2}	  &1.49E-3	  &9.73E-5   &6.16E-6	  &3.82E-7\\
rate            & --        &{\bf 1.63}	      &1.87	      &1.97       &1.99	      &2.01\\ \hline
$\eps_0/2$      &1.63E-1	&1.19E-2	  &{\bf 8.26E-4}	  &5.26E-5	  &3.30E-6	  &2.04E-7\\
rate            & --        &1.85	      &{\bf 1.92}	      &1.99	      &2.00	      &2.01\\ \hline
$\eps_0/2^2$    &1.63E-1	&3.22E-2	  &2.62E-3	  &{\bf 1.63E-4}	  &1.01E-5	  &6.28E-7\\
rate 	        & --        &1.17	      &1.81	      &{\bf 2.00}	      &2.00	      &2.01\\ \hline
$\eps_0/2^3$    &{\it 1.01E-1}	&3.68E-2	  &6.22E-3   &5.13E-4	  &{\bf 3.23E-5}	  &2.00E-6\\
rate	        & --        &0.73	      &1.27	      &1.80	      &{\bf 1.99}	      &2.01\\ \hline
$\eps_0/2^4$    &9.67E-2	&1.30E-2	  &9.62E-3	  &1.60E-3	  &1.32E-4	  &{\bf 8.26E-6}\\
rate	        & --        &1.44	      &0.23	      &1.30	      &1.80	      &{\bf 2.00}\\ \hline
$\eps_0/2^5$    &9.50E-2	&{\it 6.22E-3}	  &2.77E-3	  &2.62E-3   &5.03E-4	  &3.86E-5\\
rate	        & --        &{\it 2.00}	      &0.58	      &0.04  	  &1.19	      &1.85\\ \hline
$\eps_0/2^7$    &9.56E-2	&5.61E-3	  &{\it 4.30E-4}	  &1.19E-4	  &1.62E-4	  &1.69E-4\\
rate	        & --        &2.04	      &{\it 1.85}  	  &0.93	      &-0.22	      &-0.03\\ \hline
$\eps_0/2^9$    &9.44E-2	&5.48E-3	  &3.43E-4	  &{\it 2.06E-5}	  &1.19E-6	  &3.51E-6\\
rate	        & --        &2.06	      &2.00  	  &{\it 2.02}	      &2.05	      &-0.77\\ \hline
$\eps_0/2^{11}$ &9.67E-2	&5.60E-3	  &3.48E-4	  &2.19E-5	  &{\it 1.66E-6}	  &1.67E-7\\
rate	        & --        &2.05	      &2.00	      &1.99	      &{\it 1.86}	      &1.66\\  \hline
$\eps_0/2^{13}$    &9.50E-2	&5.48E-3	  &3.40E-4	  &2.12E-5	  &1.29E-6	  &{\it 7.35E-8}\\
rate	        & --        &2.05	      &2.00  	  &2.00	      &2.01	      &{\it 2.06}\\ \hline
$\eps_0/2^{15}$    &9.50E-2	&5.50E-3	  &3.41E-4	  &2.13E-5	  &1.33E-6	  &8.60E-8\\
rate	        & --        &2.05	      &2.00  	  &2.00	      &2.00	      &1.98\\ \hline \hline
$\fe_{\infty}^{\tau,h}$   &1.90E-1	    &3.68E-2	    &9.62E-3	    &2.62E-3	     &5.03E-4	    &1.69E-4	 \\
rate                        &--	        &1.19	        &0.97	        &0.94	         &1.19	        &0.80	      \\
\end{tabular*}
{\rule{\temptablewidth}{0.75pt}}
\end{center}
\end{table}

\begin{table}[t!]
  \caption{Temporal Error of TSF-FP1 for different $\eps$ at time $t=1$ under $h= 1/8,h_\xi=\pi/32$.}\label{tab:TSF1}
  \vspace*{-10pt}
  \begin{center}
\def\temptablewidth{1\textwidth}
{\rule{\temptablewidth}{0.75pt}}
\begin{tabular*}{\temptablewidth}{@{\extracolsep{\fill}}lllllll}
$\fe_{\eps}^{\tau,h}(t=1)$      & $\tau=0.2$ & $\tau/2^2$ & $\tau/2^4$ & $\tau/2^6$ & $\tau/2^8$ & $\tau/2^{10}$ \\[0.25em]
\hline
$\eps_0=1$      &1.07E-1	&3.05E-2	  &7.92E-3	  &2.01E-3	  &5.04E-4	  &1.26E-4\\
rate            & --        &0.91	      &0.97	      &0.99	      &1.00	      &1.00\\  \hline
$\eps_0/2$      &8.88E-2	&4.18E-2	  &1.70E-2	  &5.18E-3	  &1.38E-3	  &3.53E-4\\
rate            & --        &0.54	      &0.65	      &0.86	      &0.95	      &0.99\\ \hline
$\eps_0/2^2$    &6.39E-2	&1.70E-2	  &7.35E-3	  &4.76E-3	  &2.14E-3	  &6.96E-4\\
rate 	        & --        &0.96	      &0.60	      &0.31	      &0.58	      &0.81\\ \hline
$\eps_0/2^3$    &8.43E-2	&1.98E-2	  &5.04E-3	  &1.44E-3	  &6.90E-4	  &4.75E-4\\
rate	        & --        &1.04	      &0.99	      &0.90	      &0.53	      &0.27\\ \hline
$\eps_0/2^4$    &9.67E-2	&2.15E-2	  &5.28E-3	  &1.32E-3	  &3.38E-4	  &9.79E-5\\
rate	        & --        &1.09	      &1.01	      &1.00	      &0.98	      &0.89\\ \hline
$\eps_0/2^5$    &9.05E-2	&1.98E-2	  &4.96E-3	  &1.24E-3	  &3.11E-4	  &7.81E-5\\
rate	        & --        &1.1	      &1.00	      &1.00	      &1.00	      &1.00\\ \hline
$\eps_0/2^8$    &9.61E-2	&2.03E-2	  &5.08E-3	  &1.27E-3	  &3.22E-4	  &8.48E-5\\
rate	        & --        &1.12	      &1.00	      &1.00	      &0.99	      &0.96\\ \hline
$\eps_0/2^{11}$   &1.01E-1	&2.20E-2	  &5.45E-3	  &1.36E-3	  &3.40E-4	  &8.48E-5\\
rate	        & --        &1.10	      &1.01	      &1.00	      &1.00	      &1.00\\ \hline \hline
$\fe_{\infty}^{\tau,h}$   &1.07E-1	    &4.18E-2    &1.70E-2	    &5.18E-3     &2.14E-3	    &6.96E-4	 \\
rate                        &--	        &0.68	        &0.65	        &0.86	         &0.64	        &0.81	      \\
\end{tabular*}
{\rule{\temptablewidth}{0.75pt}}
\end{center}
\end{table}

\begin{table}[t!]
  \caption{Temporal error of TSF-FP2 for different $\eps$ at time $t=1$ under $h= 1/8,h_\xi=\pi/32$.}\label{tab:TSF2}
  \vspace*{-10pt}
\begin{center}
\def\temptablewidth{1\textwidth}
{\rule{\temptablewidth}{0.75pt}}
\begin{tabular*}{\temptablewidth}{@{\extracolsep{\fill}}lllllll}
  $\fe_{\eps}^{\tau,h}(t=1)$   & $\tau=0.2$ & $\tau/2^2$ & $\tau/2^4$ & $\tau/2^6$ & $\tau/2^8$ & $\tau/2^{10}$ \\[0.25em]
\hline
$\eps_0=1$      &1.86E-2	&1.18E-3	  &7.35E-5   &4.57E-6	  &2.84E-7	  &1.67E-8\\
rate            & --        &1.99	      &2.00	      &2.00	      &2.00	      &2.04\\ \hline
$\eps_0/2$      &3.45E-2	&5.25E-3	  &3.44E-4	  &2.15E-5	  &1.35E-6	  &8.26E-8\\
rate            & --        &1.36	      &1.97	      &2.00	      &2.00	      &2.01\\ \hline
$\eps_0/2^2$    &2.94E-2	&2.82E-3	  &9.16E-4	  &8.77E-5   &5.47E-6	  &3.39E-7\\
rate 	        & --        &1.69	      &0.81	      &1.69	      &2.00 	  &2.01\\ \hline
$\eps_0/2^3$    &2.43E-2	&1.01E-3	  &2.00E-4	  &7.07E-5	  &1.15E-5   &7.35E-7\\
rate	        & --        &2.29	      &1.17	      &0.75	      &1.31	      &1.98\\ \hline
$\eps_0/2^4$    &3.34E-2	&5.94E-4	  &7.47E-5	  &1.39E-5	  &1.28E-6	  &5.94E-7\\
rate	        & --        &2.90	      &1.50  	  &1.21	      &1.72	      &0.55\\ \hline
$\eps_0/2^5$    &3.73E-2	&5.43E-4	  &3.54E-5	  &5.12E-6	  &1.00E-6	  &8.14E-8\\
rate	        & --        &3.05	      &1.97	      &1.40  	  &1.18	      &1.81\\ \hline
$\eps_0/2^8$    &3.85E-2	&5.25E-4	  &3.11E-5	  &1.94E-6	  &1.21E-7	  &7.81E-9\\
rate	        & --        &3.09	      &2.04	      &2.00	      &2.00	      &1.98\\ \hline
$\eps_0/2^{11}$   &3.79E-2	&5.25E-4	  &2.82E-5	  &1.76E-6	  &1.10E-7	  &6.39E-9\\
rate	        & --        &3.15	      &2.04	      &2.00	      &2.00	      &2.04\\ \hline \hline
$\fe_{\infty}^{\tau,h}$   &3.85E-2	    &5.25E-3	    &9.16E-4	    &8.77E-5	     &1.15E-5	    &7.35E-7 \\
rate                        &--	        &1.44	        &1.26	        &1.70	         &1.47	        &1.98	      \\
\end{tabular*}
{\rule{\temptablewidth}{0.75pt}}
\end{center}
\end{table}

\begin{table}[t!]
  \caption{Temporal error of IEI-FP for different $\eps$ at time $t=1$ under $h= 1/8$.}\label{tab:IEI}
  \vspace*{-10pt}
\begin{center}
\def\temptablewidth{1\textwidth}
{\rule{\temptablewidth}{0.75pt}}
\begin{tabular*}{\temptablewidth}{@{\extracolsep{\fill}}lllllll}
  $\fe_{\eps}^{\tau,h}(t=1)$   & $\tau=0.2$ & $\tau/2^2$ & $\tau/2^4$ & $\tau/2^6$ & $\tau/2^8$ & $\tau/2^{10}$ \\[0.25em]
\hline
$\eps_0=1$      &5.43E-2    &3.58E-3      &2.45E-4    &1.57E-5    &9.84E-7    &6.11E-8\\
rate            & --        &1.96         &1.94       &1.98       &2.00       &2.00\\ \hline
$\eps_0/2$      &2.43E-2    &2.16E-3      &1.40E-4    &8.77E-6    &5.48E-7    &3.43E-8\\
rate            & --        &1.75         &1.97       &2.00       &2.00       &2.00\\ \hline
$\eps_0/2^2$    &1.19E-1    &2.36E-3      &1.36E-4    &8.43E-6    &5.27E-7    &3.26E-8\\
rate 	        & --        &2.83         &2.06       &2.00       &2.00       &2.01\\ \hline
$\eps_0/2^3$    &5.71E-2    &1.70E-2      &8.48E-5    &4.75E-6    &2.91E-7    &1.35E-8\\
rate	        & --        &0.88         &3.82       &2.08       &2.01       &2.21\\ \hline
$\eps_0/2^4$    &3.62E-2    &5.31E-3      &1.47E-3    &4.61E-6    &3.43E-7    &1.67E-8\\
rate	        & --        &1.39         &0.92       &4.16       &1.88       &2.18\\ \hline
$\eps_0/2^5$    &3.68E-2    &6.73E-4      &6.11E-5    &1.51E-5    &3.26E-7    &2.12E-8\\
rate	        & --        &2.89         &1.73       &1.01       &2.77       &1.97\\ \hline
$\eps_0/2^8$    &3.85E-2    &7.07E-4      &4.19E-5    &2.58E-6    &1.57E-7    &7.81E-9\\
rate	        & --        &2.88         &2.04       &2.01       &2.02       &2.17\\ \hline
$\eps_0/2^{11}$   &3.85E-2  &6.96E-4      &4.21E-5    &2.62E-6    &1.62E-7    &7.58E-9 \\
rate	        & --        &2.89	      &2.03	      &2.00	      &2.01	      &2.21\\ \hline \hline
$\fe_{\infty}^{\tau,h}$   &1.19E-1	    &1.70E-2	    &1.47E-3	    &1.57E-5	     &9.84E-7	    &6.11E-8 \\
rate                        &--	        &1.41	        &1.76	        &3.27	         &2.00	        &2.00	      \\
\end{tabular*}
{\rule{\temptablewidth}{0.75pt}}
\end{center}
\end{table}

\begin{table}[t!]
  \caption{Temporal error of MTI-FP2 for different $\eps$ at time $t=1$ under $h= 1/8$.}\label{tab:MTI2}
  \vspace*{-10pt}
\begin{center}
\def\temptablewidth{1\textwidth}
{\rule{\temptablewidth}{0.75pt}}
\begin{tabular*}{\temptablewidth}{@{\extracolsep{\fill}}lllllll}
  $\fe_{\eps}^{\tau,h}(t=1)$   & $\tau=0.2$ & $\tau/2^2$ & $\tau/2^4$ & $\tau/2^6$ & $\tau/2^8$ & $\tau/2^{10}$ \\[0.25em]
\hline
$\eps_0=1$      &5.65E-2    &3.91E-3      & 2.47E-4   &1.54E-5    &9.60E-7    &5.44E-8\\
rate            & --        &1.93         &1.99       &2.00       &2.00       &2.07\\ \hline
$\eps_0/2$      &9.35E-2    &8.88E-3      &5.40E-4    &3.34E-5    &2.08E-6    &1.31E-7\\
rate            & --        &1.70         &2.02       &2.00       &2.00       &1.99\\ \hline
$\eps_0/2^2$    &1.33E-1    &2.13E-2      &1.15E-3   &7.02E-5    &4.34E-6    &2.68E-7\\
rate 	        & --        &1.32         &2.11       &2.02       &2.01       &2.01\\ \hline
$\eps_0/2^3$    &2.10E-1    &1.35E-2      &2.00E-3    &9.72E-5    &5.83E-6    &3.59E-7\\
rate	        & --        &1.98         &1.38       &2.18       &2.03       &2.01\\ \hline
$\eps_0/2^4$    &2.45E-1    &1.55E-2      &9.77E-4    &1.38E-4    &6.66E-6    &3.97E-7\\
rate	        & --        &2.03         &1.99       &1.41       &2.18       &2.04 \\ \hline
$\eps_0/2^5$    &2.62E-1    &1.59E-2      &9.97E-4    &6.23E-5    &8.88E-6    &4.31E-7\\
rate	        & --        &2.02         &1.99       &1.98       &1.41       &2.18\\ \hline
$\eps_0/2^8$    &2.64E-1    &1.62E-2      &1.00E-3    &6.28E-5    &3.94E-6    &2.48E-7 \\
rate	        & --        &2.00         &2.01       &2.00       &2.00       &2.00 \\ \hline
$\eps_0/2^{11}$   &2.58E-1  &1.64E-2      &1.01E-3    &6.33E-5    &3.94E-6    &2.45E-7 \\
rate	        & --        &1.99         &2.01       &2.00       &2.00       &2.00\\ \hline \hline
$\fe_{\infty}^{\tau,h}$   &2.64E-1	    &2.13E-2 	    &2.00E-3 	    &1.38E-4	     &8.88E-6 	    &4.31E-7 \\
rate                        &--	        &1.82	        &1.71	        &1.92	         &1.98	        &2.18	      \\
\end{tabular*}
{\rule{\temptablewidth}{0.75pt}}
\end{center}
\end{table}

As a summary, a detailed table on the comparison of computational complexity of each method has been given in Tab. \ref{tab:complex}. The comparison of the temporal error of each method in the classical regime, i.e. $\eps=O(1)$, has been given in Tab. \ref{tab:comp1} together with the computational time. The comparisons of the temporal error in the limit regime are given in Tab. \ref{tab:comp2}. Tab. \ref{tab:comp3} shows the temporal error of different methods under the natural mesh strategy, i.e.  $\tau=O(\eps^2)$ which fully resolves the temporal wavelength of the oscillation. All methods are programmed with Matlab and run on an Intel i3-3120M 2.5GHz CPU laptop.

\subsection{Comparison of different methods}

\begin{table}[t!]
  \caption{Comparison of properties of different numerical methods. Here $N$ denotes the number of grid points in $x$-direction and $N_\xi$ denotes the number of grid point in $\xi$-direction.}\label{tab:complex}
  \vspace*{-10pt}
\begin{center}
\def\temptablewidth{1\textwidth}
{\rule{\temptablewidth}{0.75pt}}
\begin{tabular*}{\temptablewidth}{@{\extracolsep{\fill}}ccccccc}
Method          &LFFD       &SIFD         &ECFD
&EWI-FP &TS-FP
 &$\ba{c} \hbox{LI-FP1}\\
\hbox{(or\ LI-FP2)}\\
\ea$     \\[0.25em]
\hline
Time symmetric  &Yes    	&Yes	      &Yes	      &Yes	  &Yes    &Yes	  \\
Energy conservation   &No	&No	          &Yes        &No	  &No    &No\\
Unconditionally stable   &No	    &No	          &No	      &Yes	   &Yes   &Yes\\
Explicit        &Yes    	&No	          &No	      &Yes	  &Yes    &Yes\\
Temporal accuracy& 2nd      &2nd          &2nd        &2nd    &2nd    &2nd\\
Spatial accuracy& 2nd       &2nd          &2nd        &spectral  &spectral  &spectral\\
Memory cost          &$O(N)$	    &$O(N)$	      &$O(N)$	      &$O(N)$	 &$O(N)$     &$O(N)$\\
Computational cost    &$O(N)$	&$O(N)$	      &$\gg O(N)$	  &$O(N\ln N)$	&$O(N\ln N)$  &$O(N\ln N)$\\ \hline
Resolution      &$h$=O(1) &$h$=O(1)	&$h$=O(1)&$h$=O(1) &$h$=O(1) &$h$=O(1) \\
when $0<\eps\ll1$ &$\tau=O(\eps^3)$   &$\tau=O(\eps^3)$	&$\tau=O(\eps^3)$ &$\tau=O(\eps^2)$	
&$\tau=O(1)$ &$\tau=O(1)$\\
Uniformly accurate  &No      &No           &No          &No   &No    &No\\
\end{tabular*}
{\rule{\temptablewidth}{0.75pt}}
\end{center}
\end{table}
\addtocounter{table}{-1}

\begin{table}[t!]
\renewcommand\thetable{18}
  \caption{(con't)}
  \vspace*{-10pt}
\begin{center}
\def\temptablewidth{1\textwidth}
{\rule{\temptablewidth}{0.75pt}}
\begin{tabular*}{\temptablewidth}{@{\extracolsep{\fill}}cccccc}
Method             &MTI-FP    &TSF-FP1  &TSF-FP2  &IEI-FP   &MTI-FP2\\[0.25em]
\hline
Time symmetric  	     	      &No	      &No &No	&  No&  No\\
Energy conservation                   &No	      &No &No &   No&  No\\
Unconditionally stable     	     	      &Yes	      &Yes &Yes &Yes&Yes\\
Explicit        	         	      &Yes	      &Yes &Yes & Yes&Yes\\
Temporal accuracy                 &2nd        &1st &2nd &2nd&2nd\\
Spatial accuracy                &spectral    &spectral&spectral &spectral&spectral\\
Memory cost             	      &$O(N)$	      &$O(N_\xi N)$ &$O(N_\xi N)$ &$O(N)$&$O(N)$\\
Computational cost       	  &$O(N\ln N)$	  &$O(N_\xi N\ ln N)$ &$O(N_\xi N\ ln N)$ &$O(N\ln N)$&$O(N\ln N)$\\ \hline
Resolution      	&$h$=O(1) &$h,h_\xi$=O(1)  &$h,h_\xi$=O(1) &$h$=O(1) &$h$=O(1)\\
when $0<\eps\ll1$ 	&$\tau=O(1)$	&$\tau=O(1)$ &$\tau=O(1)$ &$\tau=O(1)$&$\tau=O(1)$\\
Uniformly accurate                  &Yes        &Yes &Yes & Yes& Yes\\
Optimally accurate                  &No         &No  &No  &Yes &Yes
\end{tabular*}
{\rule{\temptablewidth}{0.75pt}}
\end{center}
\end{table}

\begin{table}[t!]
  \caption{Comparison of temporal errors and their corresponding
computational time (seconds) of
different methods for the NKGE (\ref{KG}) with $\eps = 1$, $h=1/8$ and $h_\xi=\pi/16$.}\label{tab:comp1}
  \vspace*{-10pt}
\begin{center}
\def\temptablewidth{1\textwidth}
{\rule{\temptablewidth}{0.75pt}}
\begin{tabular*}{\temptablewidth}{@{\extracolsep{\fill}}lllllll}
  $\fe_{\eps}^{\tau,h}(t=1)$   & $\tau_0=0.2$ & $\tau_0/2^2$ & $\tau_0/2^4$ & $\tau_0/2^6$ & $\tau_0/2^8$ & $\tau_0/2^{10}$ \\[0.25em]
\hline
 EWI-FP  &1.41E-2     &8.14E-4	    &5.07E-5	&3.09E-6	 &1.62E-7	    &1.06E-8 \\
 {\rm time (cpu)} &4.5E-4 	 	&1.7E-3 		&6.6E-3 	&2.6E-2 		&1E-1 	&4.2E-1 \\
 \hline
 TS-FP &8.49E-3		&5.12E-4		&3.19E-5		&2.00E-6		&1.24E-7	&7.64E-9   \\
 {\rm time (cpu)} &7E-4 	&2.6E-3 		&1E-2 	&3.7E-2 	 	& 1.5E-1		&5.9E-1   \\
 \hline
 LI-FP2  &12.4	 &13.3	&14.1	 &14.1	&14.1	&14.1 \\
 {\rm time (cpu)} &2.1E-3 	 	&8.3E-3 		&3.3E-2 		&1.3E-1 		&5E-1 		&2.1 \\
 \hline
  MTI-FP &1.90E-1	    &1.98E-2	  &1.49E-3	  &9.73E-5   &6.16E-6	  &3.82E-7   \\
 {\rm time (cpu)} &2.5E-3 	&8.5E-3 		&3.4E-2 	&1.4E-1 	 	&5.5E-1 		&2.2   \\
 \hline
 TSF-FP2  &1.86E-2	&1.18E-3	  &7.35E-5   &4.57E-6	  &2.84E-7	  &1.67E-8  \\
 {\rm time (cpu)} &4E-2 	 	&1.6E-1 		&6.3E-1 		&2.5 		&9.9 		&39.1  \\
 \hline
 IEI-FP  &5.43E-2    &3.58E-3      &2.45E-4    &1.57E-5    &9.84E-7    &6.11E-8 \\
 {\rm time (cpu)} &5.6E-3 	 	&1.6E-2 		&9.4E-2 		&3.4E-1 		&1.4 		&5.6  \\
 \hline
 MTI-FP2  &5.65E-2    &3.91E-3      & 2.47E-4   &1.54E-5    &9.60E-7    &5.44E-8  \\
 {\rm time (cpu)} &1.6E-2 	 	&4.7E-2 		&1.9E-1 		&7.3E-1 		&3.0 		&11.6  \\
\end{tabular*}
{\rule{\temptablewidth}{1pt}}
\end{center}
\end{table}

\begin{table}[t!]
  \caption{Comparison of temporal errors and their corresponding
computational time (seconds) of
different methods for the NKGE (\ref{KG}) with $\eps = 2^{-11}$, $h=1/8$ and $h_\xi=\pi/4$.}\label{tab:comp2}
  \vspace*{-10pt}
\begin{center}
\def\temptablewidth{1\textwidth}
{\rule{\temptablewidth}{0.75pt}}
\begin{tabular*}{\temptablewidth}{@{\extracolsep{\fill}}lllllll}
  $\fe_{\eps}^{\tau,h}(t=1)$   & $\tau_0=0.2$ & $\tau_0/2^2$ & $\tau_0/2^4$ & $\tau_0/2^6$ & $\tau_0/2^8$ & $\tau_0/2^{10}$ \\[0.25em]
\hline
 EWI-FP  &9.97E-1 	 	&1.12 		&1.18 	&1.20 		&1.22 	&1.21 \\
 {\rm time (cpu)}& 4.5E-4 	 	&1.7E-3 		&6.6E-3 	&2.6E-2 		&1E-1 	&4.2E-1 \\
 \hline
 TS-FP &6.39E-1  &3.07E-1  &8.48E-3   &5.65E-3   &6.22E-3  &2.03E-4   \\
 {\rm time (cpu)} &7E-4 	&2.6E-3 		&1E-2 	&3.7E-2 	 	& 1.5E-1		&5.9E-1  \\
 \hline
 LI-FP2  &4.15E-2 &9.45E-4 &5.77E-5 &3.60E-6 &2.27E-7 &1.63E-8  \\
 {\rm time (cpu)} &2.1E-3 	 	&8.3E-3 		&3.3E-2 		&1.3E-1 		&5E-1 		&2.1 \\
 \hline
  MTI-FP &9.67E-2	&5.60E-3	  &3.48E-4	  &2.19E-5	  &1.66E-6	  &1.67E-7   \\
 {\rm time (cpu)} &2.5E-3 	&8.5E-2 		&3.4E-2 	&1.4E-1 	 	&5.5E-1 		&2.2 \\
 \hline
 TSF-FP2  &3.79E-2	&5.25E-4	  &2.82E-5	  &1.76E-6	  &1.10E-7	  &6.39E-9 \\
 {\rm time (cpu)} &9E-3 	 	& 3E-2		&1.2E-1 		&4.8E-1 		&1.8 		&7.4  \\
 \hline
 IEI-FP  &3.85E-2  &6.96E-4      &4.21E-5    &2.62E-6    &1.62E-7    &7.58E-9 \\
 {\rm time (cpu)} &5.6E-3 	 	&1.6E-2 		&9.4E-2 		&3.4E-1 		&1.4 		&5.6  \\
 \hline
 MTI-FP2  &2.58E-1  &1.64E-2      &1.01E-3    &6.33E-5    &3.94E-6    &2.45E-7 \\
 {\rm time (cpu)} &1.6E-2 	 	&4.7E-2 		&1.9E-1 		&7.3E-1 		&3.0 		&11.6  \\
\end{tabular*}
{\rule{\temptablewidth}{1pt}}
\end{center}
\end{table}

\begin{table}[t!]
  \caption{Comparison of temporal errors and their corresponding
computational time (seconds) of
different methods for the NKGE (\ref{KG}) with $\tau=O(\eps^2)$, $h=1/8$ and $h_\xi=\pi/16$.}\label{tab:comp3}
  \vspace*{-10pt}
\begin{center}
\def\temptablewidth{1\textwidth}
{\rule{\temptablewidth}{0.75pt}}
\begin{tabular*}{\temptablewidth}{@{\extracolsep{\fill}}lllllll}
  $\fe_{\eps}^{\tau,h}(t=1)$
  &$\begin{array}{l} \eps_0=1.0\\
  \tau_0=0.2\\
   \end{array}$
  &$\begin{array}{l} \eps_0/2\\
  \tau_0/2^2\\
   \end{array}$
  &$\begin{array}{l} \eps_0/2^2\\
  \tau_0/2^4\\
   \end{array}$
   &$\begin{array}{l} \eps_0/2^3\\
  \tau_0/2^6\\
   \end{array}$
   &$\begin{array}{l} \eps_0/2^4\\
  \tau_0/2^8\\
   \end{array}$  &$\begin{array}{l} \eps_0/2^5\\
  \tau_0/2^{10}\\
   \end{array}$  \\[0.25em]
\hline
 EWI-FP  &1.41E-2 &4.42E-3 &3.88E-3 &4.01E-3  &3.99E-3 &3.74E-3  \\
 {\rm time (cpu)} &4.5E-4 	 	&1.7E-3 		&6.6E-3 	&2.6E-2 		&1E-1 	&4.2E-1  \\
 \hline
 TS-FP &8.48E-3  &3.20E-3  &1.11E-3   &2.74E-4   &7.18E-5  &2.94E-5   \\
 {\rm time (cpu)} &7E-4 	&2.6E-3 		&1E-2 	&3.7E-2 	 	& 1.5E-1		&5.9E-1 \\
 \hline
 LI-FP2  &12.4 &3.16 &6.47E-1 &1.18E-1 &9.77E-3 &1.07E-3  \\
 {\rm time (cpu)} &2.1E-3 	 	&8.3E-3 		&3.3E-2 		&1.3E-1 		&5E-1 		&2.1 \\
 \hline
  MTI-FP &1.90E-1  &1.19E-2  &2.62E-3   &5.12E-4   &1.32E-4  &3.86E-5   \\
 {\rm time (cpu)} &2.5E-3 	&8.5E-2 		&3.4E-2 	&1.4E-1 	 	&5.5E-1 		&2.2 \\
 \hline
 TSF-FP2  &1.87E-2 &5.25E-3 &9.16E-4 &7.07E-5 &1.28E-6 &8.14E-8  \\
 {\rm time (cpu)} &4E-2 	 	&1.6E-1 		&6.3E-1 		&2.5 		&9.9		&39.1  \\ \hline

  IEI-FP  &5.41E-2 &2.16E-3 &1.36E-4 &4.75E-6 &3.44E-7 &2.14E-8  \\
 {\rm time (cpu)} &5.6E-3 	 	&2.2E-2 		&9.5E-2 		&3.6E-1 		&1.4		&5.5  \\ \hline
  MTI-FP2  &5.64E-2 &8.87E-3 &1.15E-3 &9.72E-5 &6.66E-6 &4.31E-7  \\
 {\rm time (cpu)} &1.5E-2 	 	&4.7E-2 		&1.9E-1 		&7.2E-1 		&3		&11.5  \\
\end{tabular*}
{\rule{\temptablewidth}{1pt}}
\end{center}
\end{table}

From Tabs. \ref{tab:ECFD}-\ref{tab:comp3}, we can draw the following conclusions:

(i) All FDTD methods have temporal $\eps$-scalability $\tau=O(\eps^3)$ (cf. Tabs. \ref{tab:ECFD}-\ref{tab:LFFD}). In the classical regime, LFFD is the most accurate and efficient method among the three. However, all FDTD methods become inefficient when $\eps$ becomes small.

(ii) Both EWI-FP and TS-FP have temporal $\eps$-scalability $\tau=O(\eps^2)$ (cf. Tabs. \ref{tab:GIFP}\&\ref{tab:DIFP}). While when $\tau\lesssim\eps^2$, TS-FP has an improved error bound at $\tau^2/\eps^2$ with respect to the small parameter $\eps\in(0,1]$. Both methods perform very well in the classical regime (cf. Tab. \ref{tab:comp1}), but they are unsatisfactory in the limit regime. They have very similar efficiency but TS-FP is more accurate when $\eps$ is small.

(iii) The two LI methods are accurate in the limit regime, but both of them do not have convergence in the classical or the intermediate regime (cf. Tabs. \ref{tab:APFP1}\&\ref{tab:APFP2}). LI-FP2 is more accurate than LI-FP1 due to the correction, but LI-FP2 is much less efficient.

(iv) MTI-FP, TSF-FP1, TSF-FP2, IEI-FP and MTI-FP2 have temporal $\eps$-scalability $\tau=O(1)$ and they offer uniformly correct results for all $\eps\in(0,1]$ (cf. Tabs. \ref{tab:MTI}-\ref{tab:TSF2}). All the five UA   methods have some temporal convergence order reductions in the resonance regime.  Between the two first order UA schemes, MTI-FP is more accurate than TSF-FP1 (cf. Tabs. \ref{tab:MTI}\&\ref{tab:TSF1}). TSF-FP2 is the most accurate method among the five (cf. Tabs. \ref{tab:comp1}\&\ref{tab:comp2}), while from the computational cost point of view, TSF-FP2 is more expensive than the UOA methods due to the extra dimension (cf. Tabs. \ref{tab:complex}\&\ref{tab:comp1}), especially the memory cost in high dimensions (cf. Tab. \ref{tab:complex}), and IEI-FP is found to be most efficient (cf. Tabs. \ref{tab:comp1}\&\ref{tab:comp2}\&\ref{tab:comp3}).

(v) Among all the methods, in the $\eps=O(1)$ regime, the EWI-FP method and TS-FP are the most accurate and efficient methods (cf. Tab. \ref{tab:comp1}). While in the intermediate regime and the limit regime, the two UOA methods, i.e. IEI-FP and MTI-FP2 are significantly more powerful than the others (cf. Tabs. \ref{tab:comp2}\&\ref{tab:comp3}).

\begin{remark}
The second order uniform accuracy of TSF-FP2 has been shown in \cite{Chartier} under condition $\tau\leq C$ with $C>0$ independent of $\eps$. Here in our test (cf. Tab. \ref{tab:TSF2}), we are interested in performance of the scheme with a wide range of time step. Hence the order reduction here does not conflict with the theoretical results.
\end{remark}

\section{Applications} \label{appl9}

In this section, we apply the UOA MTI-FP2 method to study numerically the convergence rates from the NKGE (\ref{KG}) to its limiting  models (\ref{nlsw}) and (\ref{nls}), and to simulate wave interaction in two dimensions (2D).

\subsection{Convergence rates of NKGE to its limiting models}

We take $d=1$ and $\lambda=1$ in the NKGE (\ref{KG}).
Let $u$ be the solution of NKGE (\ref{KG}), $z_{\rm sw}$ be the solution of
the NLSW (\ref{nlsw}) and $z_{\rm s}$ be the solution of
the NLSE (\ref{nls}). Take the initial data as
\begin{equation}\label{smooth}
 \phi_1(x)=\frac{\fe^{-x^2}}{\sqrt{\pi}},\quad \phi_2(x)=\frac{1}{2}\sech(x^2)\sin(x), \qquad x\in{\mathbb R},
\end{equation}
or
\begin{equation}\label{nonsmooth}
 \phi_1(x)=\frac{x^m|x|\fe^{-x^2}}{\sqrt{\pi}},\quad \phi_2(x)=\frac{1}{2}\sech(x^2)\sin(x),\qquad x\in{\mathbb R},
\end{equation}
where $m=1,2$.
The solutions are obtained numerically with very fine mesh on a bounded interval $\Omega=(-128,128)$ with periodic boundary conditions.
Define
\[
u_{\rm sw}(x,t):=\fe^{it/\eps^2}z_{\rm sw}(x,t)+\fe^{-it/\eps^2}\overline{z_{\rm sw}}(x,t),\qquad
u_{\rm s}(x,t):=\fe^{it/\eps^2}z_{\rm s}(x,t)+\fe^{-it/\eps^2}\overline{z_{\rm s}}(x,t),
\]
and define the error functions as
\begin{equation}\label{eta def}
\begin{split}
&\eta_{\rm sw}(t):=\|u(\cdot,t)-u_{\rm sw}(\cdot,t)\|_{H^1},\qquad
\eta_{\rm s}(t):=\|u(\cdot,t)-u_{\rm s}(\cdot,t)\|_{H^1}.
\end{split}
\end{equation}
Fig. \ref{fig:convergence} shows the errors defined in (\ref{eta def}) as  functions of time with the smooth initial data (\ref{smooth}). Fig. \ref{fig:convergence2} and Fig. \ref{fig:convergence3} show the results from the nonsmooth initial data (\ref{nonsmooth}) with $m=2$ and $m=1$, respectively.

\begin{figure}[t!]
\centerline{\psfig{figure=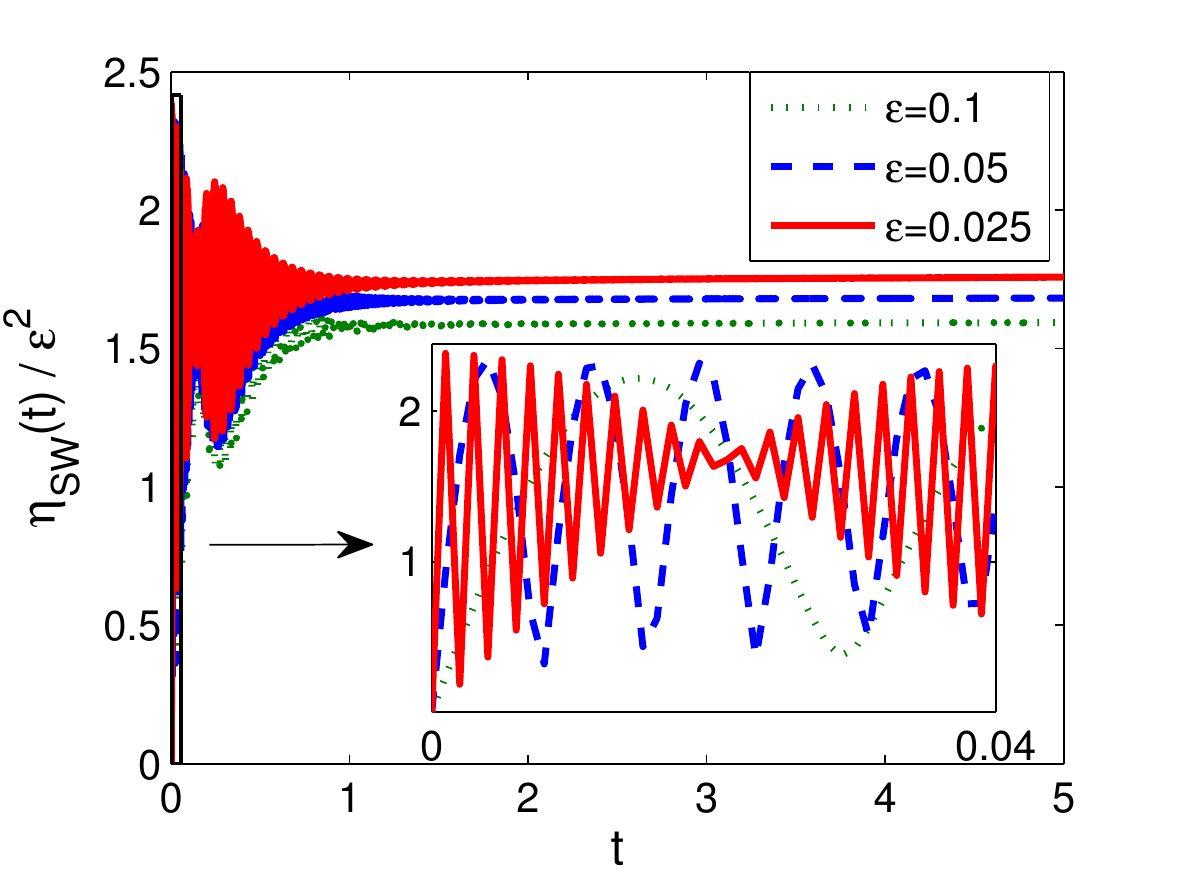,height=6cm,width=8cm}
\psfig{figure=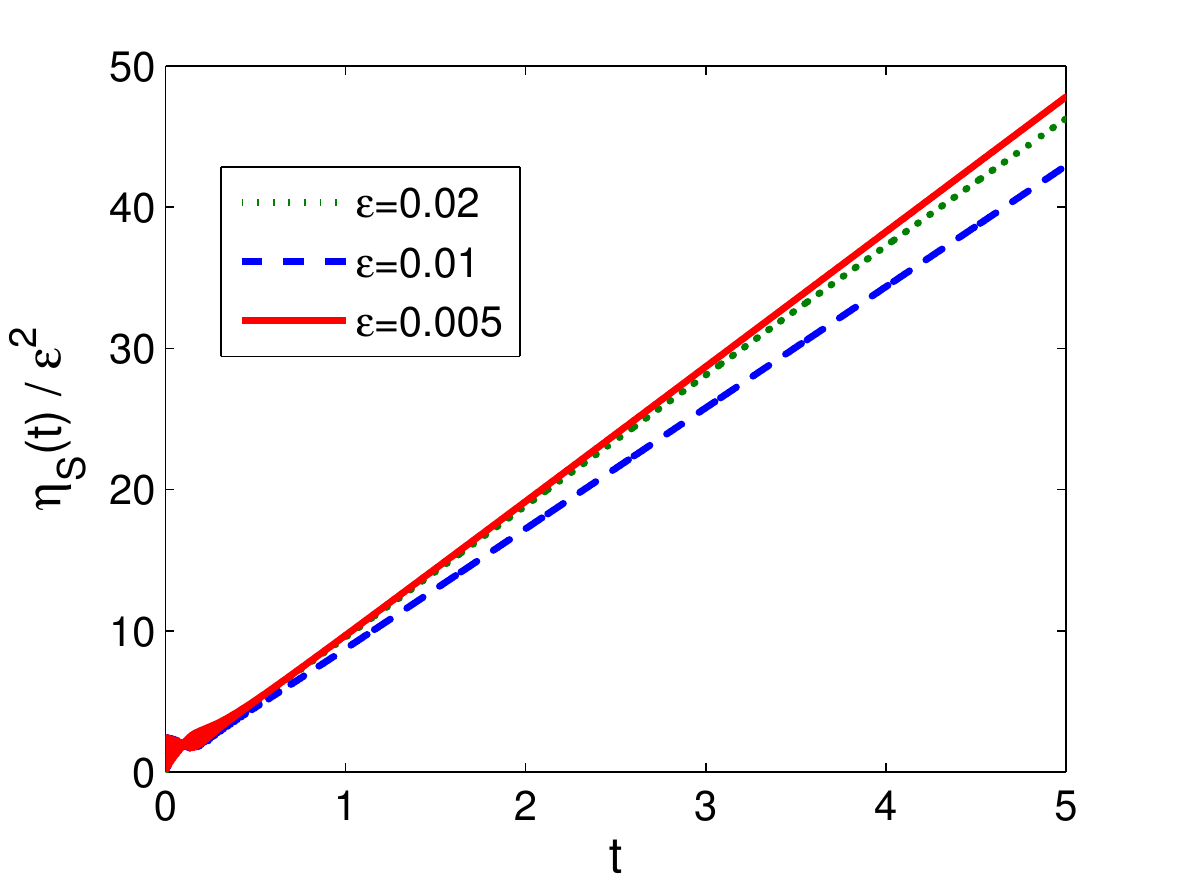,height=6cm,width=8cm}}
\caption{Time evolution of  $\eta_{\rm sw}(t)$ and $\eta_{\rm s}(t)$ for the smooth initial data (\ref{smooth})
under different $\eps$.}\label{fig:convergence}
\end{figure}

\begin{figure}[t!]
\centerline{\psfig{figure=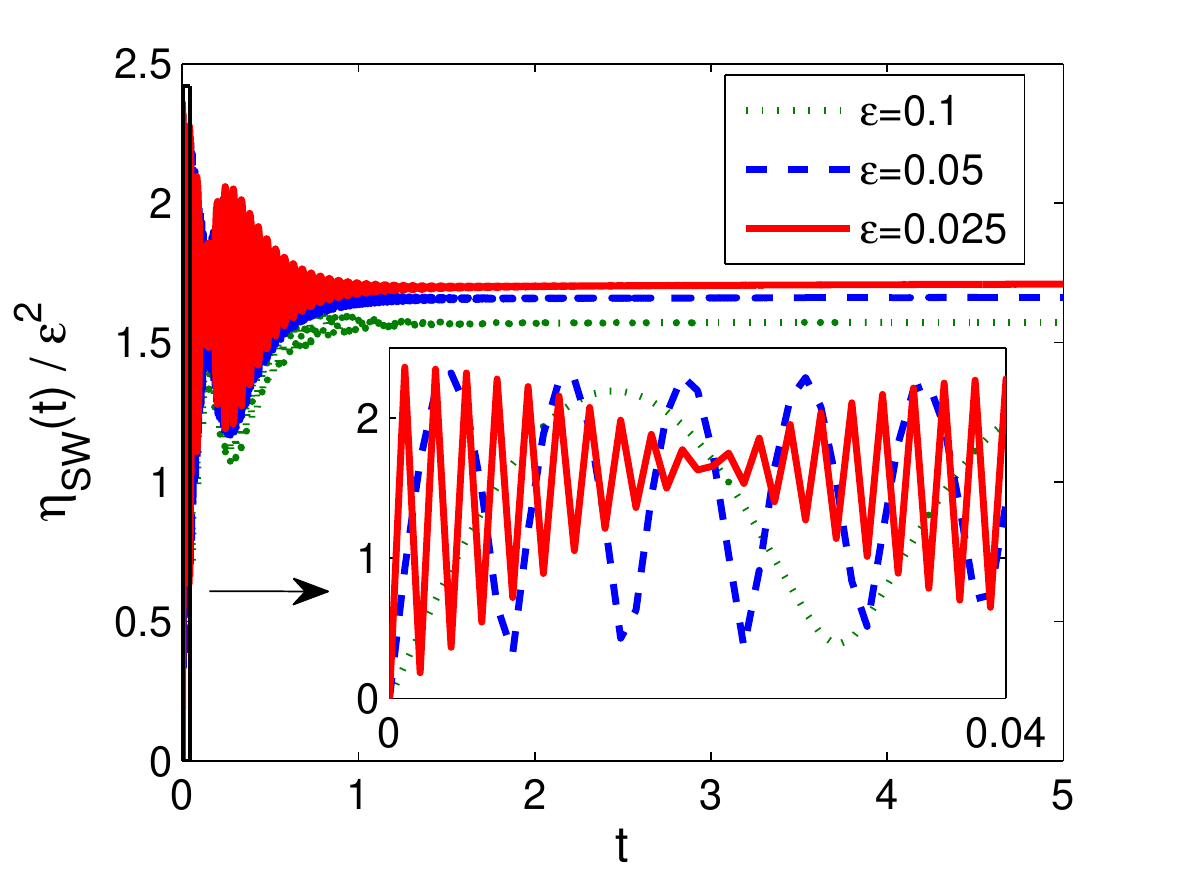,height=6cm,width=8cm}
\psfig{figure=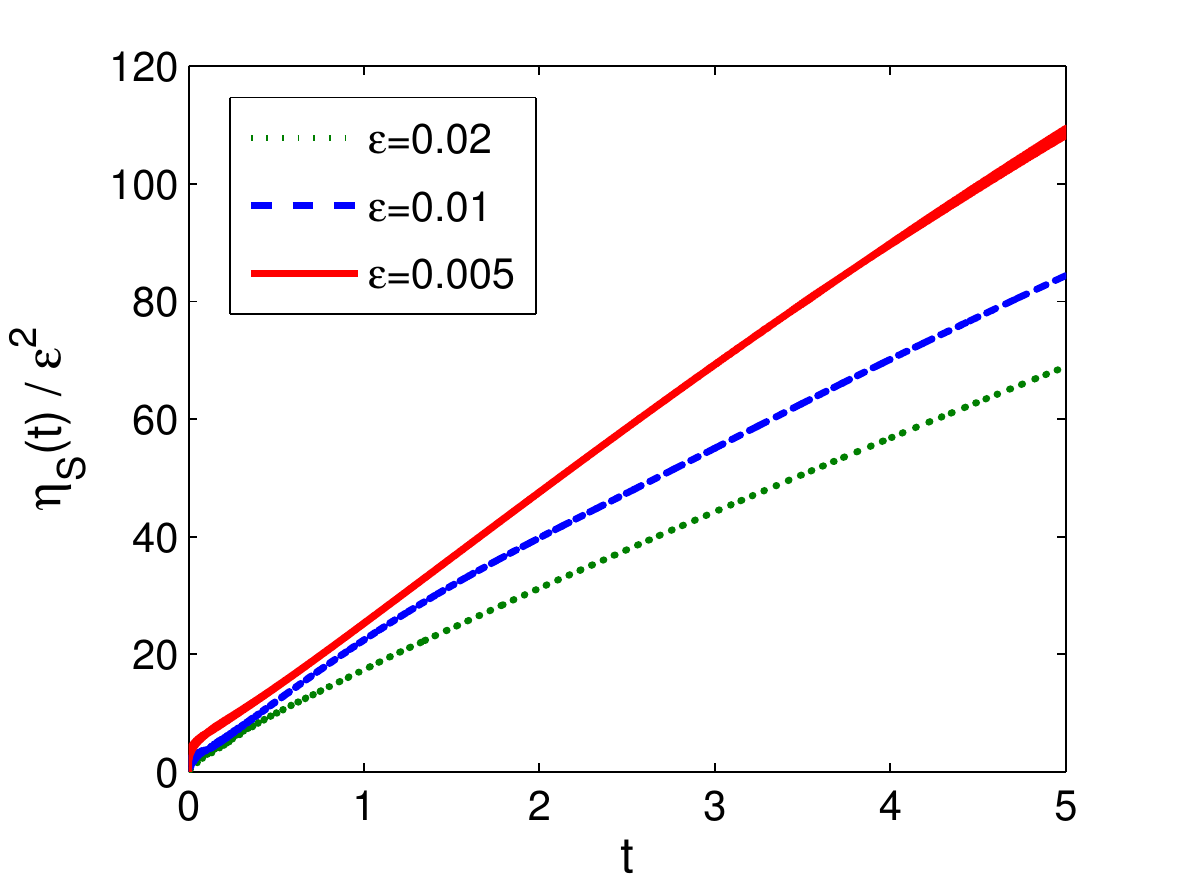,height=6cm,width=8cm}}
\caption{Time evolution of $\eta_{\rm sw}(t)$ and $\eta_{\rm s}(t)$ with nonsmooth data (\ref{nonsmooth}) for $m=2$ under different $\eps$.}\label{fig:convergence2}
\end{figure}

\begin{figure}[t!]
\centerline{\psfig{figure=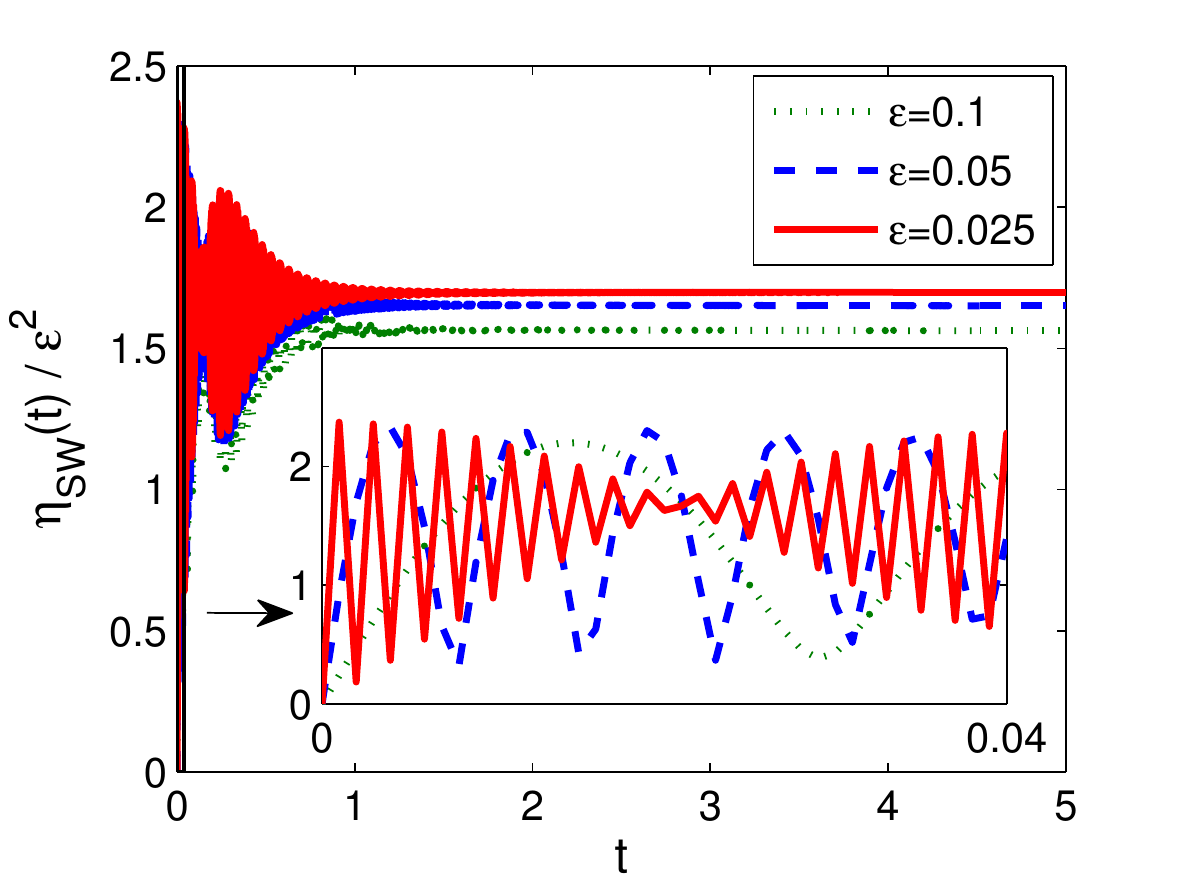,height=6cm,width=8cm}
\psfig{figure=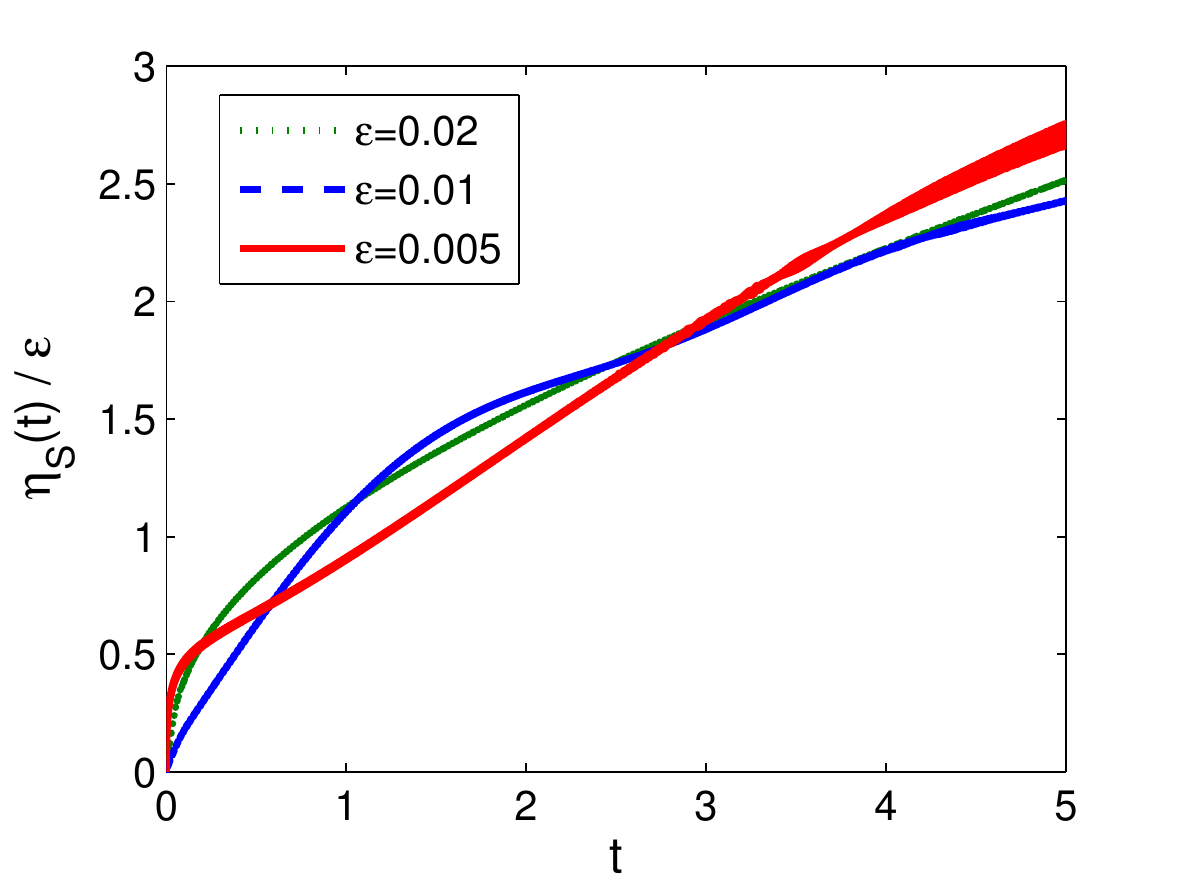,height=6cm,width=8cm}}
\caption{Time evolution of $\eta_{\rm sw}(t)$ and $\eta_{\rm s}(t)$ with nonsmooth data (\ref{nonsmooth}) for $m=1$ under different $\eps$.}\label{fig:convergence3}
\end{figure}

From Figs. \ref{fig:convergence}-\ref{fig:convergence3}, we can draw the following conclusions:

(i) The solution of the NKGE (\ref{KG}) converges to that of the
NLSW (\ref{nlsw}) quadratically in $\eps$
(and uniformly in time) provided that the initial data in (\ref{KG}) is smooth or at least satisfies
$\phi_1$ and $\phi_2\in H^2(\Omega)$, i.e.
$$\|u(\cdot,t)-u_{\rm sw}(\cdot,t)\|_{H^1}\leq C_0\eps^2,\qquad t\ge0,$$
where the constant $C_0>0$ is independent of $\eps$ and time $t\ge0$.

(ii) The solution of the NKGE (\ref{KG}) converges to that of the
NLSE  (\ref{nls}) quadratically in $\eps$
(in general, not uniformly in time) provided that the initial data in (\ref{KG}) is smooth or at least
satisfies $\phi_1$ and $\phi_2\in H^3(\Omega)$, i.e.
\[\|u(\cdot,t)-u_{\rm s}(\cdot,t)\|_{H^1}\leq (C_1+C_2T)\eps^2, \quad 0\le t\le T,
\]
where $C_1$ and $ C_2$ are two positive constants which are
independent of $\eps$ and $T$. On the contrary, if the regularity of the initial data is weaker, e.g.
$\phi_1$ and/or $\phi_2\in H^2(\Omega)$, then the convergence rate collapses to linear rate, i.e.
\begin{equation*}
\|u(\cdot,t)-u_{\rm s}(\cdot,t)\|_{H^1}\leq (C_3+C_4T)\eps, \ \ 0\le t\le T,
\end{equation*}
where  $C_3$ and $C_4$ are two positive constants which are independent of $\eps$ and $T$.
Rigorous mathematical justification for these numerical observations is on-going.

(iii) Under the same $\eps$ and at the same time $t$, the error $\eta_{\rm sw}(t)$ is much small than $\eta_{\rm w}(t)$. It indicates that the NLSW (\ref{nlsw}) would be a better choice to design the LI scheme than the limit model (\ref{nls}), especially considering the long time behavior of the approximation.

\subsection{Wave interactions in 2D}

We take $d=2$ and $\lambda=1$ in the NKGE (\ref{KG}) and choose
the initial data as
\begin{equation}\label{C4:2dnum}
\begin{split}
 &\phi_1(x,y)=\exp{(-(x+2)^2-y^2)}+\exp{(-(x-2)^2-y^2)},\\
 & \phi_2(x,y)=\exp{(-x^2-y^2)},\qquad (x,y)\in{\mathbb R}^2.
 \end{split}
 \end{equation}
The problem is solved numerically on a bounded computational domain
$\Omega=(-16,16)\times(-16,16)$ with the periodic boundary condition.
Fig. \ref{C4:fig2D} shows contour plots of the solutions of the NKGE (\ref{KG}) in 2D under different $\eps$.

\begin{figure}[h!]
\centerline{\psfig{figure=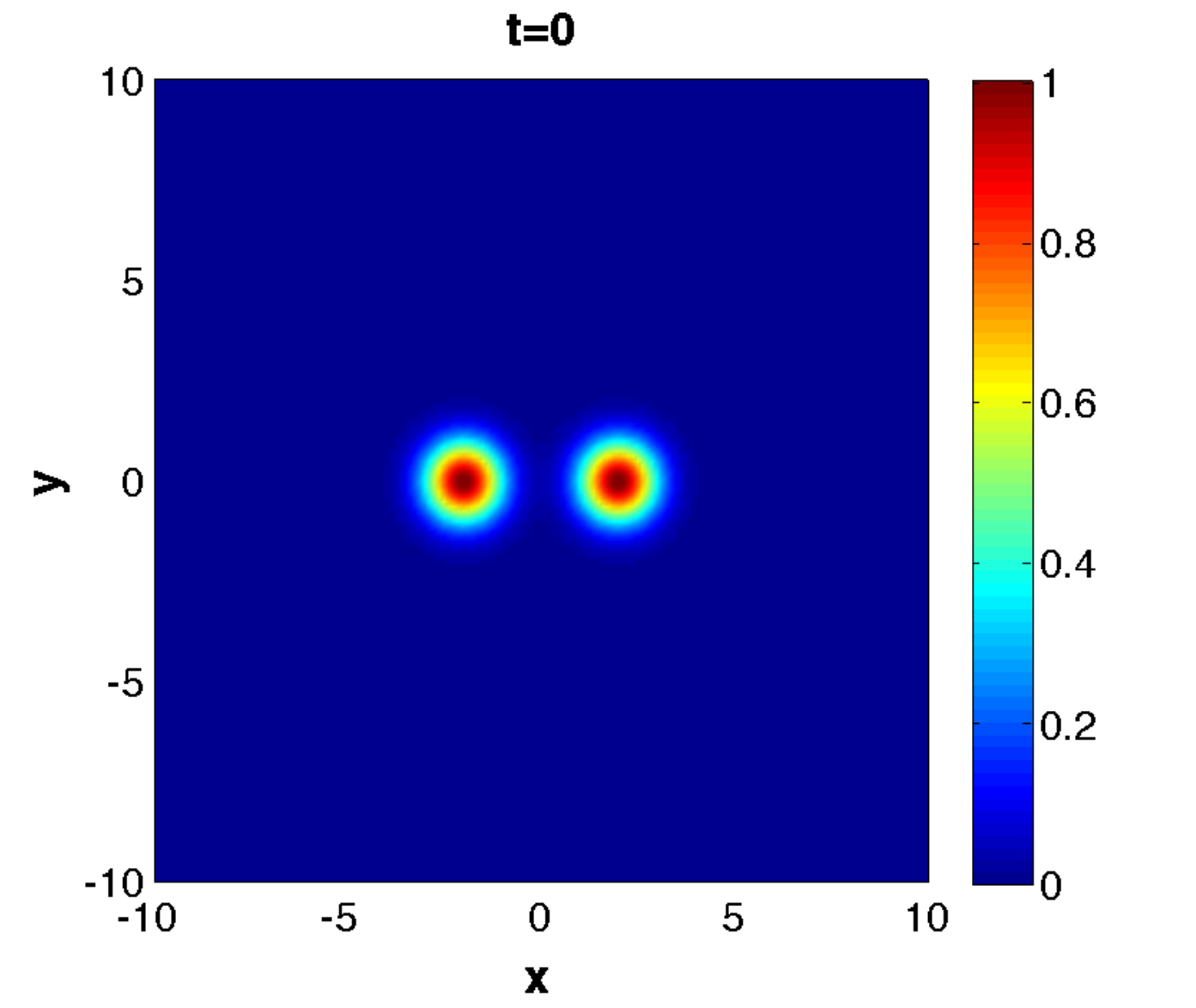,height=6cm,width=7cm}
\psfig{figure=t0-eps-converted-to.pdf,height=6cm,width=7cm}}
\centerline{\psfig{figure=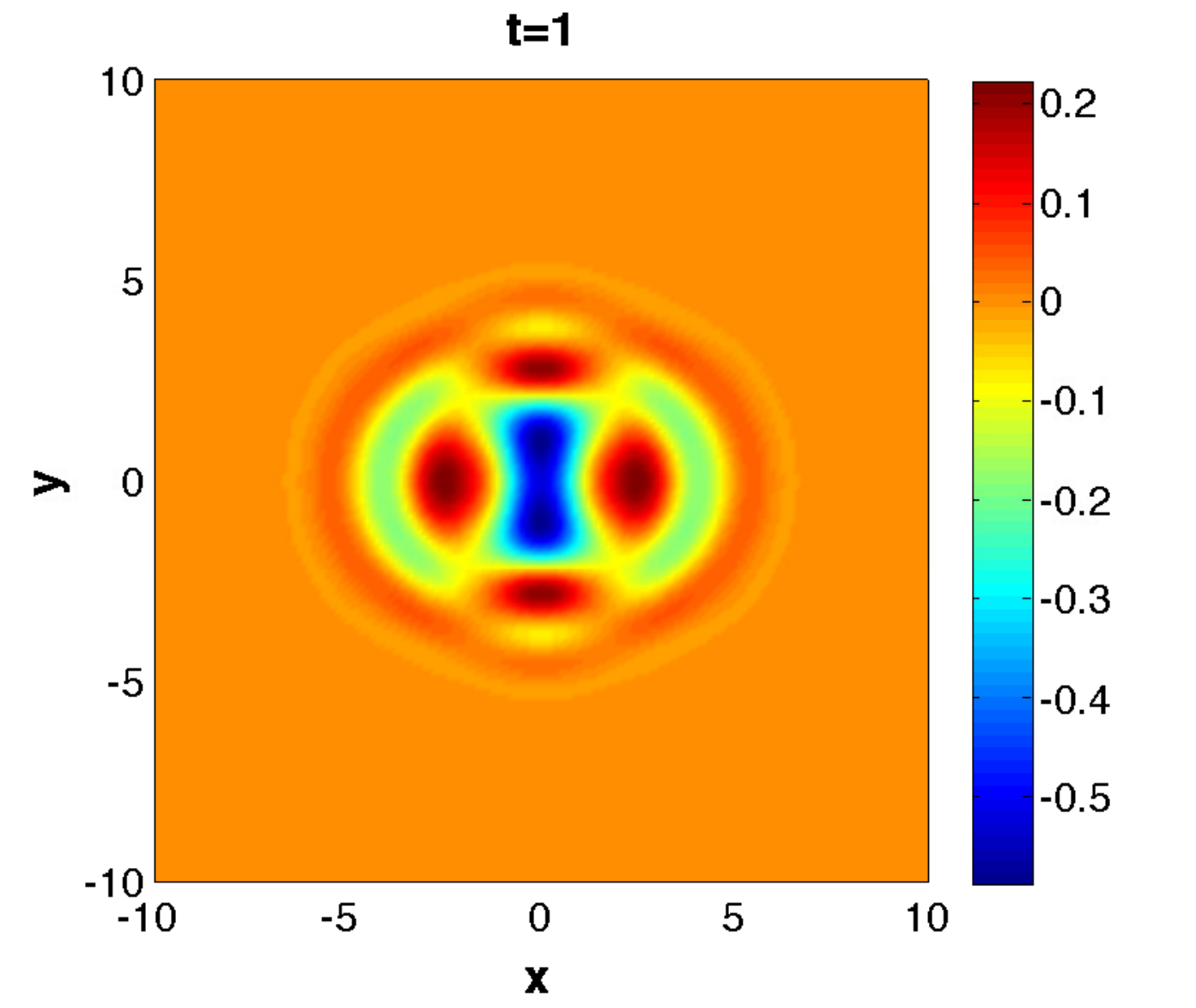,height=6cm,width=7cm}
\psfig{figure=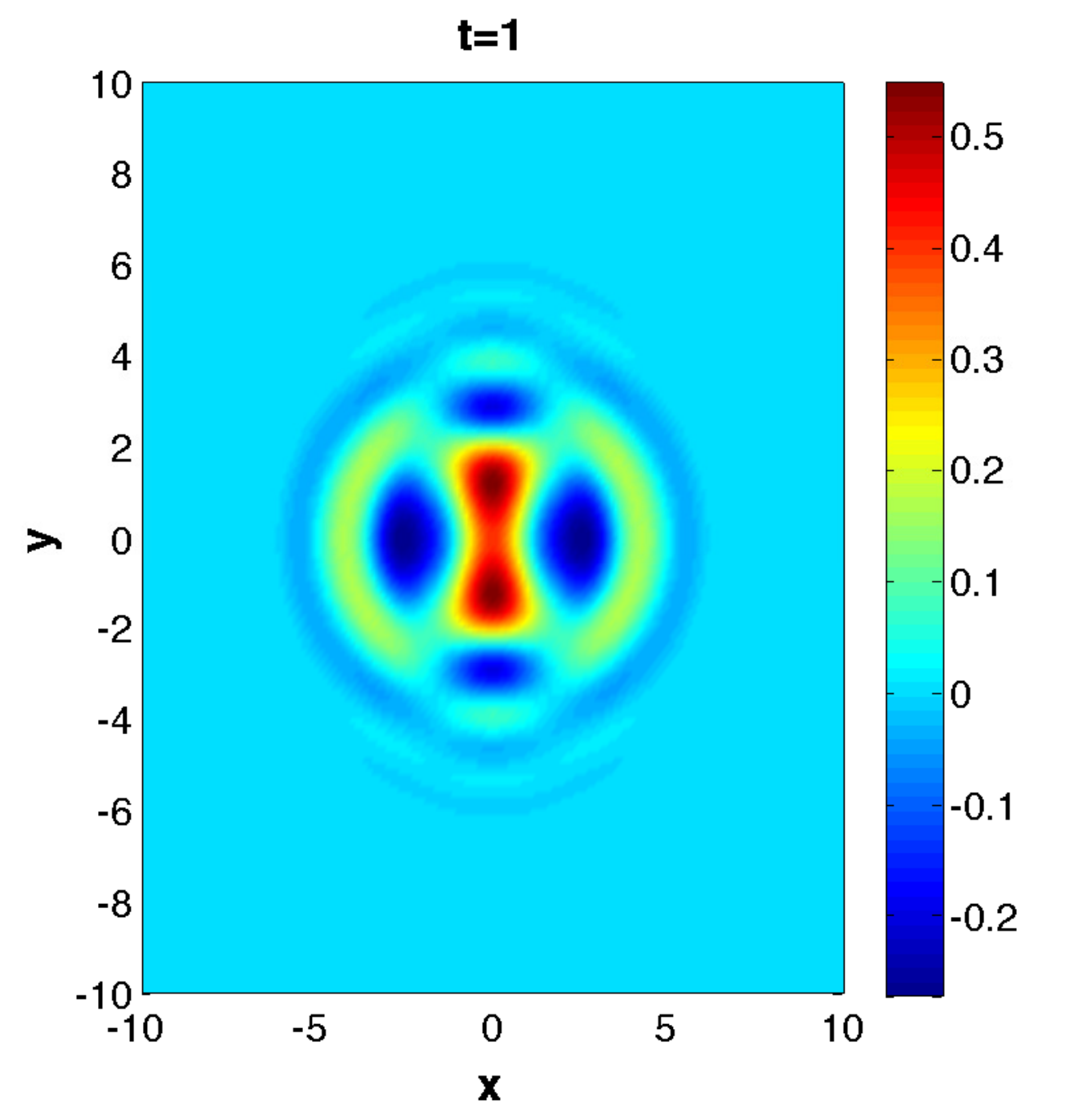,height=6cm,width=7cm}}
\centerline{\psfig{figure=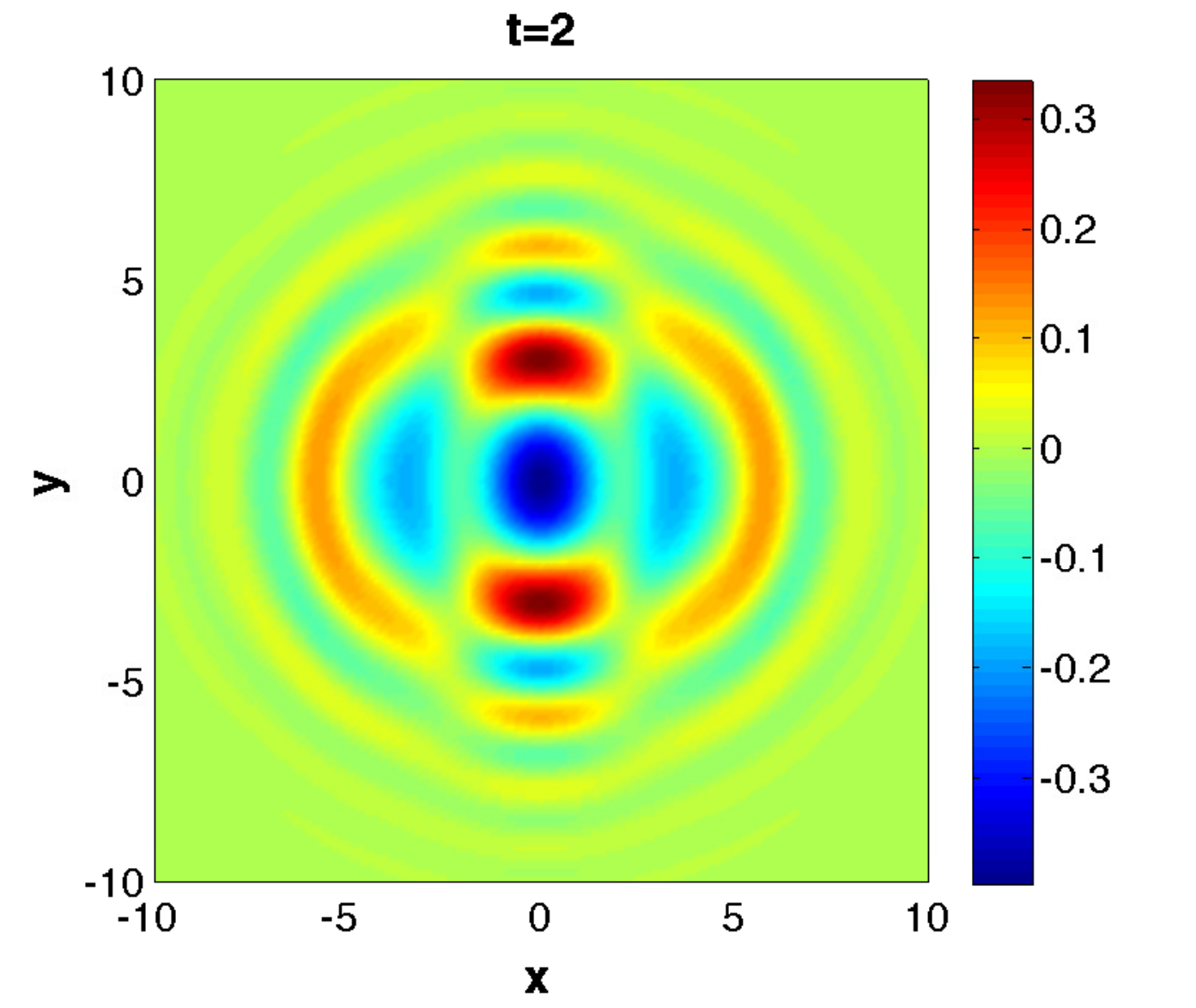,height=6cm,width=7cm}
\psfig{figure=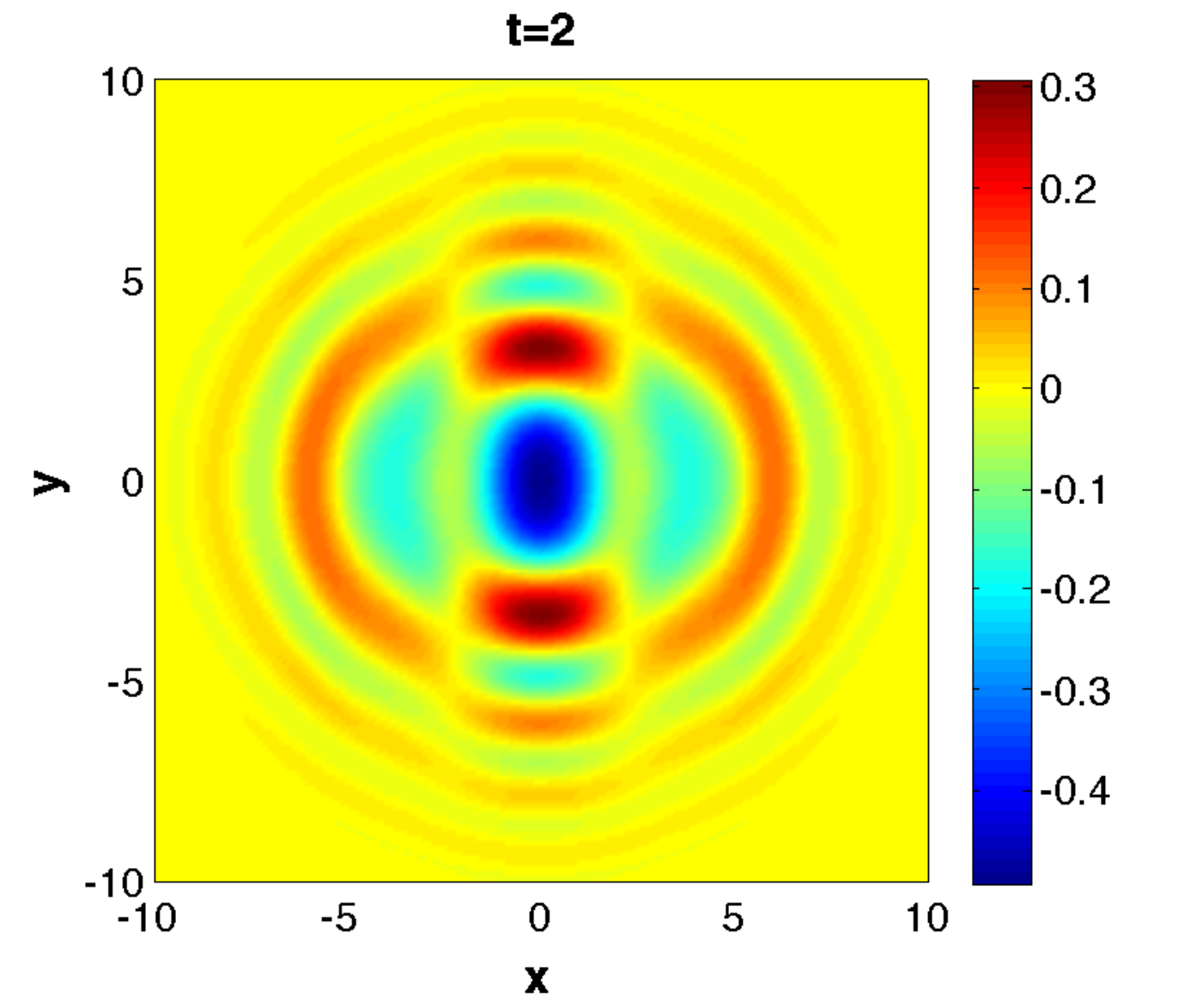,height=6cm,width=7cm}}
\caption{Contour plots of the solutions of the NKGE \eqref{KG} in 2D at different time $t$ under $\eps=0.05$ (first row) and $\eps=0.005$ (second row).}\label{C4:fig2D}
\end{figure}

\section{Conclusions}\label{sec:conc}
We systematically studied and compared different numerical methods to solve the nonlinear Klein-Gordon equation (NKGE) in the nonrelativistic limit regime, while the solution  is highly oscillatory in time in the limit regime. The numerical methods considered here include the classical finite difference time domain methods, the exponential wave integrator (EWI) spectral method, the time-splitting (TS) spectral method, the limit integrators, and the recently proposed uniformly accurate (UA) methods namely the multiscale time integrator (MTI) spectral method, the two-scale formulation (TSF) method and the iterative exponential integrator (IEI). We emphasized the finite time error bound of each method and the resolution capacity in terms of the oscillation wavelength in the limit regime. Systematical comparisons between the methods in the accuracy, computational complexity and other mathematical properties were carried out. Numerical experiments were done to show and compare the performance of each method from the classical regime to the limit regime. Our results show the EWI and TS methods are most efficient in the classical regime, while the UA methods are more powerful in the intermediate and limit regimes. Among the UA methods, the uniformly and optimally accurate methods are the most efficient and accurate for $\eps\in(0,1]$ . Finally, the UA numerical methods were applied to study numerically the convergence rates of the NKGE \eqref{KG} to its limiting models and to simulate wave interaction in two dimensions.


\section*{Acknowledgements}
This work was supported by the Ministry
of Education of Singapore grant
R-146-000-223-112 (W. Bao) and the French ANR project MOONRISE ANR-14-CE23-0007-01 (X. Zhao).

\bibliographystyle{model1-num-names}
\bibliography{<your-bib-database>}

\end{document}